\documentclass[11pt]{article}
\usepackage{amssymb,amsmath,amsthm,amsfonts,mathrsfs}
\usepackage{mathtools} 
\usepackage{bm}
\usepackage{enumerate}

\usepackage{graphicx}  
\usepackage{tikz} 
\usepackage{float}
\usepackage{wrapfig} 
\usepackage{subcaption}
\usepackage{epstopdf}

\usepackage{booktabs}
\usepackage{multirow}
\usepackage{multicol}
\usepackage{threeparttable}
\usepackage{algpseudocode}
\usepackage{array}
\usepackage[ruled,vlined,linesnumbered,commentsnumbered]{algorithm2e}

\usepackage{cite}
\usepackage[colorlinks, linktocpage, citecolor=black, linkcolor=black]{hyperref}

\usepackage{enumitem}

\usepackage{comment}  

\numberwithin{equation}{section}

\theoremstyle{remark}

\def\to{\rightarrow}

\SetKwInput{KwIn}{\bfseries Input}
\SetKwInput{KwOut}{\bfseries Output}
\SetKwFor{ForDo}{\bfseries for}{\bfseries do}{\bfseries end}

\SetCommentSty{mycommfont}

\allowdisplaybreaks 
\title{\bf A Unified  Phase-Field Fourier Neural Network Framework for Topology Optimization}
\author{Jing Li\thanks{School of Mathematical Sciences, East China Normal University, Shanghai, China. (\texttt{betterljing@163.com})} 
        \and Xindi Hu\thanks{Department of Applied Mathematics, The Hong Kong Polytechnic University, Kowloon, Hong Kong, China (\texttt{xindi.hu@polyu.edu.hk})} \and Helin Gong\thanks{Paris Elite Institute of Technology, Shanghai Jiao Tong University, Shanghai, China. (\texttt{gonghelin@sjtu.edu.cn})} \and Wei Gong\thanks{NCMIS\&SKLMS,Institute of Computational Mathematics and Scientific/Engineering Computing, Academy of Mathematics and Systems Science, Chinese Academy of Sciences, Beijing, China. (\texttt{wgong@lsec.cc.ac.cn})} \and Shengfeng Zhu\thanks{Corresponding author. Key Laboratory of Mathematics and Engineering Applications \& Shanghai Key Laboratory of Pure Mathematics and Mathematical Practice, East China Normal University, Shanghai, China. (\texttt{sfzhu@math.ecnu.edu.cn})}}
\date{\today}

\begin{document}

\maketitle

\begin{abstract}
    We propose Alternating Phase-Field Fourier Neural Networks (APF-FNNs) as a unified and physics-based framework for topology optimization. The approach decouples the design problem by representing the state, adjoint, and topology fields with three separate Fourier neural networks, which are trained via a stable collaborative alternating scheme applicable to both self-adjoint and non-self-adjoint problems. To obtain well-resolved designs, the Ginzburg--Landau energy functional is embedded in the loss of the topology network as an intrinsic regularizer, naturally enforcing smooth and distinct interfaces between the two phases. Phase-field updates are driven by adjoint-based optimality conditions, and design sensitivities are evaluated efficiently using automatic differentiation, ensuring that the gradients correspond to exact total derivatives rather than naive partial derivatives. In contrast to classical phase-field methods, APF-FNNs exploit these physically consistent design gradients directly, avoiding pseudo-time gradient-flow solvers. By formulating physics-driven losses from variational principles or strong-form PDE residuals, the framework is broadly applicable to 2D and 3D benchmark problems, including compliance minimization, eigenvalue maximization, and Stokes/Navier--Stokes flow optimization. Across these examples, APF-FNNs consistently yield competitive performance and well-resolved topologies, establishing a versatile and scalable foundation for physics-driven computational design.
\par\vspace{1ex}

\vspace{0.25cm}
\end{abstract}

\begin{figure}[h!] 
    \centering
    \includegraphics[width=1.01\linewidth]{ 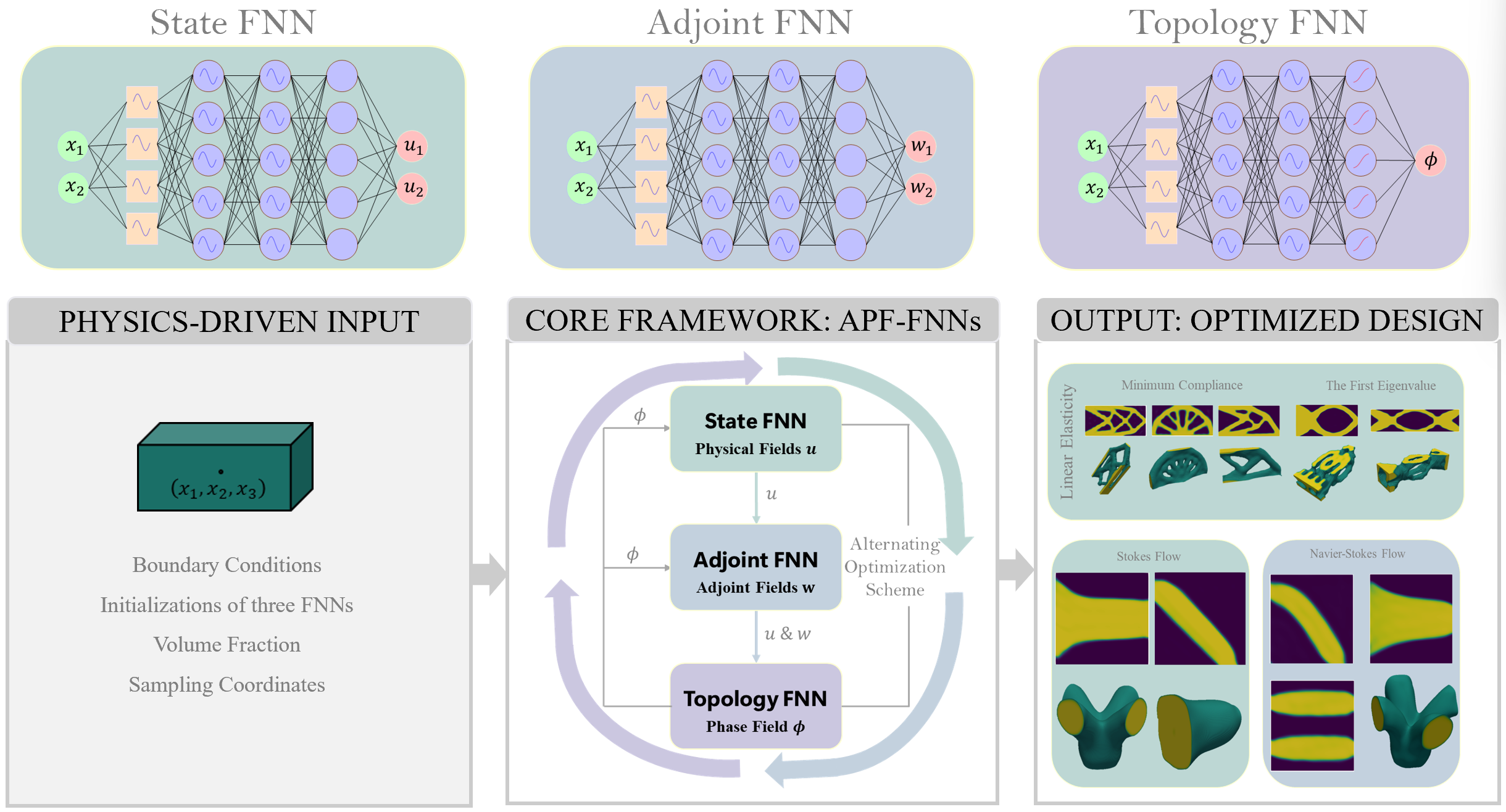} 
    \caption{\textbf{Schematic overview of the Alternating Phase-Field Fourier Neural Networks (APF-FNNs).} The process begins with the physics settings (left panel). The core framework (middle panel) features an alternating architecture with three distinct FNNs parameterizing the state \(\bm{u}\), adjoint \(\bm{w}\) and topology \(\phi\) fields. These networks are trained collaboratively through a stable alternating optimization scheme. The versatility of the framework (right panel) is demonstrated by solving a range of benchmark problems, encompassing both self-adjoint (e.g., linear elasticity and Stokes flow) and non-self-adjoint (e.g., Navier--Stokes flow) systems.} 
    \label{fig: abstract}
\end{figure}

\section{Introduction}

Topology optimization (TO) has emerged as a powerful computational paradigm in engineering design, enabling the computation of optimal material distributions within a prescribed domain under specified loads and boundary conditions \cite{bendsoe2003topology, zargham2016topology, ALLAIRE20211}. By systematically exploring a vast design space, TO produces high-performance topologies that often surpass conventional and experience-based designs. This capability has driven remarkable progress in developing lightweight and efficient components across diverse engineering fields \cite{giraldo2020multi, CHEN2024113119}. Its impact has been further amplified by additive manufacturing \cite{takezawa2020sensitivity, luo2020additive, kangzhan31122025}, which enables the realization of intricate and performance-driven topologies, establishing TO as a cornerstone of modern computational design.

Originating from the seminal homogenization method \cite{bendsoe1988generating}, practical TO \cite{sigmund2013topology} has been dominated by mesh-based density formulations such as Solid Isotropic Material with Penalization (SIMP) \cite{bendsoe1989optimal, bendsoe1999material} and Bi-directional Evolutionary Structural Optimization \cite{huang2007convergent}. Despite their success, these approaches suffer from well-known numerical pathologies, including checkerboard patterns, mesh-dependence and persistent intermediate densities \cite{diaz1995checkerboard, sigmund1998numerical}. These issues reflect a fundamental ill-posedness of the continuous problem and are typically mitigated by ad-hoc remedies such as density filtering or projection to impose a minimum length scale \cite{bourdin2001filters, bruns2001topology}. To overcome these limitations, a number of approaches decouple the geometric representation from the analysis mesh. Level-set approaches can produce crisp, well-defined boundaries but struggle to nucleate new holes inside the domain \cite{sethian2000structural, wang2003level, ALLAIRE20211}. Phase-field methods have emerged as a particularly compelling alternative, effectively circumventing the limitations of both density and level-set formulations \cite{hu2023accelerating, qian2022phase, wu2018multi, li2025efficient}. Grounded in phase-transition theory, the phase-field framework introduces a diffuse interface governed by a Ginzburg--Landau energy \cite{wang2004phase,Bourdin_Chambolle_2003,plotnikov2023geometric}. This representation naturally accommodates complex topological changes such as hole nucleation and member merging, while the Ginzburg--Landau functional provides a physics-based regularization that enforces a minimum length scale and suppresses numerical artifacts without heuristic filters \cite{chen2002phase, hu2023accelerating, qian2022phase, wu2018multi,
li2025efficient}.

Despite these conceptual advantages, conventional phase-field TO is hindered by significant computational costs. At each design iteration one must solve both the governing state/adjoint equations and a phase-field evolution equation, typically of Allen--Cahn or Cahn--Hilliard type \cite{li2024energy, bourdin2006phase}. Within a Finite Element Method (FEM), this requires repeated solutions of large and often nonlinear systems at every optimization step \cite{gogu2015improving}. A second bottleneck stems from the coupling between geometry and discretization: the sharpness of the material interface is tied to the mesh resolution, since the diffuse interface must be resolved by several elements for numerical stability \cite{Bourdin_Chambolle_2003}. High-resolution designs with crisp boundaries
therefore demand globally refined meshes. The combination of expensive inner gradient-flow iterations and fine discretizations renders high-resolution 2D and most 3D problems computationally prohibitive \cite{amir2011reducing}, and calls for more efficient and mesh-free solution strategies.

Deep neural networks have recently emerged as powerful mesh-free function approximators for PDEs, driving a paradigm shift in scientific computing \cite{lagaris1998artificial, raissi2019physics, ew2018deep, karniadakis2021physics, li2020fourier, yamada2024smo, YU2026117002}. Among these, Physics-Informed Neural Networks (PINNs) \cite{raissi2019physics} minimize strong-form residuals in the loss and can be applied to systems without a known variational structure, while the Deep Ritz Method (DRM) \cite{ew2018deep} is tailored to problems that admit an energy functional and often exhibits superior stability and accuracy in such cases. This distinction naturally suggests a hybrid strategy: employ DRM when a variational formulation is available and resort to PINN-style residual losses otherwise. In particular, the variational structure of the phase-field model makes an energy-based treatment especially attractive for TO.

Early studies of neural networks in TO are predominantly data-driven, casting the problem as image-to-image translation or generative design \cite{shin2023topology, woldseth2022use, sosnovik2019neural, banga20183d, Yamasaki2021, guoxu2026}. The first wave of research focused on surrogate models that learn a direct mapping from problem settings to final topologies, aiming to bypass costly iterative solvers \cite{banga20183d}. More recent works employ generative models, such as
GANs, to explore the design space from a latent representation and to encode complex, non-differentiable manufacturing constraints \cite{nie2021topologygan, parrott2023multidisciplinary, greminger2020generative}. Although promising, these purely data-driven paradigms require large datasets of pre-computed designs and tend to reproduce visual patterns rather than enforce physical principles. As a result, they may produce mechanically unsound or disconnected layouts, and their performance is ultimately limited by the coverage of the training data, which undermines their generality and reliability for engineering applications.

Physics-driven strategies aim to address these shortcomings by embedding physical models directly into the optimization. Initial approaches re-parameterize the design field with neural networks \cite{chandrasekhar2021tounn, sanu2025neural} but still rely on finite element solvers for the state analysis, and thus inherit their computational burden. Subsequent developments replace FEM by energy-based or PINN-style surrogates that solve for the state fields without labeled data \cite{jeong2023physics, zehnder2021ntopo}, yet the design update often remains coupled to conventional optimization schemes. More recently, mesh-free dual network architectures have been proposed, in which a density network describes the material layout and a physics network approximates the corresponding state fields \cite{singh2024dual, joglekar2024dmf}. A unified energy-based loss derived from variational principles couples the two networks and enables simultaneous, end-to-end optimization of topology and physical response while eliminating explicit finite element analysis and hand-crafted sensitivity derivations. However, existing formulations either combine all physical and design objectives into a single monolithic loss \cite{yousefpour2025simultaneous}, leading to a complex optimization landscape with delicate penalty balancing and potential training instabilities, or are restricted to specific objectives, such as linear elastic compliance in multi-material phase-field design \cite{lai2025dual}. Beyond penalty balancing and training stability, a more fundamental issue is that gradients from automatic differentiation can be inconsistent with the true total design derivatives, particularly for non-self-adjoint systems. These limitations highlight the need for a general, physics-driven framework that can handle both self-adjoint and non-self-adjoint systems while maintaining a stable and tractable optimization procedure.

In this work, we propose \textbf{Alternating Phase-Field Fourier Neural Networks (APF-FNNs)}, a unified and physics-driven framework for phase-field topology optimization. As illustrated in Fig.~\ref{fig: abstract}, the main features of APF-FNNs can be summarized as follows:
\begin{itemize}
    \item \textbf{Decoupled multi-FNN architecture.}
          In contrast to coupled formulations, APF-FNNs employ three distinct Fourier Neural Networks (FNNs) to parameterize the state, adjoint and topology fields. The use of FNNs mitigates the spectral bias of standard multilayer perceptrons and enables compact representation of sharp interfaces and complex fields without prohibitive network sizes.

    \item \textbf{Robust alternating optimization for broad applicability.}
          The three networks are trained in a cyclic alternating fashion: in each outer iteration, two networks are held fixed while the third is updated by minimizing an appropriate physics-driven loss. This decoupling separates physics solving from design update, improves stability, and enables a unified treatment of both self-adjoint and challenging non-self-adjoint systems.

    \item \textbf{Intrinsic physics-based regularization.} 
          The Ginzburg--Landau energy functional is incorporated into the topology network's loss. This acts as an intrinsic regularizer, promoting well-defined designs with smooth interfaces and eliminating the need for ad-hoc density filtering schemes.

    \item \textbf{Direct solution of optimality conditions with consistent sensitivities.}
          Instead of evolving a pseudo-time gradient flow for the phase-field, APF-FNNs directly minimize the total objective functional. Design sensitivities are obtained from adjoint-based optimality conditions and evaluated efficiently by automatic differentiation, so that the gradients used to update the topology network coincide with the correct total derivatives of the objective rather than naive partial derivatives.
\end{itemize}

This paper is organized as follows. Section~\ref{sec:pf-to} establishes the mathematical formulation for phase-field topology optimization and contrasts the conventional $L^2$ gradient-flow solution strategy with the proposed neural network-based alternating scheme. Section~\ref{sec:method} presents the APF-FNN methodology, including the FNN architecture, the construction of physics-driven loss functions for diverse multiphysics problems, and the associated gradient computation strategy. Section~\ref{sec:experiments} demonstrates the performance and versatility of APF-FNNs on a range of 2D and 3D numerical benchmarks. Concluding remarks and future research directions are provided in Section~\ref{sec:conclusion}.

\section{Phase-field TO formulations}
\label{sec:pf-to}

In the phase-field approach to topology optimization \cite{garcke2016shape, blank2011phase, jin2024adaptive}, the classical discrete material distribution problem is reformulated as a continuous one governed by a phase-field variable $\phi:\Omega \subset \mathbb{R}^d \to [0,1]$ ($d=2,3$), where $\Omega$ is an open bounded domain with Lipschitz continuous boundary $\partial\Omega$. The variable smoothly interpolates between the material phase ($\phi=1$) and the void phase ($\phi=0$) across a diffuse interface ($0 < \phi < 1$). To resolve the inherent ill-posedness of TO and control the geometric complexity of the design, the formulation augments the primary objective $\mathcal{J}_{\text{obj}}$ with a Ginzburg--Landau regularization $\mathcal{E}_{GL}(\phi)$, which penalizes interfacial area and promotes well-defined boundaries \cite{blank2011phase, takezawa2010shape}.

A general phase-field topology optimization problem can be written as the constrained minimization
\begin{equation}
\begin{aligned}
\min_{\phi} \quad & \mathcal{J}(\phi) = \mathcal{J}_{\text{obj}}(\phi, \bm{u}) + \mathcal{E}_{GL}(\phi) + \mathcal{P}_{\text{vol}}(\phi) \\
\text{s.t.} \quad & \mathcal{R}(\phi, \bm{u}) = 0 \quad \text{in } \Omega, \\
& \mathcal{B}(\phi, \bm{u}) = 0 \quad \text{on } \partial\Omega,
\end{aligned}
\label{eq:general_to}
\end{equation}
where the total objective $\mathcal{J}$ consists of the physical objective $\mathcal{J}_{\text{obj}}$, the Ginzburg--Landau regularization $\mathcal{E}_{GL}(\phi)$, and a volume constraint functional $\mathcal{P}_{\text{vol}}(\phi)$. The state variable $\bm{u}$ is constrained by the governing equations $\mathcal{R}(\phi,\bm{u})=0$ in $\Omega$ and boundary conditions $\mathcal{B}(\phi,\bm{u})=0$ on $\partial\Omega$. The term $\mathcal{P}_{\text{vol}}(\phi)$ enforces the volume constraint $\int_{\Omega} \phi \,\mathrm{d}x = \beta |\Omega|$ with target volume fraction $\beta$ and can be treated, for example, by quadratic penalties
\cite{Bourdin_Chambolle_2003}, Lagrange multipliers or augmented Lagrangian techniques \cite{li2024energy, jin2024adaptive, li2025adaptive}.

The Ginzburg--Landau functional is defined as
\begin{equation}
\mathcal{E}_{GL}(\phi) = \gamma\int_{\Omega} \left( \frac{\epsilon}{2} |\nabla\phi|^2 + \frac{1}{\epsilon} W(\phi) \right) \mathrm{d}x.
\label{eq:gl_energy}
\end{equation}
where the parameter $\gamma>0$ scales the overall regularization and $\epsilon>0$ controls the interfacial thickness. The function $W(\phi)=\tfrac{1}{4}\phi^2(1-\phi)^2$ is a double-well potential with minima at $\phi=0$ and $\phi=1$, thereby promoting an almost binary material distribution \cite{garcke2016shape, jin2024adaptive, li2024energy, li2025adaptive}.

\subsection{Classical formulation with $L^2$ gradient flow}

A standard strategy for solving PDE-constrained optimization problem in Eq.~(\ref{eq:general_to}) is to introduce a gradient-flow dynamics for the phase-field. The design variable $\phi$ is evolved in a pseudo-time $t>0$ along the steepest descent of the total functional $\mathcal{J}(\phi)$ \cite{li2024energy, li2025adaptive, takezawa2010shape}. The precise evolution equation depends on the choice of inner product in the function space of~$\phi$.

To avoid the computational cost of the fourth-order Cahn--Hilliard equation arising from an $H^{-1}$ gradient flow, many works adopt an $L^2$ gradient flow of Allen--Cahn type \cite{takezawa2020sensitivity, li2024energy}. In weak form, the evolution of $\phi$ is given by
\begin{equation}
    (\partial_t \phi, \psi)_{L^2}
    =
    - \delta_\phi \mathcal{J}(\phi)[\psi]
    \quad \forall\,\psi \in H^1(\Omega),
    \label{eq: H1gradientflow}
\end{equation}
where $(\cdot,\cdot)_{L^2}$ denotes the $L^2$ inner product. The variational derivative decomposes as
\begin{equation}
    \delta_\phi \mathcal{J}(\phi)[\psi]
    =
    \delta_\phi \mathcal{J}_{\text{obj}}(\phi, \bm{u})[\psi]
    + \delta_\phi \mathcal{E}_{GL}(\phi)[\psi]
    + \delta_\phi \mathcal{P}_{\text{vol}}(\phi)[\psi].
\end{equation}
Here $\delta_\phi \mathcal{J}_{\text{obj}}$ represents the sensitivity of the physical objective to perturbations of the phase field, while the remaining two terms are the contributions of the Ginzburg--Landau regularization and the volume constraint. The objective sensitivity $\delta_\phi \mathcal{J}_{\text{obj}}$ is typically evaluated by an adjoint method, which avoids differentiating the state variable $\bm{u}$ with respect to $\phi$ explicitly \cite{takezawa2020sensitivity}.

\begin{figure}[h!] 
    \centering
    \includegraphics[width=1\linewidth]{ 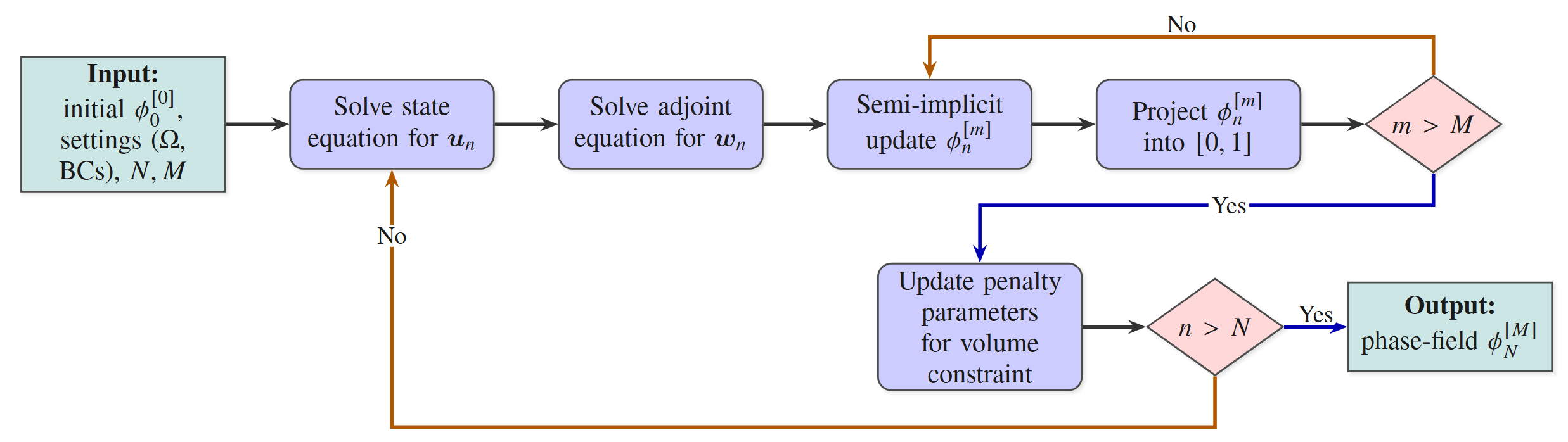} 
    \caption{Flowchart of the nested iterative scheme for conventional phase-field topology optimization.} 
    \label{fig: conventional_optimization}
\end{figure}

The resulting solution procedure follows a nested iterative scheme, illustrated in Fig.~\ref{fig: conventional_optimization}. An outer loop, indexed by $n$, updates the design, while an inner loop, indexed by $m$, performs gradient-flow steps for a fixed set of parameters. At the beginning of each outer iteration, the state and adjoint equations are solved to obtain the response $\bm{u}_n$
and the corresponding sensitivity $\delta_\phi \mathcal{J}_{\text{obj}}$. This sensitivity drives $M$ inner steps, in which the design $\phi_n$ is advanced by a semi-implicit time discretization of the $L^2$ gradient flow \eqref{eq: H1gradientflow}, followed by a projection enforcing the box constraint $\phi_n\in[0,1]$. After the inner loop, the penalty parameters associated with the volume constraint are updated. The process is repeated for $N$ outer iterations or until convergence to a locally optimal design. While effective, such gradient-flow-based formulations rely on an artificial time evolution and nested iterative solvers, which complicate the direct optimization of the objective functional and may limit scalability.

\subsection{Alternating neural parameterization with APF-FNNs}

\begin{figure}[h!] 
    \centering
    \includegraphics[width=1\linewidth]{ 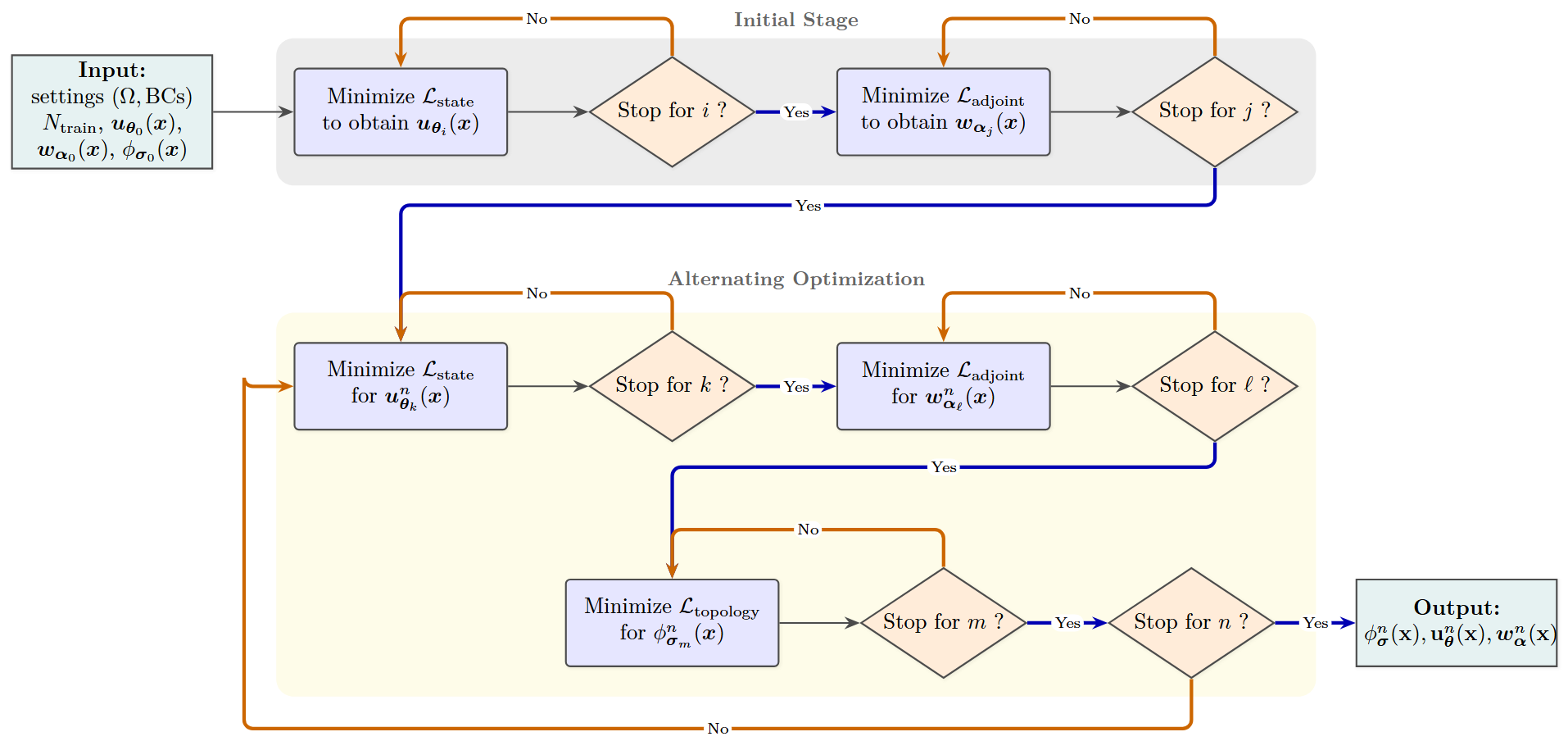} 
    \caption{Flowchart of the proposed APF-FNNs. The process begins with an Initial Stage to pre-train the state and adjoint networks. This is followed by the main Alternating Optimization stage, where the neural networks for the state, adjoint and topology are refined by sequentially minimizing their respective loss functions.} 
    \label{fig: mesh_free_optimization}
\end{figure}

The proposed APF-FNNs replace the classical mesh-based gradient-flow solver by a physics-driven, neural parameterization of the phase-field topology optimization problem. Instead of iteratively solving discretized PDEs on a mesh, we approximate the relevant fields with neural networks and directly minimize their associated loss functionals. As sketched in Fig.~\ref{fig: mesh_free_optimization}, the framework employs three distinct networks to represent the state $\bm{u}_{\boldsymbol{\theta}}(\bm{x})$, the adjoint $\bm{w}_{\boldsymbol{\alpha}}(\bm{x})$ and the phase-field $\phi_{\boldsymbol{\sigma}}(\bm{x})$.

The optimization begins with an \texttt{Initial Stage} that provides a robust warm start for the alternating scheme. In this stage, the state network $\bm{u}_{\boldsymbol{\theta}}$ and the adjoint network $\bm{w}_{\boldsymbol{\alpha}}$ are sequentially pre-trained for a fixed initial phase-field $\phi_{\boldsymbol{\sigma}}$ by minimizing their respective loss functions $\mathcal{L}_{\text{state}}$ and $\mathcal{L}_{\text{adjoint}}$ for a prescribed number of iterations. The construction of these losses is physics-driven: when the governing equations admit an energy functional, we employ a Deep Ritz formulation \cite{ew2018deep}; otherwise, a PINN-style loss based on strong-form residuals is used \cite{raissi2019physics}.

After initialization, the method proceeds to the main \texttt{Alternating Optimization} stage, where the three networks are updated in a cyclic scheme. Each outer iteration $n$ consists of three steps: (i) the state network $\bm{u}_{\boldsymbol{\theta}}$ is trained for $k$ steps by minimizing $\mathcal{L}_{\text{state}}$; (ii) given the updated state, the adjoint network $\bm{w}_{\boldsymbol{\alpha}}$ is trained for $l$ steps by minimizing $\mathcal{L}_{\text{adjoint}}$; and (iii) with the current state and adjoint fields fixed, the topology network
$\phi_{\boldsymbol{\sigma}}$ is refined for $m$ steps. In contrast to conventional approaches based on time-stepping the $L^2$ gradient flow, the topology update is formulated as a direct minimization of the total objective
\begin{equation}
  \mathcal{L}_{\text{topology}}(\boldsymbol{\sigma};
  \bm{u}_{\boldsymbol{\theta}}, \bm{w}_{\boldsymbol{\alpha}})
  :=
  \mathcal{J}(\phi_{\boldsymbol{\sigma}},
              \bm{u}_{\boldsymbol{\theta}},
              \bm{w}_{\boldsymbol{\alpha}}).
  \label{eq:nn_design_loss}
  \nonumber
\end{equation}

This neural formulation offers two key advantages:
\begin{itemize}
    \item \textbf{Mesh-free representation.}
          The networks operate directly on continuous spatial coordinates $\bm{x}$, making the framework inherently mesh-free. This eliminates the cost of mesh generation and refinement, which is especially burdensome for complex geometries and three-dimensional problems.
    \item \textbf{Efficient optimization strategy.}
          The computational efficiency is enhanced in two ways:
          \begin{itemize}
              \item \textit{Direct gradient computation.}
                    By combining adjoint-based sensitivity analysis with automatic differentiation, the gradients of the objective
                    with respect to the topology parameters $\boldsymbol{\sigma}$ are obtained efficiently via backpropagation. This avoids explicit time integration of a gradient-flow equation.
              \item \textit{Warm starts through transfer learning.}
                    The pre-trained networks from the \texttt{Initial Stage} provide a high-quality initialization for the alternating loop, and within the loop the parameters from iteration $n$ serve as warm starts for iteration $n+1$. This reuse of information allows the state and adjoint problems to be solved to moderate accuracy at intermediate steps, reducing the overall cost per outer iteration.
          \end{itemize}
\end{itemize}

%======================================================================
% SECTION 3: PROPOSED METHODOLOGY
%======================================================================
\section{APF-FNN methodology}
\label{sec:method}

As illustrated in Fig.~\ref{fig: abstract}, the proposed APF-FNN framework is built on three FNNs that parameterize the state, adjoint and topology fields. Denoting these networks by $\bm{u_\theta}(\bm{x})$, $\bm{w_\alpha}(\bm{x})$, and $\phi_{\bm{\sigma}}(\bm{x})$, the optimization proceeds in an alternating scheme: in each outer iteration, two of the networks are kept fixed while the remaining one is updated by minimizing a physics-driven loss function. A central aim of APF-FNNs is to obtain physically consistent design gradients for topology updates across both self-adjoint and non-self-adjoint systems. This decoupled design allows the state and adjoint networks to focus on accurately solving the state and adjoint PDEs, while the topology network is dedicated to exploring the design space under phase-field regularization and volume constraints. The remainder of this section details the components of the methodology: Section~\ref{subsec:FNN-architecture} introduces the FNN architecture used for all three fields, Section~\ref{subsec:physics-loss} describes the physics-based loss functions for the elasticity and flow problems, and Section~\ref{subsec:gradients} discusses gradient computation and adjoint integration.

\subsection{Neural network archtectures}
\label{subsec:FNN-architecture}

Multi-Layer Perceptrons (MLPs) are known to suffer from spectral bias, a tendency to learn low-frequency components of a target function much more readily than high-frequency details. This limitation is particularly detrimental for phase-field topology optimization, where sharp interfaces and complex geometries must be resolved accurately \cite{karniadakis2021physics, tancik2020fourier}. To mitigate this issue, all neural modules in the APF-FNN framework are implemented as FNNs, which augment a standard MLP with a high-frequency Fourier feature mapping of the input coordinates.

\subsubsection{Fourier neural networks} \label{subsub:Fourier neural networks}
Given a spatial coordinate $\bm{x} \in \mathbb{R}^{d}$, the FNN first lifts it to a higher dimensional feature space through the mapping
\begin{equation}
  \Gamma(\bm{x})
  =
  \big[\cos(2\pi B\bm{x}),\,\sin(2\pi B\bm{x})\big]^{\mathrm{T}},
  \label{eq:fourier-mapping}
\end{equation}
where $B \in \mathbb{R}^{m\times d}$ is a frequency matrix. At initialization, $B$ can be sampled from a Gaussian distribution or constructed from a structured frequency grid. This hyperparameter controls the range of frequencies present in the sinusoidal mapping and hence the level of detail that can be efficiently represented.

The resulting Fourier feature vector $\Gamma(\bm{x})\in\mathbb{R}^{2m}$ is then processed by a fully connected MLP to produce the final output $f(\bm{x};\bm{\xi})\in\mathbb{R}^{d_{\text{out}}}$. We denote the input to the MLP by $h^{(0)} = \Gamma(\bm{x})$ and define the hidden-layer updates as
\begin{equation}
  h^{(\ell)} = S\big(W^{(\ell)} h^{(\ell-1)} + b^{(\ell)}\big),
  \qquad \ell = 1,2,\dots,L,
  \label{eq:FNN-hidden}
\end{equation}
where $W^{(\ell)}\in\mathbb{R}^{n_\ell\times n_{\ell-1}}$ and $b^{(\ell)}\in\mathbb{R}^{n_\ell}$ denote the weight matrix and bias vector of the $\ell$-th hidden layer, respectively, and $S(\cdot)$ is a nonlinear activation function. The network prediction is obtained by a final linear layer
\begin{equation}
  f(\bm{x};\bm{\xi}) = W^{(L+1)} h^{(L)} + b^{(L+1)},
  \label{eq:FNN-output}
\end{equation}
with $W^{(L+1)}\in\mathbb{R}^{d_{\text{out}}\times n_L}$ and $b^{(L+1)}\in\mathbb{R}^{d_{\text{out}}}$. Collectively, the trainable parameters of a single FNN consist of the Fourier mapping matrix $B$ and the MLP parameters $\bm{\xi} = \{W^{(\ell)}, b^{(\ell)}\}_{\ell=1}^{L+1}$.

Within the APF-FNN framework, the same FNN architecture is instantiated three times to parameterize the state, adjoint, and topology fields, as schematically illustrated in Fig.~\ref{fig: abstract}, with each network having an independent set of parameters:

\begin{itemize}
  \item \textbf{State network.} The state network $\bm{u_\theta}(\bm{x})$, with parameters $\bm{\theta}=(B_{\bm{u}},\bm{\xi}_{\bm{u}})$, outputs the physical fields of interest. A multi-output formulation is adopted when multiple state variables are present (e.g., velocity and pressure in incompressible flow), allowing them to be represented jointly by a single network with shared parameters.
  \item \textbf{Adjoint network.} The adjoint network $\bm{w_\alpha}(\bm{x})$ follows the same architectural template but uses its own parameters $\bm{\alpha}=(B_{\bm{w}},\bm{\xi}_{\bm{w}})$. All adjoint variables are represented as components of a single multi-output network and are employed to compute consistent design sensitivities in non-self-adjoint problems.
  \item \textbf{Topology network.} The topology network $\phi_{\bm{\sigma}}(\bm{x})$, with parameters $\bm{\sigma}=(B_\phi,\bm{\xi}_\phi)$, encodes the phase-field design variable. It shares the same Fourier feature mapping as the other networks, while a sigmoid activation is applied at the output to enforce the box constraint $\phi_{\bm{\sigma}}(\bm{x})\in[0,1]$.
\end{itemize}

By providing a common representational backbone for state, adjoint, and topology fields, FNNs enable a unified yet modular architecture capable of capturing both smooth physical solutions and fine-scale geometric features \cite{tancik2020fourier, joglekar2024dmf}.

\subsubsection{FNNs for state, adjoint and topology}
\label{subsubsec: fnnsforall}

All benchmarks use Fourier-feature representations for the state, adjoint, and topology fields. While the framework is shared across problems, the feature construction and hyperparameters (e.g.\ bandwidths, numbers of modes, and layer widths) are problem dependent and specified in the corresponding test cases.

\paragraph{Topology network.}
The phase-field variable $\phi_{\bm{\sigma}}(\bm{x})\in[0,1]$ is parameterized by a Fourier-feature neural network. For a given spatial dimension $d$ and prescribed numbers of Fourier features $\{n_j^\phi\}_{j=1}^d$ along each coordinate direction, one-dimensional frequency samples are first constructed independently for each direction. Specifically, given a bandwidth parameter $\omega_{\max}>0$, the frequencies along the $j$-th direction are sampled symmetrically from the interval $(-\omega_{\max},0)\cup(0,\omega_{\max})$, forming the set
\begin{equation}
    \Omega_j=\{\omega^{(j)}_i\}_{i=1}^{2n_j^\phi}.
    \nonumber
\end{equation}
A $d$-dimensional frequency grid is then constructed via the tensor-product operation
\begin{equation}
    (C_1,\dots,C_d)=\mathrm{meshgrid}(\Omega_1,\dots,\Omega_d),
    \nonumber
\end{equation}
and the resulting frequency matrix is obtained by flattening each component,
\begin{equation}
    K_\phi=
    \begin{bmatrix}
    \mathrm{vec}(C_1)^\top\\
    \vdots\\
    \mathrm{vec}(C_d)^\top
    \end{bmatrix}
    \in\mathbb{R}^{d\times M_\phi}, \qquad
    M_\phi=\prod_{j=1}^d 2n_j^\phi .
    \nonumber
\end{equation}
The Fourier-feature mapping is defined as
\begin{equation}
    z_\phi(\bm{x})=\sin(\bm{x}^\top K_\phi + b_\phi),
    \label{eq:fourier_feature}
\end{equation}
where $b_\phi\in\mathbb{R}^{M_\phi}$ denotes a phase-shift vector. The phase-field output is obtained through a linear projection followed by a sigmoid activation,
\begin{equation}
    \phi_{\bm{\sigma}}(\bm{x})=\Xi\!\left(z_\phi(\bm{x})W_\phi\right),
\end{equation}
with trainable weights $W_\phi\in\mathbb{R}^{M_\phi\times 1}$ initialized to zero and no bias term,
resulting in a spatially uniform initial phase field $\phi_{\bm{\sigma}}(\bm{x})\equiv 0.5$. Here, the parameter set $\bm{\sigma}=\{K_\phi,W_\phi,b_\phi\}$ collects all trainable variables of the topology network.

\paragraph{Linear elasticity problems.}
For compliance benchmarks, the displacement field $\bm{u}_{\bm{\theta}}(\bm{x})\in\mathbb{R}^d$
is represented by a single-layer Fourier-feature network with an independent set of trainable parameters. The displacement network follows the same Fourier-feature construction as the topology network, but employs its own frequency matrix $K_{\bm{u}}$, phase shift $b_{\bm{u}}$, and feature resolution
$\{n_j^{\bm{u}}\}_{j=1}^d$.
Specifically,
\begin{equation}
    z_{\bm{u}}(\bm{x})=\sin(\bm{x}^\top K_{\bm{u}} + b_{\bm{u}}), \qquad
    \tilde{\bm{u}}_{\bm{\theta}}(\bm{x})=z_{\bm{u}}(\bm{x})\,W_{\bm{u}},
    \nonumber
\end{equation}
where $W_{\bm{u}}\in\mathbb{R}^{M_{\bm{u}}\times d}$ is a trainable weight matrix with
$M_{\bm{u}}=\prod_{j=1}^d 2n_j^{\bm{u}}$. Homogeneous Dirichlet boundary conditions are enforced analytically by multiplying the raw network output with a scalar factor $g(\bm{x})$, yielding
\begin{equation}
    \bm{u}_{\bm{\theta}}(\bm{x}) = g(\bm{x})\,\tilde{\bm{u}}_{\bm{\theta}}(\bm{x}).
\end{equation}
If no such boundary encoding is required, we simply set $g(\bm{x})\equiv 1$. The weight matrix $W_{\bm{u}}$ is initialized to zero and updated during training.

For the linear-elastic eigenvalue benchmarks, the topology network remains unchanged, while the displacement eigenmode
$\bm{u}_{\bm{\theta}}(\bm{x})$ is modeled by an FNN consisting of one Fourier layer followed by an additional hidden layer. Instead of using a tensor-product frequency grid, the Fourier layer employs a trainable Fourier kernel $K_{\bm{u}}\in\mathbb{R}^{d\times M_f}$ and a phase vector $\bm{b}^{(1)}\in\mathbb{R}^{M_f}$, where $M_f$ denotes the number of Fourier features. The kernel and phase shift are initialized as
\begin{equation}
  (K_{\bm{u}})_{ij} \sim \mathcal{U}(-\omega_{\max}^{\bm{u}},\omega_{\max}^{\bm{u}}),
  \qquad
  b^{(1)}_j \sim \mathcal{U}(0,2\pi),
  \label{eq:high_band}
\end{equation}
and kept trainable throughout optimization. The Fourier layer computes
\begin{equation}
  \bm{z}^{(1)}(\bm{x})
  =
  \sin\!\bigl(\bm{x}^\top K_{\bm{u}} + \bm{b}^{(1)}\bigr)
  \in\mathbb{R}^{M_f},
\end{equation}
followed by a dense hidden layer with sine activation,
\begin{equation}
  \bm{z}^{(2)}(\bm{x})
  =
  \sin\!\bigl(\bm{z}^{(1)}(\bm{x})\,W^{(1)} + \bm{b}^{(2)}\bigr)
  \in\mathbb{R}^{M_h},
\end{equation}
where $M_h$ denotes the number of hidden units in the hidden layer, $W^{(1)}\in\mathbb{R}^{M_f\times M_h}$ is initialized by a Glorot/Xavier scheme and $\bm{b}^{(2)}\in\mathbb{R}^{M_h}$ is initialized to zero. The displacement output is then obtained by
\begin{equation}
  \bm{u}_{\bm{\theta}}(\bm{x})
  =
  g(\bm{x})\,\bm{z}^{(2)}(\bm{x})\,W^{(2)},
\end{equation}
with $W^{(2)}\in\mathbb{R}^{M_h\times d}$ the output weight matrix.
Homogeneous Dirichlet boundary conditions on clamped boundaries are enforced by the smooth scalar factor $g(\bm{x})$,
chosen such that $g(\bm{x})=0$ on the constrained segments and $g(\bm{x})\approx 1$ in the interior (e.g.\ a scaled
$\tanh$-type transition function).

\paragraph{Stokes and Navier--Stokes problems.}
For the Stokes and Navier--Stokes benchmarks, the topology network remains identical to that introduced previously. The Stokes state network parameterizes the coupled velocity--pressure field
$\bm{s}_{\bm{\theta}}(\bm{x})=(\bm{u}(\bm{x}),p(\bm{x}))$,
which takes values in $\mathbb{R}^{d+1}$, using Gaussian random Fourier features followed by a sine-activated multilayer perceptron. Given an even feature dimension $M_{\mathrm{Fourier}}$ and a scale parameter $\sigma_{\mathrm{Fourier}}>0$, a fixed projection matrix $B\in\mathbb{R}^{d\times (M_{\mathrm{Fourier}}/2)}$ is sampled as
\begin{equation}
  B_{ij}\sim \mathcal{N}(0,\sigma_{\mathrm{Fourier}}^2),
\end{equation}
and the lifted input is defined by
\begin{equation}
  \bm{z}^{(0)}(\bm{x})
  =
  \bigl[\sin(B^{\top}\bm{x}),\;\cos(B^{\top}\bm{x})\bigr]
  \in\mathbb{R}^{M_{\mathrm{Fourier}}}.
  \label{eq:MFourierStokes}
\end{equation}
This feature vector is processed by $L$ fully connected layers with sine activation,
\begin{equation}
  \bm{z}^{(\ell)}(\bm{x})
  =
  \sin\!\bigl(\bm{z}^{(\ell-1)}(\bm{x})W^{(\ell)}+\bm{b}^{(\ell)}\bigr),
  \qquad \ell=1,\dots,L,
\end{equation}
and the final state output is
\begin{equation}
  \bm{s}_{\bm{\theta}}(\bm{x})
  =
  \bm{z}^{(L)}(\bm{x})W^{(L+1)}+\bm{b}^{(L+1)}.
\end{equation}
All dense-layer parameters are trainable, while the projection matrix $B$ is kept fixed. Velocity Dirichlet conditions (e.g.\ no-slip walls and prescribed inflow) are enforced weakly through boundary penalty terms, and stress-free boundaries are imposed via Neumann residuals.

For the Navier--Stokes benchmarks, the same Gaussian Fourier feature--MLP architecture is adopted for the state network. When adjoint-based sensitivities are required, an adjoint velocity--pressure network with identical structure and an independent parameter set is introduced, and trained jointly with the topology network within the APF-FNN framework.

\subsection{Physics-driven losses for elasticity and flow problems} \label{subsec:physics-loss}

In the APF-FNN framework, each physics problem is characterized by two complementary loss functions: a state loss $\mathcal{L}_{\text{state}}$ used to train the neural approximation of the physical fields, and a topology loss $\mathcal{L}_{\text{topology}}$ used to update the phase-field network. Depending on whether the governing equations admit a variational formulation, $\mathcal{L}_{\text{state}}$ is constructed either from an energy functional in a Deep Ritz fashion \cite{ew2018deep} or from strong-form PDE residuals in a PINN-style formulation \cite{raissi2019physics}. The topology loss has a unified structure across all cases and combines the physical objective with the Ginzburg--Landau regularization and a volume-penalty term,
\begin{equation}
  \mathcal{L}_{\text{topology}}(\bm{\sigma}; \bm{u}_{\bm{\theta}})
  =
  \mathcal{J}_{\text{obj}}(\phi_{\bm{\sigma}}, \bm{u}_{\bm{\theta}})
  + \mathcal{E}_{\mathrm{GL}}(\phi_{\bm{\sigma}})
  + \mathcal{P}_{\mathrm{vol}}(\phi_{\bm{\sigma}}),
  \label{eq:Ltopology-generic}
\end{equation}
where $\phi_{\bm{\sigma}}$ is the phase field produced by the topology network and $\bm{u}_{\bm{\theta}}$ denotes the state network. While the loss constructions differ across physics problems, their implications for design sensitivity computation vary fundamentally, as discussed in Section~\ref{subsec:gradients}. In the following, we detail the specific forms of $\mathcal{L}_{\text{state}}$ and $\mathcal{J}_{\text{obj}}$ for four representative problems: compliance minimization and first eigenvalue maximization in linear elasticity, Stokes flow optimization, and Navier--Stokes flow optimization.

\subsubsection{Compliance minimization in linear elasticity}
\label{subsubsec:physics-comp}

The linear elasticity problem is self-adjoint and derives from the principle of minimum potential energy. The material distribution is modeled via a SIMP-type interpolation of the elasticity tensor,
\begin{equation}
  \mathbb{C}(\phi_{\bm{\sigma}})
  =
  \rho(\phi_{\bm{\sigma}})\,\mathbb{C}_0,
  \qquad
  \rho(\phi_{\bm{\sigma}})
  =
  (1-\phi_{\min})\,\phi_{\bm{\sigma}}^{\,p} + \phi_{\min},
  \nonumber
\end{equation}
where $p>1$ is the penalization exponent and $\phi_{\min}>0$ is a small stiffness lower bound. The tensor $\mathbb{C}_0$ denotes the isotropic linear-elastic stiffness corresponding to Young's modulus $E$ and Poisson ratio $\nu$, acting on a symmetric strain tensor $\boldsymbol{\varepsilon}$ as
\begin{equation}
  \mathbb{C}_0 \boldsymbol{\varepsilon}
  =
  2\mu\,\boldsymbol{\varepsilon}
  + \lambda\,\mathrm{tr}(\boldsymbol{\varepsilon})\,\mathrm{I},
  \qquad
  \mu = \frac{E}{2(1+\nu)}, \quad
  \lambda = \frac{E\nu}{(1+\nu)(1-2\nu)},
  \nonumber
\end{equation}
with Lam\'e parameters $(\lambda,\mu)$ and identity tensor $\mathrm{I}$.

For a fixed topology $\phi_{\bm{\sigma}}$, the displacement field $\bm{u}_{\bm{\theta}}$ is obtained by minimizing the total potential energy of the elastic body. This leads to the Deep Ritz-type state loss
\begin{equation}
  \mathcal{L}_{\text{state}}(\bm{\theta}; \phi_{\bm{\sigma}})
  =
  \frac{1}{2}\int_\Omega
    \mathbb{C}(\phi_{\bm{\sigma}})
    \,\boldsymbol{\varepsilon}(\bm{u}_{\bm{\theta}}) : \boldsymbol{\varepsilon}(\bm{u}_{\bm{\theta}})
    \,\mathrm{d}x
  - \int_\Omega \mathbf{b} \cdot \bm{u}_{\bm{\theta}} \,\mathrm{d}x
  - \int_{\Gamma_N} \mathbf{t} \cdot \bm{u}_{\bm{\theta}} \,\mathrm{d}s
  + \lambda_{\mathrm{Dir}} \int_{\Gamma_D} |\bm{u}_{\bm{\theta}} - \bm{u}_0|^2 \,\mathrm{d}s,
  \label{eq:Lstate-elasticity}
\end{equation}
where $\boldsymbol{\varepsilon}(\bm{u}_{\bm{\theta}}) = (\nabla \bm{u}_{\bm{\theta}}
+ \nabla \bm{u}_{\bm{\theta}}^\top)/2$ is the strain tensor, $\mathbf{b}$ and $\mathbf{t}$ denote body forces and surface tractions, and the last term penalizes deviations from the prescribed displacement $\bm{u}_0$ on the Dirichlet boundary $\Gamma_D$.

The topology network is trained to minimize the structural compliance under a volume constraint. Using the converged state field $\bm{u}_{\bm{\theta}}$, the topology loss
specializes \eqref{eq:Ltopology-generic} to
\begin{equation}
  \mathcal{L}_{\text{topology}}(\bm{\sigma}; \bm{u}_{\bm{\theta}})
  =
  \left(
    \int_\Omega \bm{b} \cdot \bm{u}_{\bm{\theta}} \,\mathrm{d}x
    + \int_{\Gamma_N} \bm{t} \cdot \bm{u}_{\bm{\theta}} \,\mathrm{d}s
  \right)
  + \mathcal{E}_{\mathrm{GL}}(\phi_{\bm{\sigma}})
  + \mathcal{P}_{\mathrm{vol}}(\phi_{\bm{\sigma}}),
  \label{eq:Ltopology-elasticity}
\end{equation}
where the first term is the classical compliance objective and the regularization terms are
defined as in Section~\ref{sec:pf-to}. The volume functional takes the quadratic penalty form
\begin{equation}
  \mathcal{P}_{\mathrm{vol}}(\phi_{\bm{\sigma}})
  =
  \lambda_{\text{penal}}
  \left(
    \int_\Omega \phi_{\bm{\sigma}} \,\mathrm{d}x - \beta |\Omega|
  \right)^2,
\end{equation}
with target volume fraction $0<\beta<1$ and penalty parameter $\lambda_{\text{penal}}$.

\subsubsection{First eigenvalue maximization in linear elasticity}
\label{subsubsec:physics-eig}

The second problem aims at maximizing the fundamental natural frequency of a linear-elastic structure, which is equivalent to maximizing the first eigenvalue $\lambda_1$ of the generalized eigenvalue problem. The formulation is based on the Rayleigh quotient, i.e., the ratio between the elastic strain energy and the kinetic energy. For a fixed topology $\phi_{\bm{\sigma}}$, the fundamental mode $\bm{u_\theta}$ is obtained as the minimizer of this quotient under a unit mass
normalization.

To approximate this eigenmode, the state network is trained to minimize
\begin{equation}
  \mathcal{L}_{\text{state}}(\boldsymbol{\theta}; \phi_{\boldsymbol{\sigma}}) = \frac{\int_{\Omega} \mathbb{C}(\phi_{\boldsymbol{\sigma}})\boldsymbol{\varepsilon}(\bm{u}_{\boldsymbol{\theta}})
    :\boldsymbol{\varepsilon}(\bm{u}_{\boldsymbol{\theta}}) \, \mathrm{d}x}{\int_{\Omega} \varrho(\phi_{\boldsymbol{\sigma}}) |\bm{u}_{\boldsymbol{\theta}}|^2 \, \mathrm{d}x}+ \lambda_{\text{norm}} \left( \int_{\Omega} \varrho(\phi_{\boldsymbol{\sigma}}) |\bm{u}_{\boldsymbol{\theta}}|^2 \, \mathrm{d}x - 1 \right)^2,
  \label{eq:Lstate-eig}
\end{equation}
where \(\varrho(\phi_{\boldsymbol{\sigma}}) = (1-\phi_{\text{min}})\phi_{\boldsymbol{\sigma}} + \phi_{\text{min}}\) is the interpolated mass density and $\lambda_{\text{norm}}$ is the penalty weight enforcing the mass constraint $\int_{\Omega} \varrho(\phi_{\boldsymbol{\sigma}}) \bm{u}_{\boldsymbol{\theta}}^2 \, \mathrm{d}x=1$.

Given the converged eigenmode $\bm{u_\theta}$, the topology network is updated to maximize $\lambda_1$ by minimizing the topology loss
\begin{equation}
  \mathcal{L}_{\text{topology}}(\boldsymbol{\sigma}; \bm{u}_{\boldsymbol{\theta}}) = - \frac{\int_{\Omega} \mathbb{C}(\phi_{\boldsymbol{\sigma}}) \boldsymbol{\varepsilon}(\bm{u}_{\boldsymbol{\theta}}) : \boldsymbol{\varepsilon}(\bm{u}_{\boldsymbol{\theta}}) \, \mathrm{d}x}{\int_{\Omega} \varrho(\phi_{\boldsymbol{\sigma}}) |\bm{u}_{\boldsymbol{\theta}}|^2 \, \mathrm{d}x} + \mathcal{E}_{GL}(\phi_{\boldsymbol{\sigma}}) + \mathcal{P}_{\text{vol}}(\phi_{\boldsymbol{\sigma}}),
  \label{eq:Ltop-eig}
\end{equation}
which combines the negative Rayleigh quotient with the Ginzburg--Landau regularization and the volume-penalty term. In the alternating optimization, $\mathcal{L}_{\text{state}}$ and $\mathcal{L}_{\text{topology}}$ are minimized successively to obtain both the eigenmode and the optimal material layout.

%---------------------------------------------------------------------%
\subsubsection{Stokes flow optimization}
\label{subsubsec:physics-stokes}

Topology optimization for Stokes flow seeks to minimize the total power dissipation of a viscous incompressible fluid, defined as the sum of viscous dissipation and the work done by body forces \cite{jin2025crouzeix}. The optimization is constrained by the steady Stokes equations, a saddle-point system for the velocity and pressure. A PINN formulation provides a natural way to enforce these equations and the associated boundary conditions.

For a fixed phase field $\phi_{\bm{\sigma}}$, a single state network with parameters $\bm{\theta}$ is used to approximate both the velocity and pressure fields $(\bm{u_\theta},p_{\bm{\theta}})$. The network is trained by minimizing a loss composed of four residual terms corresponding to the momentum balance, incompressibility, Dirichlet and Neumann boundary conditions,
\begin{equation}
\begin{aligned}
    \mathcal{L}_{\text{state}}(\boldsymbol{\theta}; \phi_{\boldsymbol{\sigma}}) = & 
    \int_{\Omega} |-\rho\Delta\bm{u}_{\boldsymbol{\theta}} + \nabla p_{\boldsymbol{\theta}} + \Pi(\phi_{\boldsymbol{\sigma}})\bm{u}_{\boldsymbol{\theta}} - \mathbf{f}|^2 \, \mathrm{d}x \\
    & + \lambda_{\text{div}} \int_{\Omega} (\nabla \cdot \bm{u}_{\boldsymbol{\theta}})^2 \, \mathrm{d}x + \lambda_{\text{Dir}} \int_{\Gamma_D} |\bm{u}_{\boldsymbol{\theta}} - \mathbf{g}|^2 \, \mathrm{d}s \\
    & + \lambda_{\text{Neu}} \int_{\Gamma_N} |\left(-p_{\boldsymbol{\theta}}\mathbf{I} + \rho(\nabla\bm{u}_{\boldsymbol{\theta}} + \nabla\bm{u}_{\boldsymbol{\theta}}^T)\right)\mathbf{n} - \mathbf{h}|^2 \, \mathrm{d}s,
\end{aligned}
\end{equation}
where $\rho$ is the viscosity, $\mathbf{f}$ is a body force, and $\Pi(\phi_{\bm{\sigma}})=\eta(1-\phi_{\bm{\sigma}})^2$ is the Darcy-type interpolation modeling the solid phase. The prescribed velocity $\mathbf{g}$ and traction $\mathbf{h}$ are imposed on the Dirichlet and Neumann boundaries, respectively, and ($\lambda_{\text{div}}$, $\lambda_{\text{Dir}}$, $\lambda_{\text{Neu}}$) weight the different loss components.

The topology network is trained to minimize the power dissipation evaluated with the converged state $(\bm{u_\theta},p_{\bm{\theta}})$. The corresponding topology loss is
\begin{equation}
    \mathcal{L}_{\text{topology}}(\boldsymbol{\sigma}; \bm{u}_{\boldsymbol{\theta}}) = \int_{\Omega} \left( \frac{1}{2}\Pi(\phi_{\boldsymbol{\sigma}})|\bm{u}_{\boldsymbol{\theta}}|^2 + \frac{\rho}{2}|\nabla\bm{u}_{\boldsymbol{\theta}}|^2 - \mathbf{f} \cdot \bm{u}_{\boldsymbol{\theta}} \right) \mathrm{d}x + \mathcal{E}_{GL}(\phi_{\boldsymbol{\sigma}}) + \mathcal{P}_{\text{vol}}(\phi_{\boldsymbol{\sigma}}),
    \label{eq:Ltop-stokes}
\end{equation}
which couples the physical objective with the phase-field regularization and the
volume constraint. Using the pre-computed state field in
\eqref{eq:Ltop-stokes} decouples the gradient calculation from the state
parameters and improves the stability and efficiency of the alternating scheme.

%---------------------------------------------------------------------%
\subsubsection{Navier--Stokes flow optimization}
\label{subsubsec:physics-NS}

When inertial effects become significant, particularly at moderate to high Reynolds numbers, the linear Stokes model is no longer a valid approximation of the fluid dynamics. The flow is then governed by the steady incompressible Navier--Stokes equations. Compared with the Stokes case, the momentum balance acquires the nonlinear convective term, which describes the transport of momentum by the flow itself.

For a fixed phase field $\phi_{\bm{\sigma}}$, we again use a single state network with parameters $\bm{\theta}$ to approximate both the velocity and the pressure fields $(\bm{u_\theta},p_{\bm{\theta}})$. The state loss is obtained by augmenting the Stokes loss with the convective contribution in the momentum residual

\begin{equation}
\begin{aligned}
    \mathcal{L}_{\text{state}}(\boldsymbol{\theta}; \phi_{\boldsymbol{\sigma}}) = & 
    \int_{\Omega} |(\bm{u}_{\boldsymbol{\theta}} \cdot \nabla)\bm{u}_{\boldsymbol{\theta}} -\rho\Delta\bm{u}_{\boldsymbol{\theta}} + \nabla p_{\boldsymbol{\theta}} + \Pi(\phi_{\boldsymbol{\sigma}})\bm{u}_{\boldsymbol{\theta}} - \mathbf{f}|^2 \, \mathrm{d}x \\
    & + \lambda_{\text{div}} \int_{\Omega} (\nabla \cdot \bm{u}_{\boldsymbol{\theta}})^2 \, \mathrm{d}x + \lambda_{\text{Dir}} \int_{\Gamma_D} |\bm{u}_{\boldsymbol{\theta}} - \mathbf{g}|^2 \, \mathrm{d}s \\
    & + \lambda_{\text{Neu}} \int_{\Gamma_N} |\left(-p_{\boldsymbol{\theta}}\mathbf{I} + \rho(\nabla\bm{u}_{\boldsymbol{\theta}} + \nabla\bm{u}_{\boldsymbol{\theta}}^T)\right)\mathbf{n} - \mathbf{h}|^2 \, \mathrm{d}s.
\end{aligned}
\label{eq:Lstate-NS}
\end{equation}
The divergence and boundary terms are identical to those in the Stokes formulation; only the momentum residual is modified. This seemingly small change makes the state equations nonlinear and non self-adjoint, and significantly increases the difficulty of the state solution.

The topology network is trained to minimize the regularized power dissipation evaluated with the converged Navier--Stokes solution. The topology loss is therefore
\begin{equation}
  \mathcal{L}_{\text{topology}}(\boldsymbol{\sigma}; \bm{u}_{\boldsymbol{\theta}}) = \int_{\Omega} \left( \frac{1}{2}\Pi(\phi_{\boldsymbol{\sigma}})|\bm{u}_{\boldsymbol{\theta}}|^2 + \frac{\rho}{2}|\nabla\bm{u}_{\boldsymbol{\theta}}|^2 - \mathbf{f} \cdot \bm{u}_{\boldsymbol{\theta}} \right) \mathrm{d}x + \mathcal{E}_{GL}(\phi_{\boldsymbol{\sigma}}) + \mathcal{P}_{\text{vol}}(\phi_{\boldsymbol{\sigma}}),
  \label{eq:Ltop-NS}
\end{equation}
which has the same structure as in the Stokes case but uses the Navier--Stokes state. During the topology update step, $\mathcal{L}_{\text{topology}}$ is minimized with respect to $\bm{\sigma}$ while treating $\bm{u_\theta}$ as fixed. However, because of the non-self-adjoint nature of the Navier--Stokes equations, the design gradient of $\mathcal{L}_{\text{topology}}$ cannot be obtained from automatic differentiation alone. The adjoint formulation used to compute consistent sensitivities is described in Section~\ref{subsubsec:grad-NS}.

\subsection{Gradient computation and adjoint integration}
\label{subsec:gradients}

In the APF\mbox{-}FNN framework, the topology update is driven by a loss
\(
\mathcal{L}_{\text{topology}}(\bm{\sigma})
\)
that couples the physics objective with phase\mbox{-}field regularization and volume control.
Let $\bm{\sigma}$ denote the parameters of the topology network, which produce the phase\mbox{-}field
$\phi_{\bm{\sigma}}(\bm{x})$. For fixed state network, we write
\begin{equation}
  \mathcal{L}_{\text{topology}}(\bm{\sigma})
  =
  \mathcal{J}_{\text{obj}}(\phi_{\bm{\sigma}},\bm{u}_{\bm{\theta}})
  + \mathcal{E}_{\mathrm{GL}}(\phi_{\bm{\sigma}})
  + \mathcal{P}_{\mathrm{vol}}(\phi_{\bm{\sigma}}),
  \label{eq:Ltopo}
\end{equation}
where $\mathcal{J}_{\text{obj}}$ is the physics objective, $\mathcal{E}_{\mathrm{GL}}$ is the Ginzburg--Landau energy,
and $\mathcal{P}_{\mathrm{vol}}$ is the penalty enforcing the volume fraction.

\subsubsection{Explicit automatic differentiation versus total design derivative}
Standard automatic differentiation (AD) treats the state $\bm{u}_{\bm{\theta}}$ as an output of a
black\mbox{-}box PDE solver that is independent of the design parameters $\bm{\sigma}$. Consequently, AD only differentiates through the explicit computational graph of
$\mathcal{L}_{\text{topology}}$ and yields the gradient
\begin{equation}
  \nabla^{\mathrm{AD}}_{\bm{\sigma}} \mathcal{L}_{\text{topology}}
  =
  \frac{\partial \mathcal{J}_{\text{obj}}(\phi_{\bm{\sigma}},\bm{u_\theta})}{\partial \bm{\sigma}}
  +
  \frac{\mathrm{d}\mathcal{E}_{\mathrm{GL}}(\phi_{\bm{\sigma}})}{\mathrm{d}\bm{\sigma}}
  +
  \frac{\mathrm{d}\mathcal{P}_{\mathrm{vol}}(\phi_{\bm{\sigma}})}{\mathrm{d}\bm{\sigma}}.
  \label{eq:AD-grad-sigma}
\end{equation}
This quantity captures the explicit dependence of $\mathcal{L}_{\text{topology}}$ on $\bm{\sigma}$
through $\phi_{\bm{\sigma}}$, but it ignores the fact that the state $\bm{u_\theta}$ itself is an
implicit function of the design, since it solves a PDE parameterized by $\phi_{\bm{\sigma}}$.

The physically correct update direction is given by the total derivative of
$\mathcal{L}_{\text{topology}}$ with respect to $\bm{\sigma}$, which, by the chain rule, decomposes as
\begin{equation}
  \frac{\mathrm{d}\mathcal{L}_{\text{topology}}}{\mathrm{d}\bm{\sigma}}
  =
  \nabla^{\mathrm{AD}}_{\bm{\sigma}} \mathcal{L}_{\text{topology}}
  +
  \frac{\partial \mathcal{J}_{\text{obj}}}{\partial \bm{u_\theta}}
  \frac{\mathrm{d}\bm{u_\theta}}{\mathrm{d}\bm{\sigma}}.
  \label{eq:total-deriv-sigma}
\end{equation}
The second term encodes the implicit sensitivity through the state field.
A naive application of AD that ignores $\mathrm{d}\bm{u_\theta}/\mathrm{d}\bm{\sigma}$ therefore
produces, in general, an incomplete and potentially misleading gradient.
Below we show how this discrepancy is resolved for each physics problem considered in this work.

\subsubsection{Design gradient for compliance minimization}
For compliance minimization, the objective can be written as
\begin{equation}
  \mathcal{J}_{\text{obj}}(\phi_{\bm{\sigma}},\bm{u_\theta})
  =
  \int_\Omega \boldsymbol{\varepsilon}(\bm{u_\theta}) : \mathbb{C}(\phi_{\bm{\sigma}}) : \boldsymbol{\varepsilon}(\bm{u_\theta})\,\mathrm{d}x,
\end{equation}
where $\boldsymbol{\varepsilon}(\bm{u_\theta})$ is the strain tensor and $\mathbb{C}(\phi_{\bm{\sigma}})$ is the
phase\mbox{-}field\mbox{-}dependent elasticity tensor.
The adjoint method provides the true sensitivity of the compliance with respect to the phase\mbox{-}field,
which is proportional to the negative strain energy density and can be expressed as
\begin{equation}
  \frac{\mathrm{d}\mathcal{J}_{\text{obj}}}{\mathrm{d}\phi_{\bm{\sigma}}}
  =
  - \int_\Omega
  \frac{\partial \mathbb{C}(\phi_{\bm{\sigma}})}{\partial \phi_{\bm{\sigma}}}
  \,\boldsymbol{\varepsilon}(\bm{u_\theta}) : \boldsymbol{\varepsilon}(\bm{u_\theta})\,\mathrm{d}x.
  \label{eq:compl-true}
\end{equation}
In contrast, if one applies AD directly to the compliance term in $\mathcal{L}_{\text{topology}}$,
the resulting derivative with respect to $\phi_{\bm{\sigma}}$ is
\begin{equation}
  \frac{\partial \mathcal{J}_{\text{obj}}}{\partial \phi_{\bm{\sigma}}}
  =
  \int_\Omega
  \frac{\partial \mathbb{C}(\phi_{\bm{\sigma}})}{\partial \phi_{\bm{\sigma}}}
  \,\boldsymbol{\varepsilon}(\bm{u_\theta}) : \boldsymbol{\varepsilon}(\bm{u_\theta})\,\mathrm{d}x,
  \label{eq:compl-AD}
\end{equation}
which has the opposite sign. Combining \eqref{eq:compl-true} and
\eqref{eq:compl-AD} reveals the relationship
\begin{equation}
  \frac{\mathrm{d}\mathcal{J}_{\text{obj}}}{\mathrm{d}\phi_{\bm{\sigma}}}
  =
  -\,\frac{\partial \mathcal{J}_{\text{obj}}}{\partial \phi_{\bm{\sigma}}}.
  \label{eq:compl-sign}
\end{equation}
Hence, a naive AD gradient would actually drive the design towards maximizing
compliance, moving the topology in the wrong direction.

In APF\mbox{-}FNNs, we adopt a hybrid strategy.
AD is used as an efficient tool to evaluate the strain energy density
appearing in \eqref{eq:compl-AD}, but the gradient passed to the topology
network is replaced by the physically derived expression \eqref{eq:compl-true}
via the sign flip in \eqref{eq:compl-sign}.
The phase\mbox{-}field and volume-penalty contributions are still computed by AD.
This custom gradient enforces a physically consistent descent direction while retaining
the numerical efficiency of AD.

\subsubsection{Design gradient for first eigenvalue maximization}
For the eigenvalue problems, the first eigenpair $(\lambda_1,\bm{u_\theta})$ is obtained
from a Rayleigh\mbox{-}quotient formulation with phase\mbox{-}field\mbox{-}dependent stiffness
$\mathbb{C}(\phi_{\bm{\sigma}})$ and mass density $\varrho(\phi_{\bm{\sigma}})$.
For a mass\mbox{-}normalized eigenmode, satisfying
\(
\int_\Omega \varrho(\phi_{\bm{\sigma}}) |\bm{u_\theta}|^2\,\mathrm{d}x = 1,
\)
the physical sensitivity of $\lambda_1$ with respect to the phase field reads
\begin{equation}
  \frac{\mathrm{d}\lambda_1}{\mathrm{d}\phi_{\bm{\sigma}}} = \int_{\Omega} \frac{\partial \mathbb{C}(\phi_{\bm{\sigma}})}{\partial \phi_{\bm{\sigma}}} \boldsymbol{\varepsilon}(\bm{u}_{\bm{\theta}}) : \boldsymbol{\varepsilon}(\bm{u}_{\bm{\theta}}) \, \mathrm{d}x - \lambda_1 \int_{\Omega} \frac{\partial \varrho(\phi_{\bm{\sigma}})}{\partial \phi_{\bm{\sigma}}} |\bm{u}_{\bm{\theta}}|^2 \, \mathrm{d}x.
  \label{eq:eig-true}
\end{equation}
A key property of the Rayleigh quotient is that $\lambda_1$ is stationary with respect
to first\mbox{-}order variations of the eigenmode $\bm{u_\theta}$; in other words,
the functional derivative of $\lambda_1$ with respect to $\bm{u_\theta}$ vanishes at the
converged eigenfunction. Consequently, the implicit term in
\eqref{eq:total-deriv-sigma} drops out, and the total derivative coincides with
the partial derivative,
\begin{equation}
  \frac{\mathrm{d}\lambda_1}{\mathrm{d}\phi_{\bm{\sigma}}}
  =
  \frac{\partial \lambda_1}{\partial \phi_{\bm{\sigma}}}.
  \label{eq:eig-equal}
\end{equation}
This identity has an important practical implication: the gradient of the Rayleigh
quotient computed by standard AD is already the true physical sensitivity.
Therefore, for eigenvalue maximization we directly use the AD gradient of the
Rayleigh\mbox{-}quotient loss without any custom modification, and the topology network
is updated with this exact ascent direction.

\subsubsection{Design gradient for Stokes flow optimization}
For Stokes flow, the objective is chosen as the power dissipation, which admits a
variational formulation. The governing equations can be recovered as the
Euler--Lagrange equations of $\mathcal{J}_{\text{obj}}$. In particular, for a converged state
$(\bm{u_\theta},p_{\bm{\theta}})$, the functional derivative of the objective with respect to the
velocity field is proportional to the Stokes residual and vanishes identically,
\begin{equation}
  \frac{\delta \mathcal{J}_{\text{obj}}}{\delta \bm{u_\theta}}
  =
  -\rho \Delta \bm{u_\theta}
  + \nabla p_{\bm{\theta}}
  + \Pi(\phi_{\bm{\sigma}}) \bm{u_\theta}
  - \mathbf{f}
  = 0
  \quad \text{in } \Omega,
  \label{eq:stokes-EL}
\end{equation}
with suitable boundary conditions. Here $\rho$ is the viscosity, $\mathbf{f}$ is the body force
and $\Pi(\phi_{\bm{\sigma}})$ is the Darcy\mbox{-}type interpolation.

Substituting \eqref{eq:stokes-EL} into the chain rule for the derivative
$\mathrm{d}\mathcal{J}_{\text{obj}}/\mathrm{d}\phi_{\bm{\sigma}}$ yields
\begin{equation}
  \frac{\mathrm{d}\mathcal{J}_{\text{obj}}}{\mathrm{d}\phi_{\bm{\sigma}}}
  =
  \frac{\partial \mathcal{J}_{\text{obj}}}{\partial \phi_{\bm{\sigma}}}
  +
  \Bigg\langle
    \frac{\delta \mathcal{J}_{\text{obj}}}{\delta \bm{u_\theta}},
    \frac{\mathrm{d}\bm{u_\theta}}{\mathrm{d}\phi_{\bm{\sigma}}}
  \Bigg\rangle
  =
  \frac{\partial \mathcal{J}_{\text{obj}}}{\partial \phi_{\bm{\sigma}}},
  \label{eq:stokes-total}
\end{equation}
where $\langle\cdot,\cdot\rangle$ denotes the appropriate $L^2$ inner product.
Thus, the entire implicit term vanishes and the total sensitivity is exactly equal to the
explicit partial derivative. In practice this means that, for Stokes flow,
the gradient delivered by AD coincides with the true design gradient.
We therefore use the raw AD gradient of the power\mbox{-}dissipation loss to update the
topology network, which greatly simplifies the implementation.

\subsubsection{Adjoint-based design gradient for Navier--Stokes flow}
\label{subsubsec:grad-NS}

For Navier--Stokes flow, the power-dissipation objective in \eqref{eq:Ltop-NS} does not admit a variational structure whose Euler--Lagrange equations coincide with the steady Navier--Stokes system. Even at a converged state $(\bm{u_\theta},p_{\bm{\theta}})$, the functional derivative
of the objective with respect to the velocity field does not vanish, 
\begin{equation}
  \frac{\delta \mathcal{J}_{\text{obj}}}{\delta \bm{u_\theta}}
  =
  -\,(\bm{u_\theta} \cdot \nabla)\bm{u_\theta}
  \neq 0.
  \label{eq:NS-functional}
\end{equation}
Consequently, the implicit contribution in the total derivative of $\mathcal{J}_{\text{obj}}$ with respect to the phase field is nonzero,
\begin{equation}
  \frac{\mathrm{d}\mathcal{J}_{\text{obj}}}{\mathrm{d}\phi_{\bm{\sigma}}}
  =
  \frac{\partial \mathcal{J}_{\text{obj}}}{\partial \phi_{\bm{\sigma}}}
  -
  \Bigg\langle
    (\bm{u_\theta} \cdot \nabla)\bm{u_\theta},
    \frac{\mathrm{d}\bm{u_\theta}}{\mathrm{d}\phi_{\bm{\sigma}}}
  \Bigg\rangle,
  \label{eq:NS-total-chain}
\end{equation}
which means that the gradient provided by standard automatic differentiation $\partial \mathcal{J}_{\text{obj}}/\partial \phi_{\bm{\sigma}}$ is incomplete. Recovering the missing term by differentiating through the PDE solver directly would be prohibitively expensive, so we resort to a continuous adjoint formulation.

Let $R_{\mathrm{NS}}(\phi_{\bm{\sigma}},\bm{u_\theta},p_{\bm{\theta}})=0$ denote the steady Navier--Stokes residual, including momentum balance, incompressibility and boundary conditions. We introduce an adjoint velocity--pressure pair $(\bm{w_{\alpha}},q_{\bm{\alpha}})$ and define the Lagrangian
\begin{equation}
  \mathcal{L}_{\text{NS}}(\phi_{\bm{\sigma}},\bm{u_\theta},\bm{w_\alpha},q_{\bm{\alpha}})
  =
  \mathcal{J}_{\text{obj}}(\phi_\sigma,\bm{u_\theta})
  + \big\langle R_{\mathrm{NS}}(\phi_{\bm{\sigma}},u_\theta), (\bm{w_\alpha},q_{\bm{\alpha}}) \big\rangle.
\end{equation}
where $\langle\cdot,\cdot\rangle$ denotes the $L^2$ inner product over the domain
and boundary. Stationarity of $\mathcal{L}_{\mathrm{NS}}$ with respect to
$(\bm{u_\theta},p_{\bm{\theta}})$ yields the linear adjoint Navier--Stokes system
\begin{equation}
  \begin{cases}
    (\nabla \bm{u_\theta})^{\!\top} \bm{w_{\alpha}}
    - (\bm{u_\theta}\cdot\nabla) \bm{w_{\alpha}}
    - \rho \Delta \bm{w_{\alpha}}
    + \nabla q_{\bm{\alpha}}
    + \Pi(\phi_{\bm{\sigma}})\,\bm{w_{\alpha}}
    + \displaystyle\frac{\delta \mathcal{J}_{\text{obj}}}{\delta \bm{u_\theta}}
    = 0 & \text{in } \Omega, \\[0.4em]
    \nabla\!\cdot \bm{w_{\alpha}} = 0 & \text{in } \Omega,
  \end{cases}
  \label{eq:NS-adjoint-cont}
\end{equation}
supplemented with homogeneous boundary conditions for $(\bm{w_{\alpha}},q_{\bm{\alpha}})$.

The adjoint fields is approximated by a second neural network with
parameters $\bm{\alpha}$, producing $(\bm{w_\alpha},q_{\bm\alpha})$. This adjoint network is
trained by minimizing a physics-informed loss that penalizes the residuals of
\eqref{eq:NS-adjoint-cont} and their boundary conditions,
\begin{align}
  \mathcal{L}_{\text{adjoint}}(\bm{\alpha}; \bm{u_\theta},\phi_{\bm{\sigma}})
  &=
  \int_\Omega
    \lvert
      (\nabla \bm{u_\theta})^{\!\top} \bm{w_\alpha}
      - (\bm{u_\theta}\cdot\nabla) \bm{w_\alpha}
      - \rho \Delta \bm{w_\alpha}
      + \nabla q_{\bm{\alpha}}
      + \Pi(\phi_{\bm{\sigma}}) \bm{w_\alpha}
      + \frac{\delta \mathcal{J}_{\text{obj}}}{\delta \bm{u_\theta}}
    \rvert^2 \mathrm{d}x
  \nonumber\\[0.2em]
  &\quad
  + \lambda_{\text{adj,div}}
    \int_\Omega (\nabla\!\cdot \bm{w_\alpha})^2 \mathrm{d}x
  + \lambda_{\text{adj,BC}}
    \int_{\Gamma_D} \lvert\bm{w_\alpha}\rvert^2 \mathrm{d}s 
    \\
    & + \lambda_{\text{adj,Neu}} \int_{\Gamma_N} \lvert\left(-q_{\boldsymbol{\alpha}}\mathbf{I}+ \rho(\nabla\bm{w}_{\boldsymbol{\alpha}} + \nabla\bm{w}_{\boldsymbol{\alpha}}^{\!\top})\right)\mathbf{n}\rvert^2 \, \mathrm{d}s.
  \label{eq:Ladj-NS}
\end{align}
Once the state and adjoint networks are converged for a given design, the adjoint equations and the chain rule combine to give an explicit expression for the phase-field sensitivity of the physical objective. For the Darcy-type interpolation $\Pi(\phi_{\bm{\sigma}}) = \eta(1-\phi_{\bm{\sigma}})^2$ used in the flow problems, we obtain
\begin{equation}
  \frac{\mathrm{d}\mathcal{J}_{\text{obj}}}{\mathrm{d}\phi_{\bm{\sigma}}}
  =
  \frac{\partial \Pi}{\partial \phi_{\bm{\sigma}}}
  \left(
    \frac{1}{2}\,|\bm{u_\theta}|^2
    - \bm{w}_{\bm{\alpha}} \cdot \bm{u_\theta}
  \right)
  =
  -2\eta(1-\phi_{\bm{\sigma}})
  \left(
    \frac{1}{2}\,|\bm{u_\theta}|^2
    - \bm{w}_{\bm{\alpha}} \cdot \bm{u_\theta}
  \right).
  \label{eq:NS-phi-sens}
\end{equation}

Finally, the gradient of the topology loss with respect to the design parameters $\sigma$ is assembled by combining the objective sensitivity \eqref{eq:NS-phi-sens} with the phase-field regularization and volume-penalty terms,
\begin{equation}
  \nabla_{\bm{\sigma}} \mathcal{L}_{\text{topology}}
  =
  \Bigg(
    \frac{\mathrm{d}\mathcal{J}_{\text{obj}}}{\mathrm{d}\phi_{\bm{\sigma}}}
    +
    \frac{\delta \mathcal{E}_{\mathrm{GL}}}{\delta \phi_{\bm{\sigma}}}
    +
    \frac{\delta \mathcal{P}_{\mathrm{vol}}}{\delta \phi_{\bm{\sigma}}}
  \Bigg)
  \frac{\partial \phi_{\bm{\sigma}}}{\partial {\bm{\sigma}}}.
  \label{eq:topo-grad-final}
\end{equation}
In practice, the derivatives of $\mathcal{E}_{\mathrm{GL}}$ and $\mathcal{P}_{\mathrm{vol}}$ with respect to $\phi_{\bm{\sigma}}$ are evaluated by automatic differentiation, while the adjoint network provides the missing implicit contribution of the Navier--Stokes state. This hybrid strategy combines the efficiency of AD with the physical consistency of adjoint-based sensitivity analysis. These case-dependent analyses illustrate that APF-FNNs provide a unified optimization framework, while enabling problem-specific treatments of design sensitivities to ensure physically correct descent directions.

%======================================================================
% SECTION 4: NUMERICAL EXPERIMENTS AND RESULTS
%======================================================================
\section{Numerical Experiments}
\label{sec:experiments}

We assess the performance and generality of the proposed APF-FNN framework on four classes of PDE-constrained topology optimization problems introduced in Section~\ref{subsec:physics-loss}. The same alternating training schedule in Algorithm~\ref{alg: tri-level to} is used throughout, with only the network roles adapting to the underlying physics. For non-self-adjoint systems, such as Navier--Stokes flow, the full three-network architecture is employed, with separate FNNs for the state, adjoint and topology fields. For self-adjoint problems (compliance, eigenvalue and Stokes), the adjoint equations reduce to the state equations, so the adjoint network and its update stage are omitted, yielding a simpler two-network alternating scheme between the state and topology fields.

All networks are trained with the Adam optimizer following the two-stage alternating procedure summarized in Algorithm~\ref{alg: tri-level to}. Given interior training coordinates
$X\subset\Omega$ and boundary samples $X_b\subset\partial\Omega$, we first perform an initial fitting stage in which only the state network $\bm{u}_{\bm{\theta}}$ is trained for $I$ epochs, with a fixed topology network $\phi_{\bm{\sigma}}$, to approximately satisfy the governing equations using the physics-driven loss $\mathcal{L}_{\text{state}}$. When an adjoint network $\bm{w}_{\bm{\alpha}}$ is required (e.g.\ for non-self-adjoint flow problems), its parameters are initialized by copying the trained weights of the state network, taking advantage of the identical architecture, rather than through a separate pretraining phase. The method then enters the main alternating optimization loop of $N$ outer iterations. At each iteration $n$, all networks are warm-started from iteration $n\!-\!1$ and we successively perform $K$, $L$ and $M$ Adam steps on $\mathcal{L}_{\text{state}}(\bm{\theta};\phi_{\bm{\sigma}})$, $\mathcal{L}_{\text{adjoint}}(\bm{\alpha};\phi_{\bm{\sigma}},\bm{u}_{\bm{\theta}})$ and $\mathcal{L}_{\text{topology}}(\bm{\sigma};\bm{u}_{\bm{\theta}},\bm{w}_{\bm{\alpha}})$, respectively, using problem-specific learning rates for the state, adjoint and topology networks. After each outer loop, the penalty parameter associated with the volume constraint is updated, and convergence is monitored in terms of the objective value and volume fraction.

All domain and boundary integrals in the loss functionals are approximated by equal-weight Monte Carlo quadrature over the sampled collocation points. The training set consists of interior points for enforcing the governing PDEs and boundary points for imposing Dirichlet or Neumann conditions. For the elasticity and eigenvalue benchmarks, interior points are taken from a structured tensor-product grid over the domain, and boundary points are sampled uniformly along the prescribed loaded or clamped segments. For the Stokes and Navier--Stokes benchmarks, interior points are generated by Latin hypercube sampling in the reference square or cube, and boundary points are sampled uniformly on each edge or face.

All experiments are implemented in Python using \texttt{TensorFlow} and executed on a single NVIDIA RTX~4090 GPU. For consistency across benchmarks and ease of reproducibility, we run a fixed number of outer iterations $N$ in all experiments. In all cases, the objective value and volume fraction exhibit stable behavior well before termination. An overview of the network configurations and the associated training budgets is provided in Table~\ref{tab:net_summary}.

\begin{table}[t]
\centering
\caption{Architecture summary of the state/adjoint networks and the topology networks. $d_{\mathrm{in}}$ and $d_{\mathrm{out}}$ are input and output dimensions. $N_F$ is the number of Fourier features. Hidden is reported as width$\times$depth (neurons per layer$\times$number of hidden layers). $N_\theta$ is the number of trainable parameters.}
\label{tab:net_summary}
\setlength{\tabcolsep}{6pt}
\renewcommand{\arraystretch}{1.15}

\begin{tabular}{cccccc cccc}
\toprule
\multirow{2}{*}{\centering Case} &
\multicolumn{5}{c}{State/Adjoint network} &
\multicolumn{4}{c}{Topology network} \\
\cmidrule(lr){2-6} \cmidrule(lr){7-10}
& $d_{\mathrm{in}}$ & $d_{\mathrm{out}}$ & $N_F$ & Hidden (w$\times$d) & $N_\theta$
& $d_{\mathrm{in}}$ & $d_{\mathrm{out}}$ & $N_F$ & $N_\theta$ \\
\midrule
(a)--(c) & 2 & 2 & 2304 & --            & 9216  & 2 & 1 &  576 &  1728 \\
(d)--(f) & 3 & 3 & 4096 & --            & 24576 & 3 & 1 & 3456 & 13824 \\
(g)--(h) & 2 & 2 &   32 & 32$\times$1   & 1216  & 2 & 1 & 1152 &  3456 \\
(i)--(j) & 3 & 3 &   64 & 64$\times$1   & 4608  & 3 & 1 & 3456 & 13824 \\
(k)--(m) & 2 & 3 &  512 & 64$\times$4   & 45507 & 2 & 1 &  144 &   576 \\
(n)--(o) & 3 & 4 &  256 & 128$\times$4  & 82948 & 3 & 1 &  216 &  1080 \\
(p)--(r) & 2 & 3 &  512 & 64$\times$5   & 49667 & 2 & 1 &  144 &   576 \\
(s)      & 3 & 4 &  256 & 128$\times$4  & 82948 & 3 & 1 &  216 &  1080 \\
\bottomrule
\end{tabular}
\end{table}

\begin{algorithm}[htbp]
\caption{APF-FNNs: Alternating Phase-Field Fourier Neural Networks}
\label{alg: tri-level to}

\KwIn{
    \\
    \Indp
    % $\beta$: Target volume fraction\;
    $X$: Training coordinates within the design domain $\Omega$\;
    $X_b$: Boundary coordinates where the penalty terms are applied\;
    $\text{LR}_{\bm{\theta},\text{init}}, \text{LR}_{\bm{\theta},\text{opt}}, 
    \text{LR}_{\bm{\alpha}},
    \text{LR}_{\bm{\sigma}}$: Learning rates\;
    $\lambda_{\text{penal}}^{(0)}, \zeta$: Initial penalty parameter and its scaling factor\;
    $\beta$, $(\epsilon, \gamma)$: Volume fraction and parameters for the phase-field model\;
    $I, N$: Epochs for initial fitting and total optimization\;
    $K, L, M$: Inner-loop iterations for state, adjoint and design updates\;
}
\BlankLine

% --- Initialize Section (un-numbered) ---
% Use \nonl to create an unnumbered line
% \nonl \textbf{Initialize:}\;
\textbf{Initialize:}\\
\Indp
    State network $\bm{u}_{\bm{\theta}}(\bm{x})$ with parameters $\bm{\theta}^{(0)}$\;
    Topology network $\phi_{\bm{\sigma}}(\mathbf{x})$ with parameters $\bm{\sigma}^{(0)}$\;
    %Adam optimizers for three networks\;
\Indm
\BlankLine

% --- Procedure Section ---
% \nonl \textbf{Procedure:}\;
% The first actual step gets a number with \nl
\nl Fix $\bm{\sigma}^{(0)}$ and obtain the initial design $\phi_{\bm{\sigma}^{(0)}}(\mathbf{x})$\;
\tcp{Stage 1: Initial field estimations}
\ForDo{$i = 1$ \bfseries to $I$}{
    Compute loss $\mathcal{L}_{\text{state}}(\boldsymbol{\theta}^{(i-1)}; \phi_{\boldsymbol{\sigma}^{(0)}})$\;
    Update state network: $\bm{\theta}^{(i)} \leftarrow \text{Adam}(\bm{\theta}^{(i-1)}, \nabla_{\bm{\theta}} \mathcal{L}_{\text{state}})$\;
}
\nl $\bm{\theta}^{(0)}_{0} \leftarrow \bm{\theta}^{(I)}$, $\bm{\alpha}^{(0)}_{0} \leftarrow \bm{\theta}^{(I)}$ and $\bm{\sigma}^{(0)}_{0} \leftarrow \bm{\sigma}^{(0)}$\;
\BlankLine

\tcp{Stage 2: Alternating optimization}
\ForDo{$n = 1$ \bfseries to $N$}{
    Update penalty: $\lambda_{\text{penal}}^{(\text{n})} \leftarrow \lambda_{\text{penal}}^{(\text{n}-1)}/\zeta$\;
    
    \tcp{Optimize state $\bm{u}_{\bm{\theta}}$}
    Fix design parameters $\bm{\sigma}^{(0)}_{n-1}$\;
    \ForDo{$k = 1$ \bfseries to $K$}{
        Compute loss $\mathcal{L}_{\text{state}}(\boldsymbol{\theta}^{(k-1)}_{n-1}; \phi_{\bm{\sigma}^{(0)}_{n-1}})$\;
        Update state network parameters:
        $\bm{\theta}^{(k)}_{n-1} \leftarrow \text{Adam}(\bm{\theta}^{(k-1)}_{n-1}, \nabla_{\bm{\theta}} \mathcal{L}_{\text{state}})$\;
        % \Indm
    }
    $\bm{\theta}^{(0)}_{n} \leftarrow \bm{\theta}^{(K)}_{n-1}$\;

    \tcp{Optimize adjoint $\bm{w}_{\bm{\alpha}}$}
    Fix design parameters $\bm{\sigma}^{(0)}_{n-1}, \bm{\theta}^{(0)}_{n}$\;
    \ForDo{$l = 1$ \bfseries to $L$}{
        Compute loss $\mathcal{L}_{\text{adjoint}}(\boldsymbol{\alpha}^{(l-1)}_{n-1}; \phi_{\bm{\sigma}^{(0)}_{n-1}},\bm{\theta}^{(0)}_{n})$\;
        Update adjoint network parameters:
        $\bm{\alpha}^{(l)}_{n-1} \leftarrow \text{Adam}(\bm{\alpha}^{(l-1)}_{n-1}, \nabla_{\bm{\alpha}} \mathcal{L}_{\text{adjoint}})$\;
        % \Indm
    }
    $\bm{\alpha}^{(0)}_{n} \leftarrow \bm{\alpha}^{(L)}_{n-1}$\;
    
    \tcp{Optimize phase field $\phi_{\bm{\sigma}}$}
    \ForDo{$m = 1$ \bfseries to $M$}{
        Compute loss $\mathcal{L}_{\text{topology}}(\boldsymbol{\sigma}^{(m-1)}_{n-1}; \bm{u}_{\boldsymbol{\theta}^{(0)}_{n}},\bm{w}_{\boldsymbol{\alpha}^{(0)}_{n}})$\;
        Update design network via adjoint method:
        $\bm{\sigma}^{(m)}_{n-1} \leftarrow \text{Adam}(\bm{\sigma}^{(m-1)}_{n-1}, \nabla_{\bm{\sigma}} \mathcal{L}_{\text{topology}})$\;
        % \Indm
    }
    $\bm{\sigma}^{(0)}_{n} \leftarrow \bm{\sigma}^{(M)}_{n-1}$\;
}
\BlankLine

\KwOut{The optimized parameters $\bm{\theta}^{(0)}_{N}$, $\bm{\alpha}^{(0)}_{N}$ and $\bm{\sigma}^{(0)}_{N}$.}
\end{algorithm}

\subsection{Compliance minimization in linear elasticity}
\label{subsec: numerical compliance}

In this subsection, we evaluate APF-FNNs on standard compliance-minimization benchmarks in linear elasticity, including three 2D and three 3D test cases. We summarize the common problem settings and implementation details, report the optimized topologies obtained by
APF-FNNs, and compare against a classical FEM-based phase-field solver as a baseline.

\subsubsection{Benchmark setup} \label{sec: linearsetup}

All benchmarks use an isotropic material with Young's modulus $E=1.0$ and Poisson's ratio $\nu=0.3$, SIMP penalization with exponent $p=3$, target volume fractions $\beta=0.5$ in the 2D cases and $\beta=0.3$ in the 3D cases, and a small stiffness floor $\phi_{\min}=10^{-4}$. Body forces are set to zero and any boundary not explicitly constrained is traction-free. The phase-field regularization parameters $(\epsilon,\gamma)$ are both $0.01$, and the volume-penalty scheme follows the general settings described in Section~\ref{subsubsec:physics-comp}.

\begin{figure}[htbp] 
    \centering
    \includegraphics[width=1\linewidth]{ 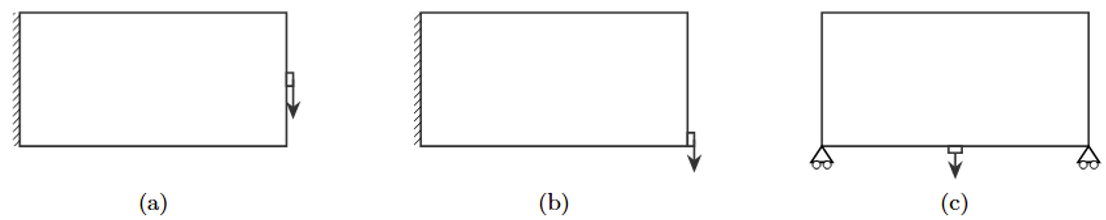} 
    \caption{Problem setups for the 2D compliance minimization benchmarks: (a) Cantilever beam, (b) Offset-loaded beam, and (c) MBB beam.} 
    \label{fig: compliance2dsettings}
\end{figure}

The 2D cases share the rectangular domain $\Omega = (-1,1)\times(-0.5,0.5)$ and differ only by their support and loading conditions as shown in Fig.~\ref{fig: compliance2dsettings}:
\begin{itemize}
    \item[(a)] \textit{Cantilever beam:} A standard cantilever setup where the left boundary is clamped on $\Gamma_D = \{0\} \times (-0.5, 0.5)$, and a downward point load $\mathbf{t} = (0, -1)^{\mathrm{T}}$ is applied to the center of the right boundary on $\Gamma_{N} = \{1\} \times (-0.05, 0.05)$.
    
    \item[(b)] \textit{Offset-loaded beam:} This case modifies the cantilever setup by relocating the load to the bottom-right corner, on the segment $\Gamma_{N} = \{1\} \times (-0.5, -0.4)$, while keeping the left boundary clamped.
    
    \item[(c)] \textit{MBB beam:} A classic Messerschmitt-Bölkow-Blohm (MBB) beam, supported at the two bottom corners on $\Gamma_{D} = (-1, -0.95) \times \{-0.5\} \cup (0.95, 1) \times \{-0.5\}$ and subjected to a downward load at the center of its bottom edge with $\Gamma_{N} = (-0.05, 0.05) \times \{-0.5\}$.
\end{itemize}

The 3D cases are direct spatial extensions of these configurations. All use the same domain $\Omega = (-1,1)\times(-0.5,0.5)\times(-0.2,0.2)$, obtained by extruding the 2D domain in the thickness direction, and adopt clamped and loaded patches that are the natural three-dimensional counterparts of the 2D segments in (a)--(c).

\begin{itemize}
    \item[(d)] \textit{3D Cantilever beam:} The left face $\Gamma_D = \{-1\} \times (-0.5, 0.5) \times (-0.2, 0.2)$ of the domain is clamped, and a downward traction of $\mathbf{t} = (0, -1, 0)^{\mathrm{T}}$ is applied on the point $(1,0,0)$.

    \item[(e)] \textit{3D Offset-loaded beam:} The cantilever setup is modified by relocating the load to the bottom-right corner of the structure. The traction is applied on the point $(1,-0.5,0)$ while the left face remains clamped.

    \item[(f)] \textit{3D MBB beam:} The structure is supported at two regions on the bottom face $\Gamma_{D} = ((-1, -0.9) \cup (0.9, 1)) \times \{-0.5\} \times (-0.05, 0.05)$ and subjected to a downward traction on a central segment of the top face $\Gamma_{N} = (-0.05, 0.05) \times \{0.5\} \times (-0.05, 0.05)$.
\end{itemize}

\subsubsection{APF-FNN implementation details and results}

\begin{table*}[t]
\centering
\caption{Hyperparameters in Algorithm \ref{alg: tri-level to} for the compliance benchmarks cases (a)--(f).}
\label{tab:compliance_hparams_full}
\small
\setlength{\tabcolsep}{4pt}
\renewcommand{\arraystretch}{1.15}
\resizebox{\textwidth}{!}{
\begin{tabular}{c c c c c c c c c c c c}
\toprule
Case &
$N$ &
$K$ &
$|X|$ &
$|X_b|$ &
$(n^{\bm{u}}_1,n^{\bm{u}}_2,n^{\bm{u}}_3)$ &
$\text{LR}_{\bm{\theta},\text{init}}$ &
$\text{LR}_{\bm{\theta},\text{opt}}$ &
$\omega_{\max}^{\bm{u}}$ &
$(n^{\phi}_1,n^{\phi}_2,n^{\phi}_3)$ &
$\omega_{\max}^{\phi}$ &
$\text{LR}_{\bm{\sigma}}$ \\
\midrule
(a) & 500 & 20 & 20000 & 500 & (24,24,0) & $1{e}{-3}$ & $1{e}{-4}$ & 35 & (12,12,0) & 30 & $1{e}{-3}$ \\
(b) & 500 & 20 & 20000 & 500 & (24,24,0) & $1{e}{-3}$ & $1{e}{-4}$ & 35 & (12,12,0) & 30 & $1{e}{-3}$ \\
(c) & 500 & 20 & 20000 & 500 & (24,24,0) & $1{e}{-3}$ & $1{e}{-4}$ & 35 & (12,12,0) & 30 & $5{e}{-4}$ \\
(d) & 500 & 30 & 54000 & 500 & (8,8,8)   & $5{e}{-2}$ & $1{e}{-2}$ & 35 & (12,6,6)  & 25 & $1{e}{-3}$ \\
(e) & 500 & 30 & 54000 & 500 & (8,8,8)   & $5{e}{-2}$ & $1{e}{-2}$ & 35 & (12,6,6)  & 20 & $1{e}{-3}$ \\
(f) & 500 & 30 & 54000 & 500 & (8,8,8)   & $5{e}{-2}$ & $1{e}{-2}$ & 35 & (12,6,6)  & 25 & $1{e}{-3}$ \\
\bottomrule
\end{tabular}%
}
\end{table*}

The network architectures are summarized in Table~\ref{tab:net_summary}. The hyperparameter settings in Algorithm~\ref{alg: tri-level to} are reported in Table~\ref{tab:compliance_hparams_full}. In the initial fitting stage, the state network is trained for at least $I_{\min}=10^3$ and at most $I_{\max}=5\times10^3$ Adam steps with a fixed topology field, and is terminated early once the ratio between the weak-form error and the compliance falls below $10^{-3}$. The subsequent alternating optimization employs $N=500$ outer iterations in all cases. Within each outer iteration, the state network is updated $K$ times and the topology network $M$ times. In the 2D benchmarks, we use $K=20$ and start with $M=1$, increasing to $M=10$ once $n>200$, except for case (b), where $M$ is increased to 20; in the 3D benchmarks we use $K=30$ and increase $M$ from $1$ to $5$ after $n>400$. Interior and boundary collocation points follow the sampling strategy described in previous Section~\ref{sec:experiments}. Concretely, each 2D benchmark uses $|X|= 2\times 10^4$ interior points on a tensor-product grid and $|X_b|=500$ boundary points on the loaded and clamped segments. The 3D benchmarks employ $|X|=5.4\times 10^4$ interior points (a $60\times 30\times 30$ grid) and again $|X_b|=500$ boundary points.

\begin{figure}[htbp] 
    \centering
    \newlength{\rowlabelwidth}
    \settowidth{\rowlabelwidth}{ (a))} 

    \begin{tabular}{
        >{}l 
        @{\hspace{1em}} 
        >{\centering\arraybackslash}m{0.22\textwidth} 
        @{\hspace{0.5em}} 
        >{\centering\arraybackslash}m{0.22\textwidth}
        @{\hspace{0.5em}}
        >{\centering\arraybackslash}m{0.22\textwidth}
        @{\hspace{0.5em}}
        >{\centering\arraybackslash}m{0.22\textwidth}
    }
        & \multicolumn{1}{c}{Epoch 100} 
        & \multicolumn{1}{c}{Epoch 200} 
        & \multicolumn{1}{c}{Epoch 300} 
        & \multicolumn{1}{c}{Epoch 500} \\
        \addlinespace[5pt]
        
        % 第一行图片
        (a) &
        \includegraphics[width=\linewidth, height=3.5cm, keepaspectratio]{ 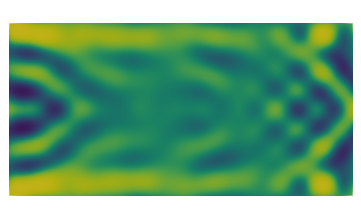} &
        \includegraphics[width=\linewidth, height=3.5cm, keepaspectratio]{ 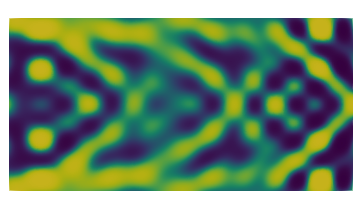} &
        \includegraphics[width=\linewidth, height=3.5cm, keepaspectratio]{ 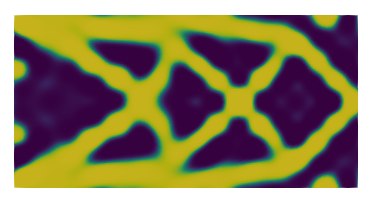} &
        \includegraphics[width=\linewidth, height=3.5cm, keepaspectratio]{ 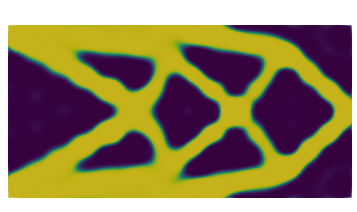} \\
        % \addlinespace[10pt]

        % 第二行图片
        (b) &
        \includegraphics[width=\linewidth, height=3.5cm, keepaspectratio]{ 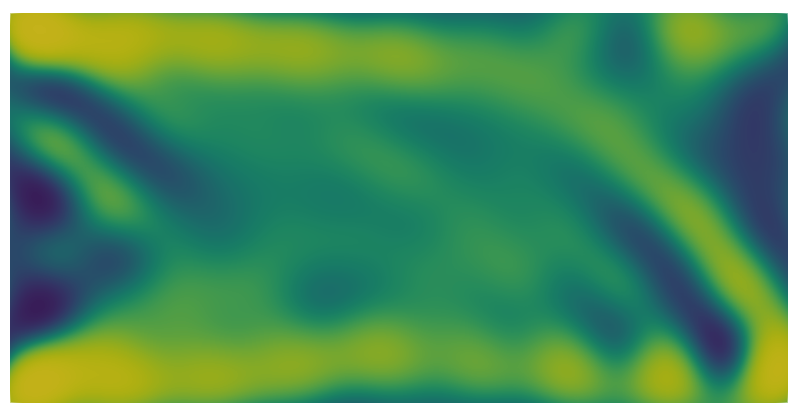} &
        \includegraphics[width=\linewidth, height=3.5cm, keepaspectratio]{ 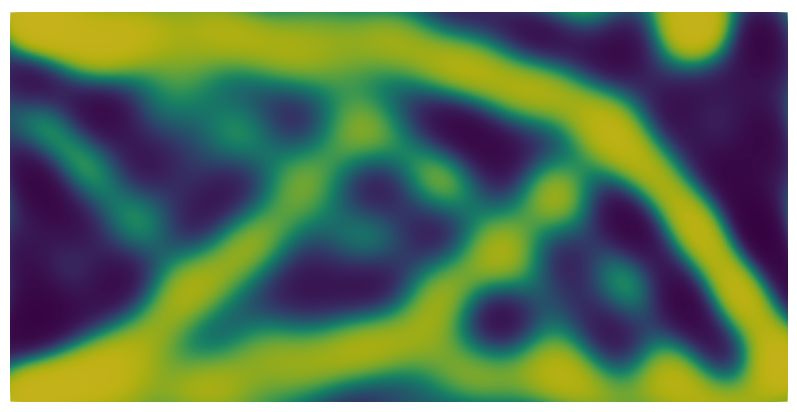} &
        \includegraphics[width=\linewidth, height=3.5cm, keepaspectratio]{ 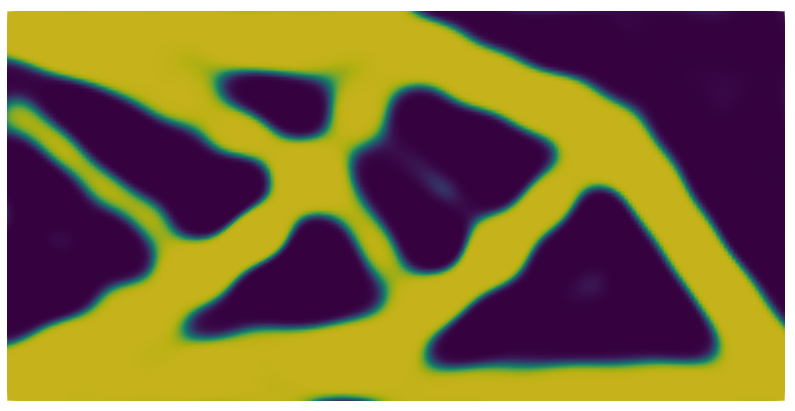} &
        \includegraphics[width=\linewidth, height=3.5cm, keepaspectratio]{ 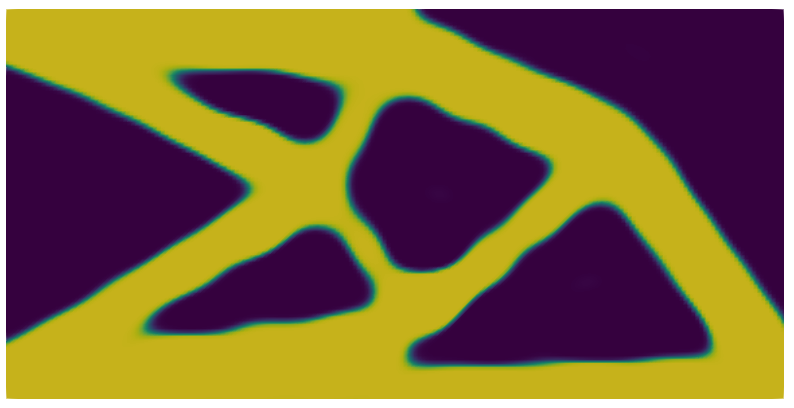} \\

        % 第三行图片
        (c) &
        \includegraphics[width=\linewidth, height=3.5cm, keepaspectratio]{ 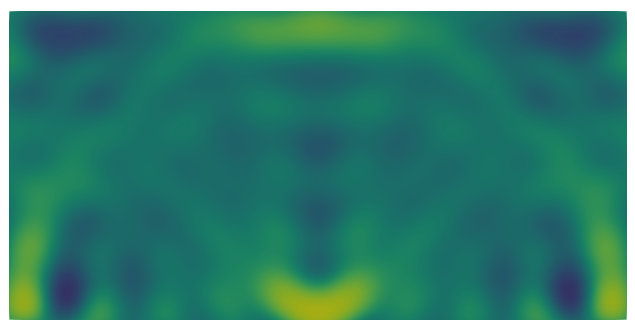} &
        \includegraphics[width=\linewidth, height=3.5cm, keepaspectratio]{ 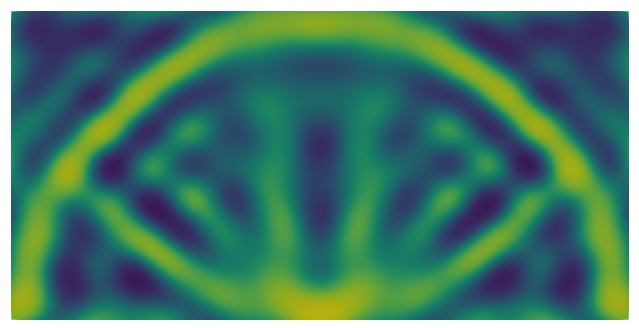} &
        \includegraphics[width=\linewidth, height=3.5cm, keepaspectratio]{ 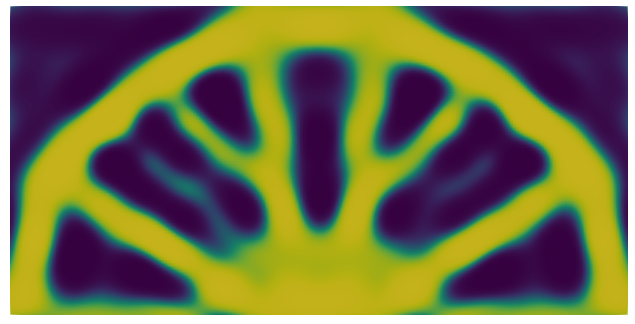} &
        \includegraphics[width=\linewidth, height=3.5cm, keepaspectratio]{ 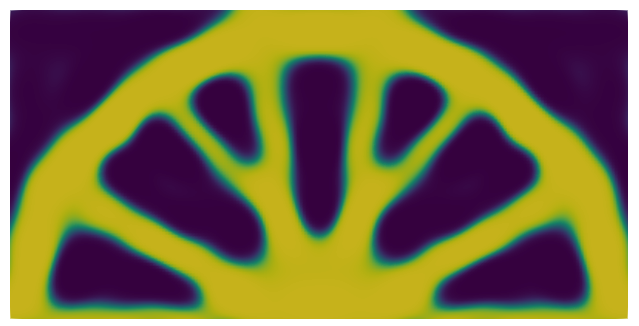}
    \end{tabular}
    
    \caption{Topological evolution for the 2D compliance minimization benchmarks. Each row illustrates the progression from a diffuse and ambiguous configurations (left) to a crisp and well-defined final topologies (right) for cases (a)-(c), respectively.} 
    \label{fig: 2d_compliance_design}
\end{figure}

The band-limited Fourier layers of the state networks in cases (a)–(f) use the same frequency band $\omega^{\bm{u}}_{\max}=35$. For the topology networks, the band $\omega^{\phi}_{\max}$ is $30$ in the 2D cases, while in the 3D cases the upper bounds are $25$ in cases (d)-(e), and $20$ in case (f). In 2D, the state FNN uses $24\times 24$ modes, while the topology FNN uses a slightly narrower basis with $12\times 12$ modes. In 3D, the state network uses $8\times 8\times 8$ modes, and the topology network uses $12\times 12\times 6$ modes in cases (d)–(f). The learning rate of the state network is $10^{-3}$ in the pretraining stage and $10^{-4}$ in the alternating stage for 2D (respectively $5\times 10^{-2}$ and $10^{-2}$ in 3D), while the topology network uses learning rates $10^{-3}$ except for case (c), where the associated learning rate is $5\times10^{-4}$. The penalty parameter $\lambda_{\text{penal}}$ in the volume functional $\mathcal{P}_{\text{vol}}$ is ramped from an $\mathcal{O}(1)$ initial value to $\mathcal{O}(10^3\!-\!10^4)$ using a multiplicative factor $\zeta=0.98$ in 2D and $\zeta=1/1.1$ in 3D, so that the volume constraint is enforced progressively as the topology stabilizes. A soft Dirichlet penalty is required only in the 2D and 3D MBB beams to impose the simply supported boundary with penalty parameters being 5 and 10, respectively ; all other clamped boundaries are enforced analytically through the boundary factor $g(\bm{x})$ in the state network, chosen as a linear function that vanishes on the Dirichlet boundary.

\begin{figure}[htbp] 
    \centering
    \includegraphics[width=0.7\linewidth]{ 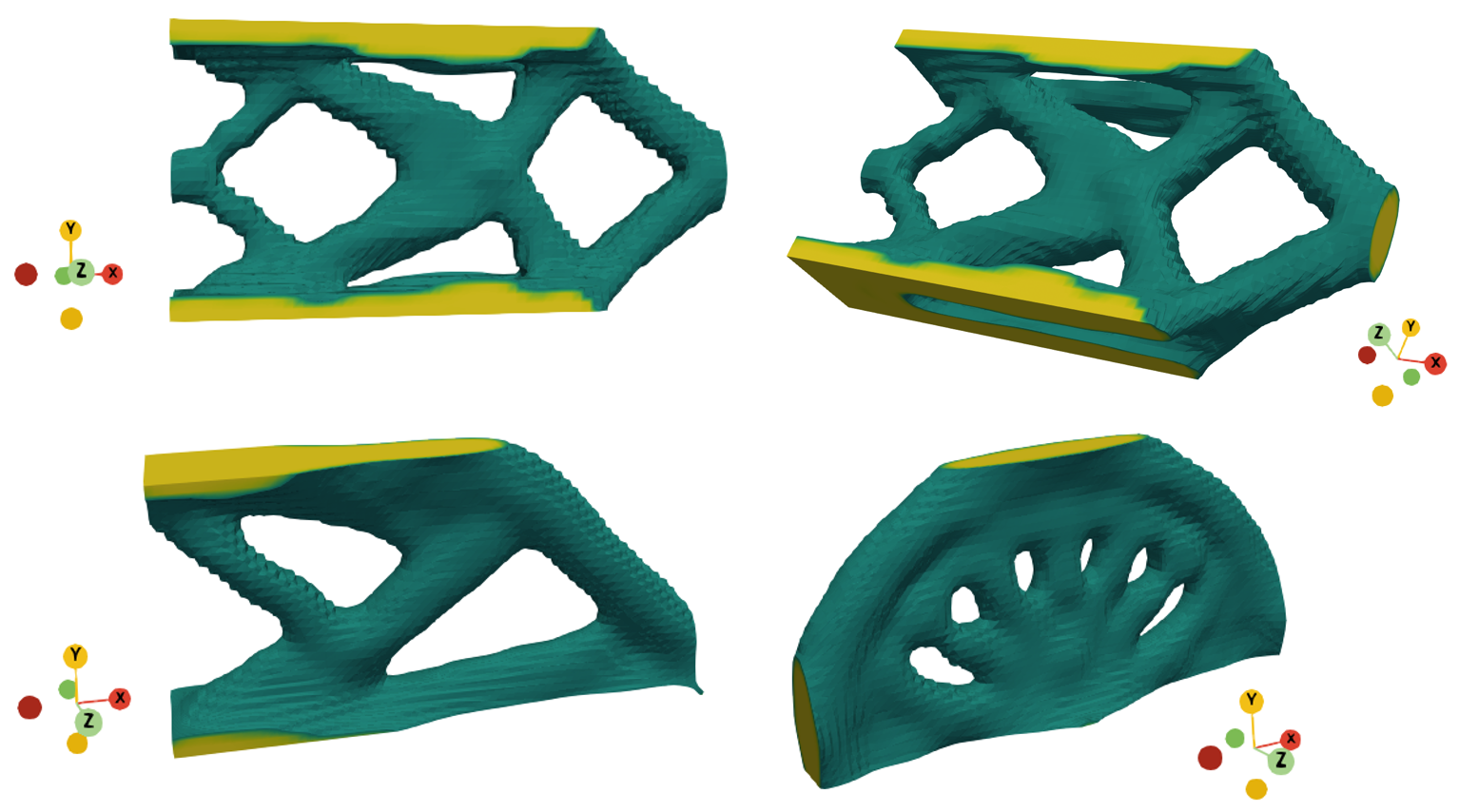}
    \caption{Optimized topologies for the 3D compliance minimization benchmarks. Top row: Two views of case (d). Bottom row: The final designs for cases (e) (left) and (f) (right).} 
    \label{fig: 3d_compliance_design}
\end{figure}

\begin{figure}[htbp] 
    \centering
    \includegraphics[width=0.36\linewidth]{ 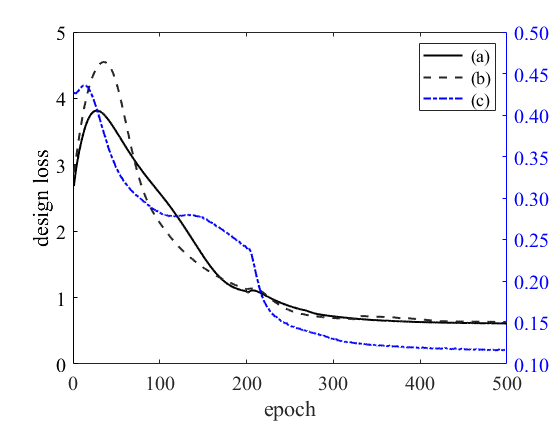}
    \includegraphics[width=0.36\linewidth]{ 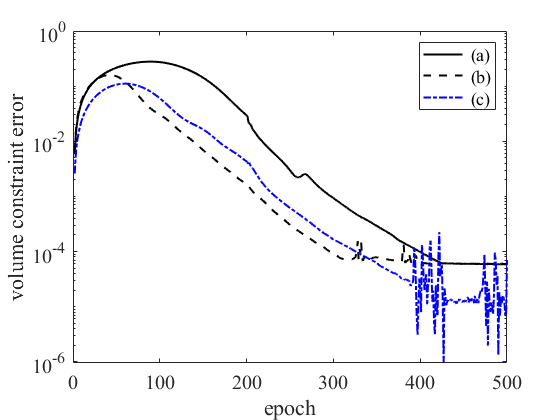}
    \\
    \includegraphics[width=0.36\linewidth]{ 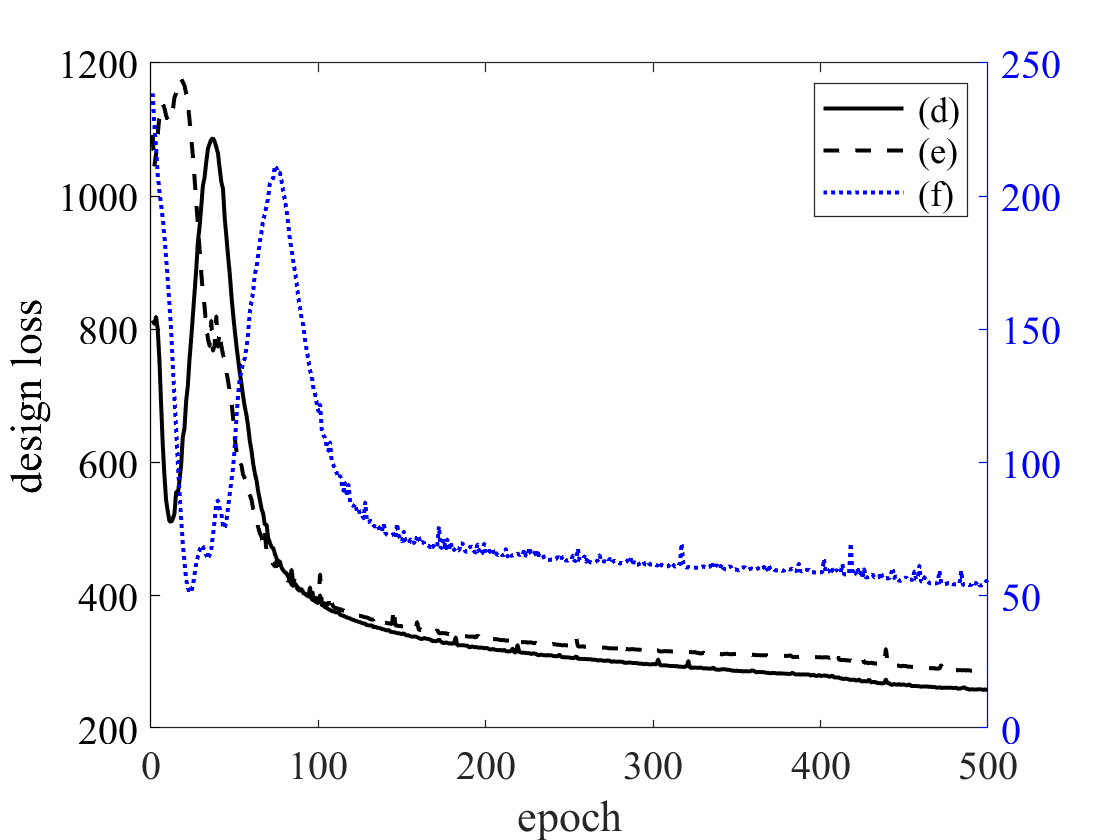}
    \includegraphics[width=0.36\linewidth]{ 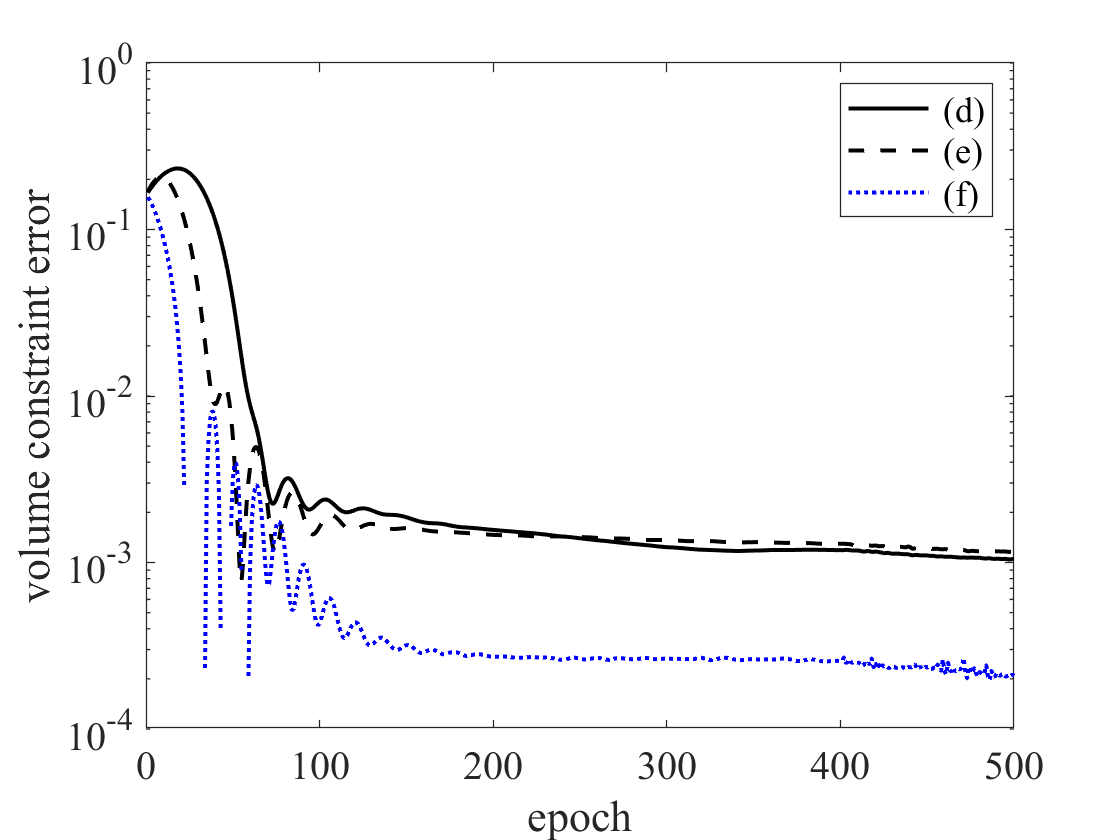}
    \caption{Convergence history of the 2D cases (a)–(c) and 3D cases (d)–(f) compliance minimization benchmarks, showing the evolution of the design loss (left) and the volume constraint error (right). The loss for cases (c) and (f) are plotted against the right-hand $y$–axis.} 
    \label{fig: compliance_2d3d_convergence_curve}
\end{figure}

The numerical behavior of the APF-FNN framework for all six compliance benchmarks is illustrated in Figs.~\ref{fig: 2d_compliance_design}--\ref{fig: 3d_compliance_design}. In the 2D tests, the phase field functions evolve from a spatially uniform initialization into well-organized truss-like layouts with almost binary values and sharp yet smooth interfaces, cf. Fig.~\ref{fig: 2d_compliance_design}. The final designs closely reproduce classical reference solutions, including the characteristic tension/compression paths and load-transfer symmetries, while avoiding checkerboarding or pixel-level artifacts despite the relatively compact FNN parameterization. The 3D extensions in Fig.~\ref{fig: 3d_compliance_design} exhibit the same qualitative behavior: starting from a uniform distribution, the optimizer progressively forms spatially intricate, load-bearing frame-like structures with well-resolved interfaces.

The convergence histories in Fig.~\ref{fig: compliance_2d3d_convergence_curve} indicate that these topologies are reached through a stable alternating optimization process. Across all cases, the design loss $\mathcal{L}_{\text{topology}}$ typically undergoes a short transient in the early iterations, where mild non-monotonicity (peaks in 2D and peak--valley patterns in 3D) can appear as the state solver catches up with rapidly changing designs and the volume-penalty parameter is updated. After this initial adjustment, the design loss decreases steadily followed by a gradual refinement phase, while the volume constraint violation is quickly reduced and subsequently maintained within a narrow band around zero. Minor oscillations in
the 3D runs are correlated with scheduled updates of the volume-penalty parameter, but do not affect the overall downward trend. Collectively, these results demonstrate that the APF-FNN scheme produces well-resolved, volume-constrained designs in both planar and fully three-dimensional settings under a unified choice of network architectures and training schedule.

\subsubsection{Comparison with the classical FEM phase-field solver}

The comparison is restricted to a classical FEM-based phase-field solver, since most existing neural-network-based topology optimization approaches do not provide publicly available implementations, making a fair and reproducible comparison difficult without re-implementing those methods. We compare APF-FNNs with a conventional FEM implementation of phase-field topology optimization based on the $L^2$ gradient flow described in Section~\ref{eq: H1gradientflow} (Fig.~\ref{fig: conventional_optimization}). Both methods use identical problem setups and parameters in Section~\ref{sec: linearsetup}. The initial phase field is set to a spatially uniform value $0.5$ for both approaches. For the FEM baseline, we implemented a classical phase-field topology optimization solver in \texttt{FEniCS} on a laptop (Intel(R) Core(TM) Ultra 5 225H @ 1.70\,GHz, 32\,GB RAM, 64-bit OS). The phase-field evolution equation is discretized following the semi-implicit scheme in \cite{qian2022phase}. For the 2D benchmarks, we use a $120\times 60$ mesh and run $100$ outer optimization iterations, with $20$ inner iterations for the phase-field update per outer iteration. For the 3D benchmarks, we use a $50\times 25\times 10$ mesh and run $200$ outer optimization iterations, with $5$ inner iterations for the phase-field update per outer iteration.

\begin{table}[t]
\centering
\caption{Quantitative comparison with the classical FEM-based phase-field method in cases (a)-(f).}
\label{tab:pf_fem_compare}
\begin{tabular}{c|cc|cc|cc}
\toprule
\multirow{2}{*}{Case} &
\multicolumn{2}{c|}{$\mathcal{J}_{obj}$} &
\multicolumn{2}{c|}{Time (s)} &
\multirow{2}{*}{$\Delta \mathcal{J}_{obj}$ (\%)} &
\multirow{2}{*}{Speedup} \\
\cmidrule(lr){2-3}\cmidrule(lr){4-5}
& APF-FNN & FEM-PF & APF-FNN & FEM-PF & & \\
\midrule
(a) & 0.5953 & 0.6930 & 204 & 202  & $-14.10$ & 0.99 \\
(b) & 0.6187 & 0.7235 & 211 & 202  & $-14.49$ & 0.96 \\
(c) & 0.1015 & 0.1486 & 174 & 203  & $-31.70$ & 1.17 \\
\midrule
(d) & 309.24 & 316.36 & 788 & 1461 &  $-2.25$ & 1.85 \\
(e) & 345.58 & 445.32 & 757 & 1681 & $-22.40$ & 2.22 \\
(f) &  67.78 &  86.50 & 1083 & 1527 & $-21.64$ & 1.41 \\
\bottomrule
\end{tabular}
\end{table}

\begin{figure}[t] 
\centering
\setlength{\tabcolsep}{3pt}
\renewcommand{\arraystretch}{1.0}
\newcommand{\imgw}{0.18\linewidth} 

\begin{tabular}{c|cc|c|cc}
\toprule
 & \multicolumn{2}{c|}{\text{2D}} & & \multicolumn{2}{c}{\text{3D}} \\
\cmidrule(lr){2-3}\cmidrule(lr){5-6}
 & \text{APF-FNN} & \text{FEM-PF} & & \text{APF-FNN} & \text{FEM-PF} \\
\midrule

(a) &
\includegraphics[width=\imgw]{ 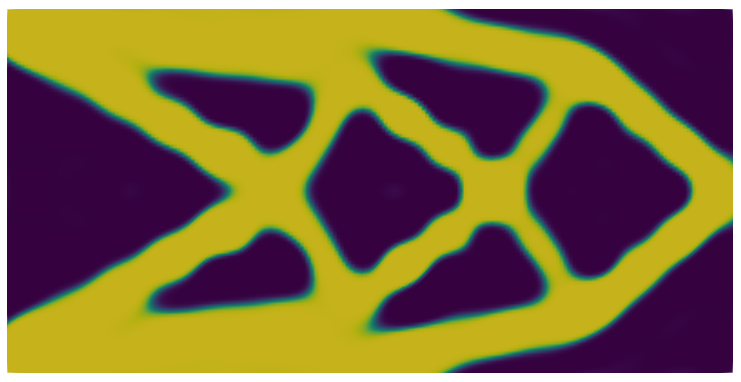} &
\includegraphics[width=\imgw]{ 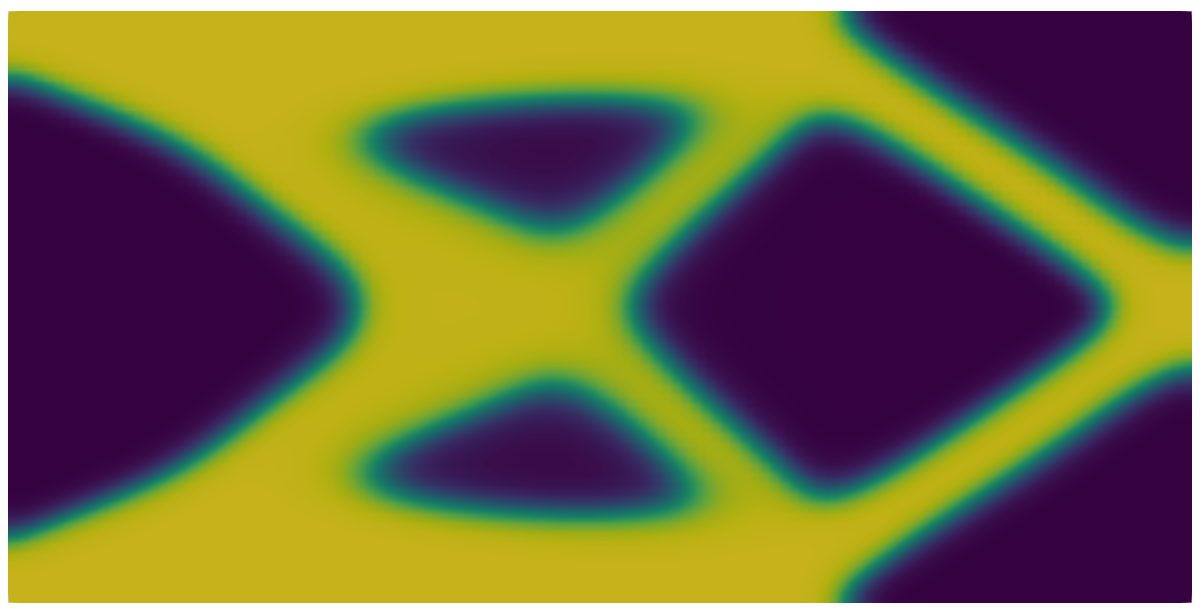} &
(d) &
\includegraphics[width=\imgw]{ 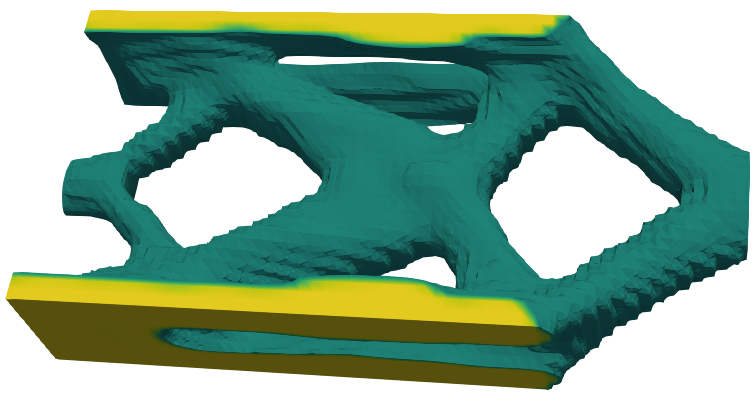} & 
\includegraphics[width=\imgw]{ 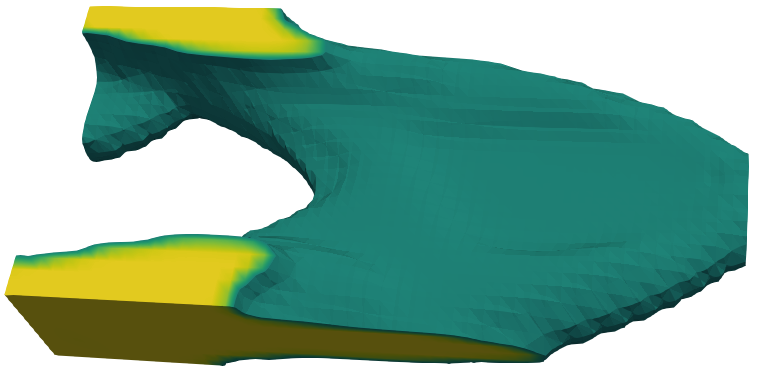} \\[3pt]

(b) &
\includegraphics[width=\imgw]{ 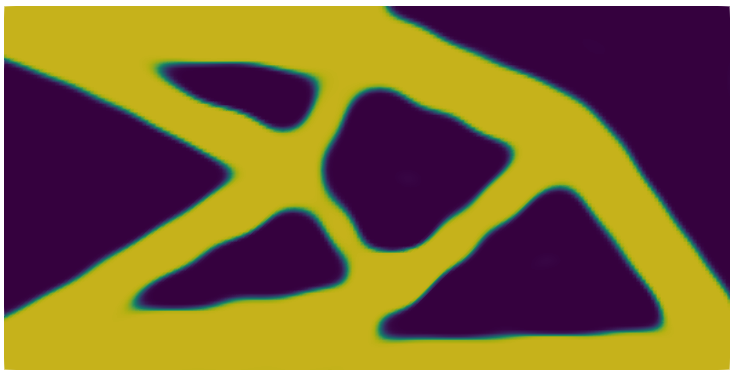} &
\includegraphics[width=\imgw]{ 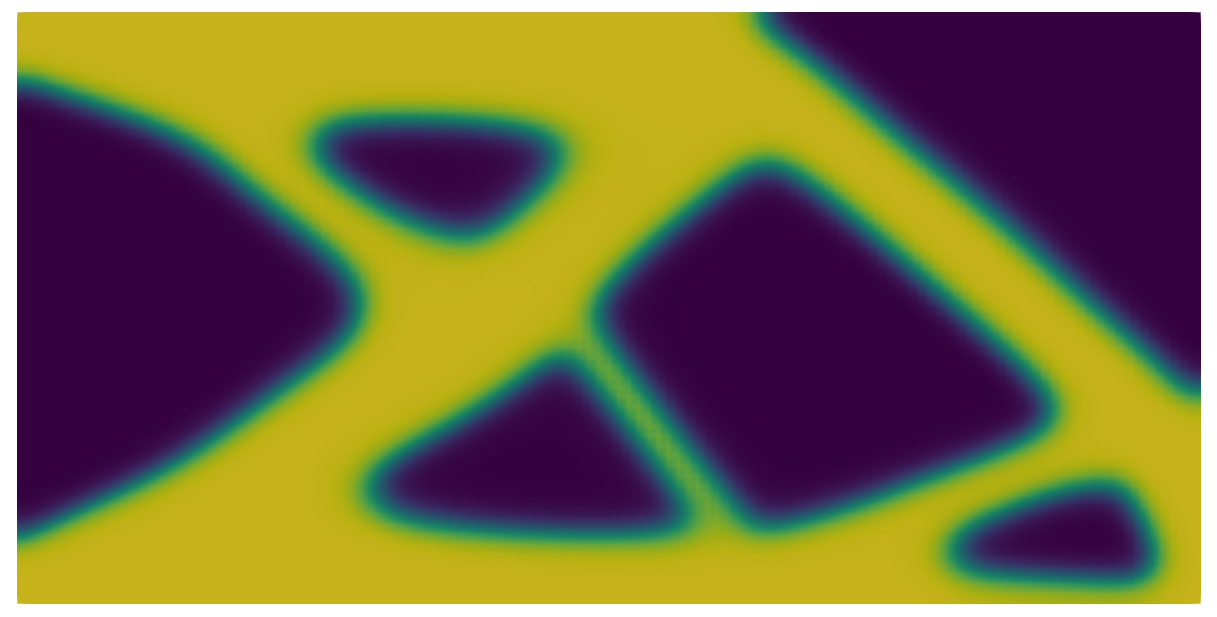} &
(e) &
\includegraphics[width=\imgw]{ 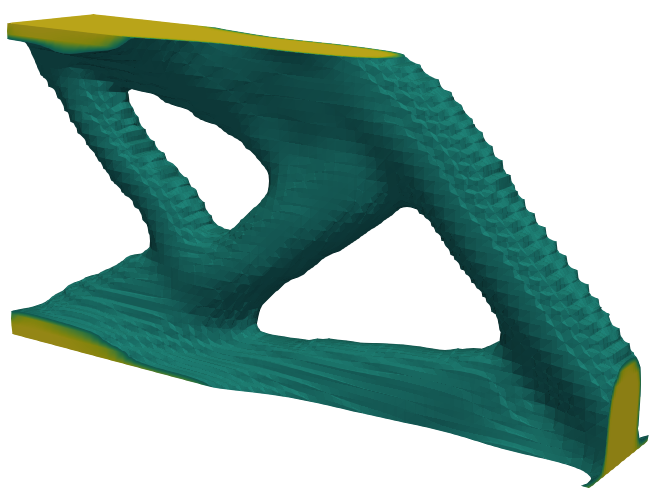} &   
\includegraphics[width=\imgw]{ 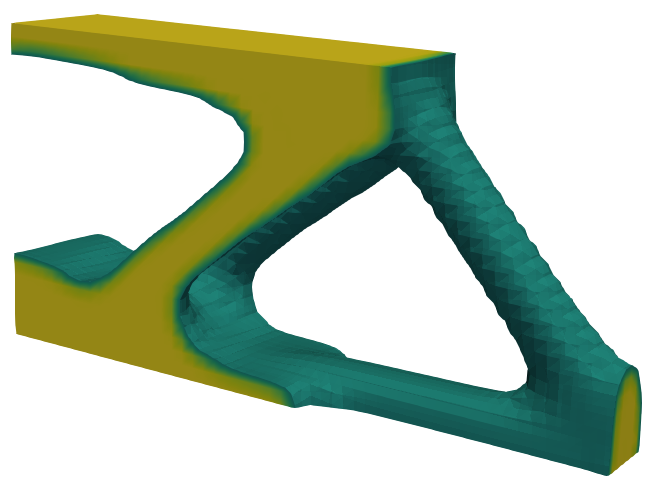} \\[3pt]

(c) &
\includegraphics[width=\imgw]{ 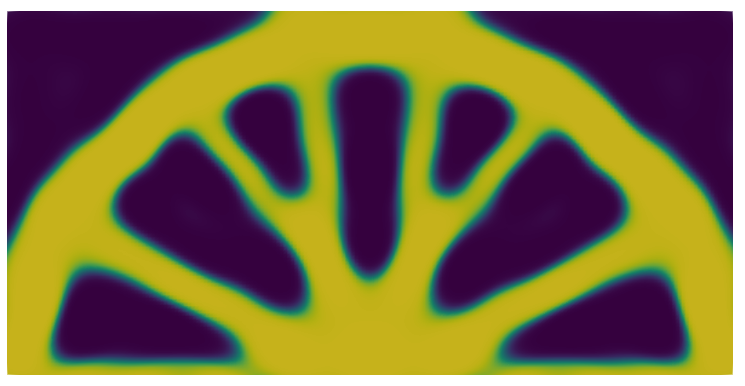} &
\includegraphics[width=\imgw]{ 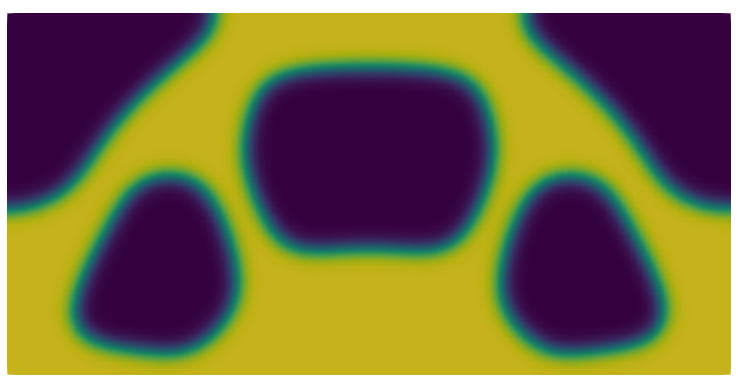} &
(f) &
\includegraphics[width=\imgw]{ 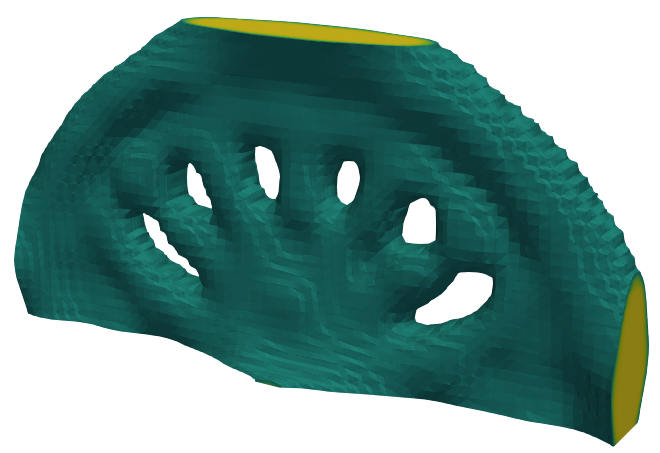} &   
\includegraphics[width=\imgw]{ 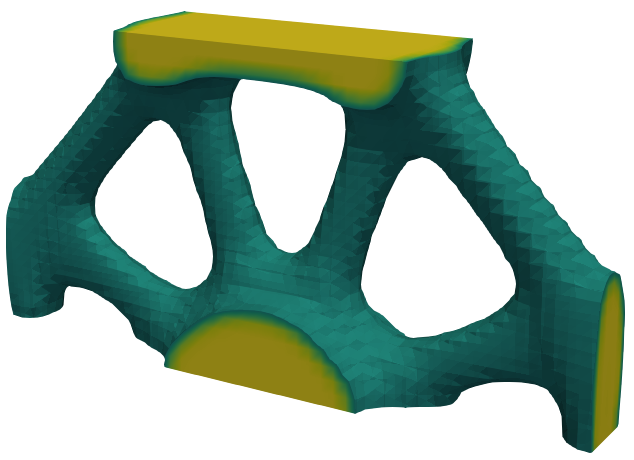} \\

\bottomrule
\end{tabular}

\caption{Benchmark comparison: APF-FNNs vs. FEM phase-field baseline for cases (a)-(f).}
\label{fig:pf_fem_compare}
\end{figure}

Fig.~\ref{fig:pf_fem_compare} provides a side-by-side visualization of the final optimized designs obtained by APF-FNNs and the classical FEM-based phase-field (FEM-PF) baseline. The left block reports the three 2D benchmarks (cases (a)--(c)), where each row corresponds to a benchmark and the two columns compare APF-FNN and FEM-PF. The right block presents the corresponding 3D extensions (cases (d)--(f)) in the same layout. Overall, both approaches recover physically meaningful load-carrying layouts, while APF-FNNs tend to yield clearer structural members and more distinct interfaces under the same uniform initialization.

Table~\ref{tab:pf_fem_compare} summarizes the quantitative comparison in terms of the compliance objective $\mathcal{J}_{\mathrm{obj}}$ and wall-clock time, where $\Delta\mathcal{J}_{\mathrm{obj}}$ denotes the relative change with respect to the FEM-PF baseline and the speedup is the ratio of wall-clock times. APF-FNNs achieve lower compliance in all six benchmarks, with relative reductions ranging from $14.10\%$ to $31.70\%$ in 2D and from $2.25\%$ to $22.40\%$ in 3D. The lower compliance achieved by APF-FNNs stems from the elimination of the $L^2$ gradient-flow approximation. While the FEM-PF approaches optimality through an artificial time evolution, APF-FNNs directly minimize the objective, resulting in designs that are closer to first-order optimality. In terms of runtime, the two methods are comparable for the 2D problems (speedup $\approx 0.96$--$1.17$), whereas in 3D the FEM-PF baseline becomes substantially more expensive and APF-FNNs provide consistent speedups ($\approx 1.41$--$2.22$) while still attaining better objective values. This runtime trend is expected since the FEM-PF baseline is a high-fidelity nested gradient-flow scheme that repeatedly solves the state/adjoint systems and updates the phase field on the full 3D mesh, making its cost grow rapidly in 3D. By contrast, APF-FNNs use neural surrogates for the state/adjoint and warm-start the networks via transfer learning, thereby avoiding repeatedly solving the full PDE systems from scratch at each iteration.

\subsubsection{Scalability analysis}
\label{subsec:scalability}

We analyze scalability on the representative benchmark case (a). Specifically, we investigate how the optimization behavior varies with the phase-field regularization parameters, the number of collocation points, the use of displacement-network pretraining, the inner-loop budgets for the state and topology networks in the alternating optimization, and the spectral bandwidth of the Fourier-feature representations.

\begin{figure}[t]
  \centering

  \begin{subfigure}[t]{0.24\linewidth}
    \centering
    \includegraphics[width=\linewidth]{ 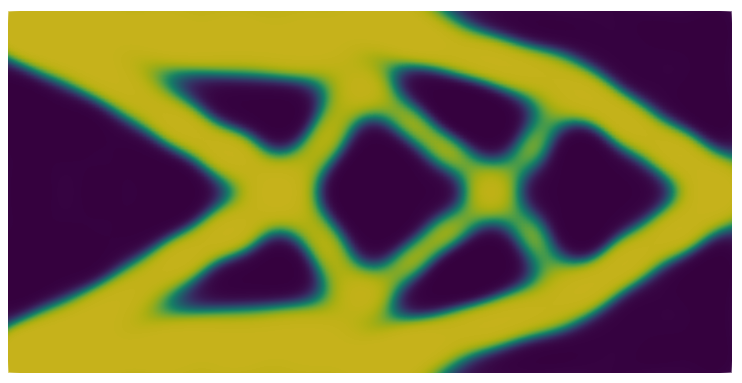}
    \caption*{(i) $\gamma=10^{-2},\,\epsilon=10^{-1}$}
    \label{fig:pf_g1e-2_e1e-1}
  \end{subfigure}
  \begin{subfigure}[t]{0.24\linewidth}
    \centering
    \includegraphics[width=\linewidth]{ 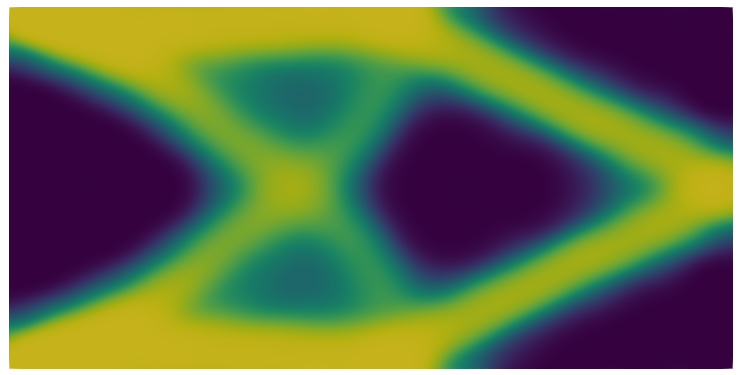}
    \caption*{(ii) $\gamma=10^{-2},\,\epsilon=1$}
    \label{fig:pf_g1e-2_e1e0}
  \end{subfigure}
  \begin{subfigure}[t]{0.24\linewidth}
    \centering
    \includegraphics[width=\linewidth]{ 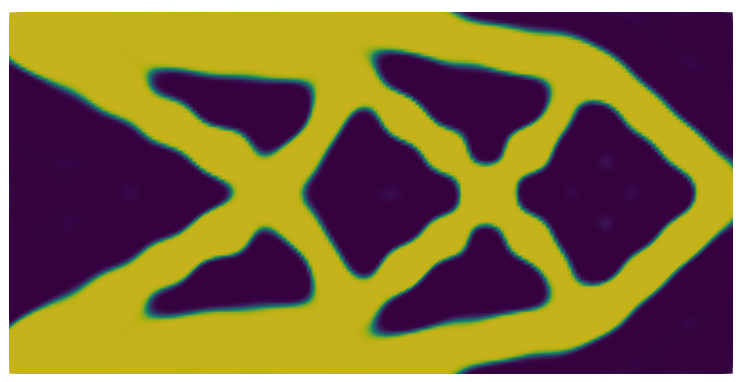}
    \caption*{(iii) $\gamma=10^{-3},\,\epsilon=10^{-2}$}
    \label{fig:pf_g1e-3_e1e-2}
  \end{subfigure}
  \begin{subfigure}[t]{0.24\linewidth}
    \centering
    \includegraphics[width=\linewidth]{ 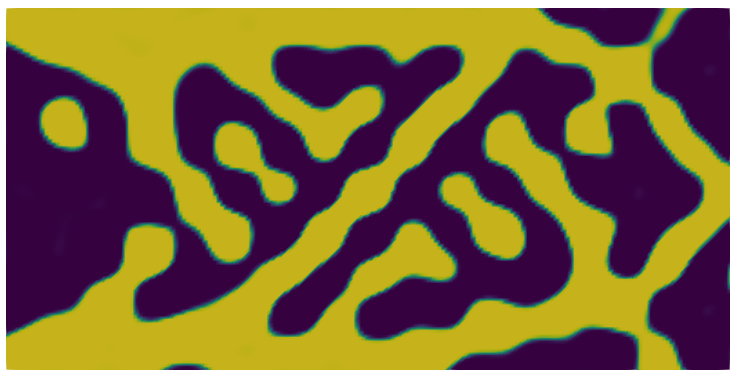}
    \caption*{(iv) $\gamma=10^{-4},\,\epsilon=10^{-2}$}
    \label{fig:pf_g1e-4_e1e-2}
  \end{subfigure}

  \caption{Influence of phase-field parameters $(\gamma,\epsilon)$ on the optimized topology fields in case (a).}
  \label{fig:pf_sensitivity}
\end{figure}

\paragraph{Influence of phase-field parameters.}
To examine the influence of the phase-field parameters on the optimized topology, we perform a controlled sensitivity study on $\gamma$ and $\epsilon$ while keeping all other settings identical. As illustrated in Fig.~\ref{fig:pf_sensitivity}, the optimized topologies exhibit clear and systematic changes as $\gamma$ and $\epsilon$ vary. With $\gamma=10^{-2}$ fixed, increasing $\epsilon$ from $10^{-2}$ (Fig.~\ref{fig: 2d_compliance_design}) to $10^{-1}$ and $1$ (Fig.~\ref{fig:pf_sensitivity}(i)-(ii)) clearly thickens the solid--void transition: boundaries become more diffuse, members thicken, and fine truss-like branches are gradually smoothed out; at $\epsilon=1$ the design is overly coarse with smeared interfaces. With $\epsilon=10^{-2}$ fixed, decreasing $\gamma$ weakens regularization: $\gamma=10^{-3}$ still yields a clear truss-like layout (Fig.~\ref{fig:pf_sensitivity}(iii)), whereas $\gamma=10^{-4}$ produces fragmented, noisy microstructures with many small branches (Fig.~\ref{fig:pf_sensitivity}(iv)). Overall, $\epsilon$ controls interface thickness (minimum length scale) and $\gamma$ controls the suppression of spurious fine-scale artifacts, confirming that both parameters take effect as predicted by phase-field model in our implementation.

\paragraph{Scalability w.r.t. the number of sampling points.}
To study the effect of sampling density, we vary the numbers of interior collocation points in the $x$- and $y$-directions and the number of sampling points on the non-homogeneous Neumann boundary, denoted by $(n_x,n_y,n_\Gamma)$. As the sampling density increases from $(50,25,125)$ to $(100,50,250)$ and $(200,100,500)$, the final compliance decreases from $0.6915$ to $0.6350$ and $0.5953$, indicating that a denser sampling provides a more accurate approximation of the domain and boundary integrals and yields a more reliable sensitivity for topology updates. This accuracy gain comes at an increased computational cost: the corresponding wall-clock times are $36\,\mathrm{s}$, $67\,\mathrm{s}$, and $204\,\mathrm{s}$, respectively. When the sampling density is further increased to $(400,200,1000)$, the compliance does not decrease further under the same training budget; instead, the late-stage compliance reaches $0.6054$, which is comparable to, or slightly higher than, that obtained with $(200,100,500)$, while the runtime increases substantially to $746\,\mathrm{s}$. This behavior suggests that simply increasing sampling points does not automatically translate into a better final design unless the optimization budget (e.g., number of training iterations and/or learning-rate schedule) is adjusted accordingly. Given the fixed settings in this study, $(200,100,500)$ provides the most favorable trade-off between solution quality and computational cost.

% \begin{figure}[t]
%   \centering
%   \includegraphics[width=0.4\linewidth]{ sampling_data_compliance3zoomin.png}
%   \caption{Compliance histories for Case (a) with different sampling settings $(n_x,n_y,n_\Gamma)$, where $n_x$ and $n_y$ are the numbers of interior collocation points in the $x$- and $y$-directions, and $n_\Gamma$ is the number of sampling points on the non-homogeneous Neumann boundary.}
%   \label{fig:casea_sampling_compliance}
% \end{figure}

\begin{figure}[htbp]
  \centering

  \begin{minipage}[t]{0.62\textwidth}
    \vspace{0pt}
    \centering
    \setlength{\tabcolsep}{3pt} 
    \renewcommand{\arraystretch}{1.0}
    \begin{tabular}[t]{ccc}
      \includegraphics[width=0.32\linewidth]{ 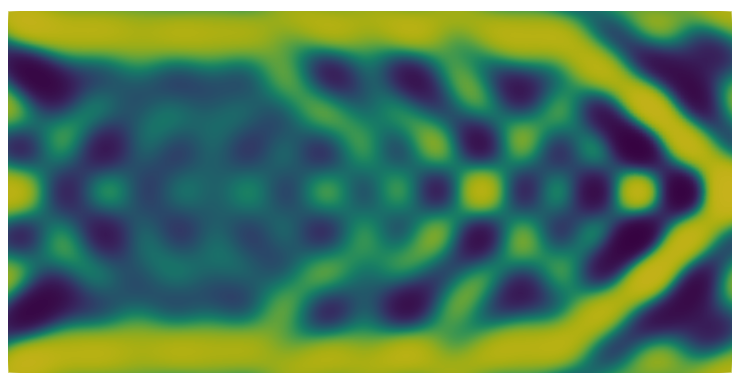} &
      \includegraphics[width=0.32\linewidth]{ 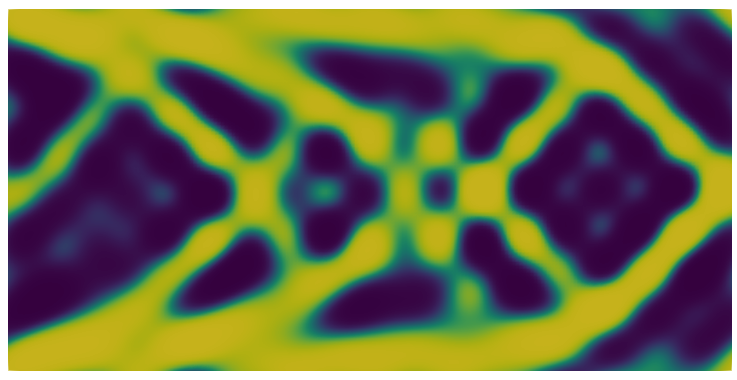} &
      \includegraphics[width=0.32\linewidth]{ 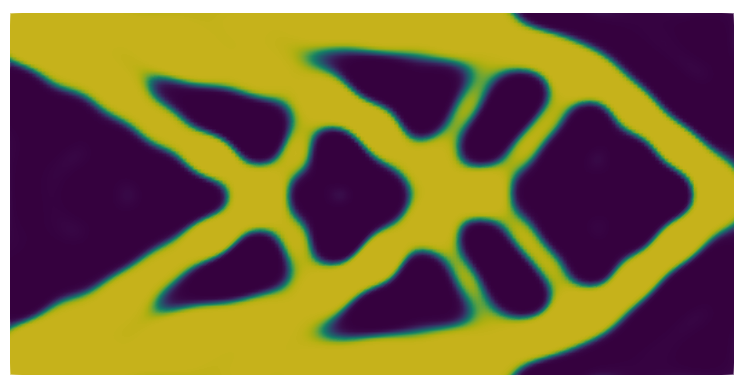} \\
      \small (i) & \small (ii) & \small (iii) \\[0.8ex]
      \includegraphics[width=0.32\linewidth]{ 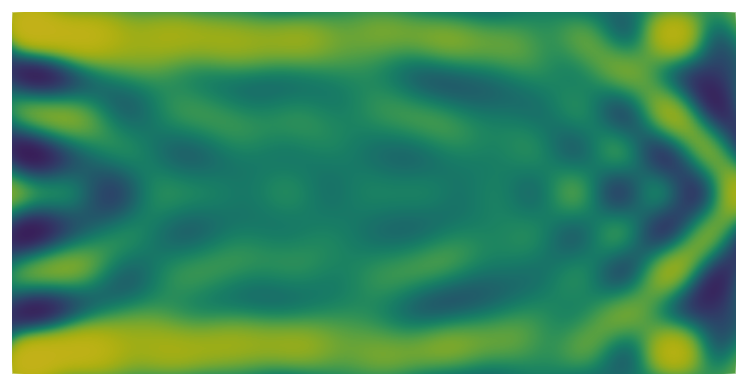} &
      \includegraphics[width=0.32\linewidth]{ 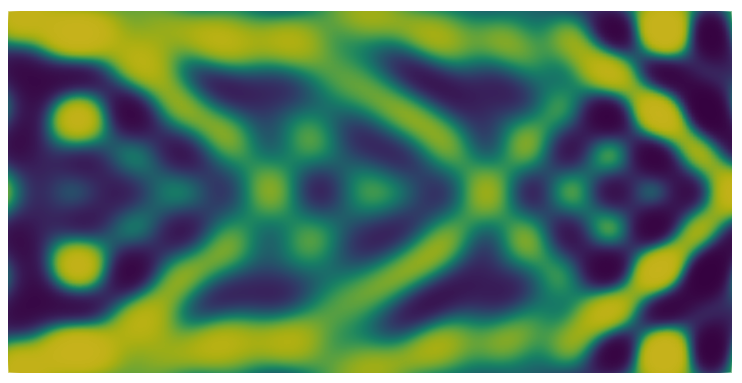} &
      \includegraphics[width=0.32\linewidth]{ 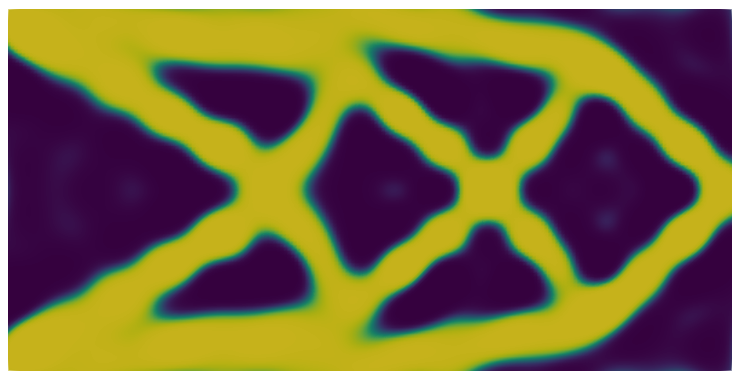} \\
      \small (iv) & \small (v) & \small (vi)
    \end{tabular}
  \end{minipage}\hfill
  \begin{minipage}[t]{0.36\textwidth}
    \vspace{0pt}
    \centering
    \includegraphics[width=\linewidth]{ 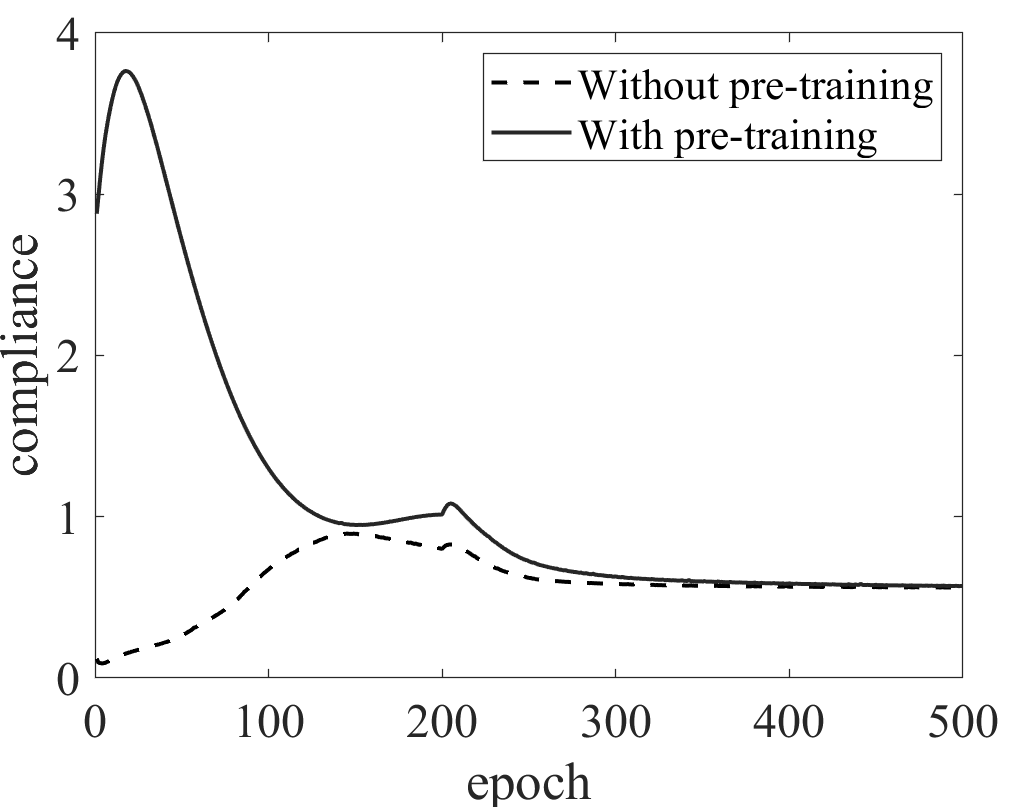} %\\[0.5ex]
    %\small (g)
  \end{minipage}

  \caption{
  Effect of state-network pretraining on the optimization dynamics in case (a). (i)--(iii) topology fields without pretraining at epochs 100, 200, and 300; (iv)--(vi) topology fields with pretraining at epochs 100, 200, and 300. The right panel shows the corresponding compliance convergence histories.}
  \label{fig:pretrain_scalability}
\end{figure}

\paragraph{Effect of state-network pretraining.}
We evaluate the scalability of APF-FNNs by examining how state-network pretraining affects the optimization dynamics. As shown in Fig.~\ref{fig:pretrain_scalability}, pretraining does not alter the final optimized topology or compliance, but it markedly improves early-stage stability and accelerates the emergence of coherent load-carrying structures. Without pretraining, the state approximation is inaccurate in the initial iterations, which contaminates design sensitivities and leads to oscillatory topology updates and delayed structural organization; in contrast, pretraining provides a physically consistent initialization of the state field, yielding reliable design gradients from the outset. This benefit is increasingly important as the problem scale grows: for high-resolution settings with substantially more design degrees of freedom, reaching an effective optimization regime within a limited computational budget becomes critical. Overall, these results indicate that state-network pretraining is an enabling component for the practical scalability of APF-FNNs.

\begin{figure*}[t]
  \centering

  % ---- Left: Figure ----
  \begin{minipage}[t]{0.6\textwidth}
    \vspace{0pt}
    \centering
    \includegraphics[width=0.6\linewidth]{ 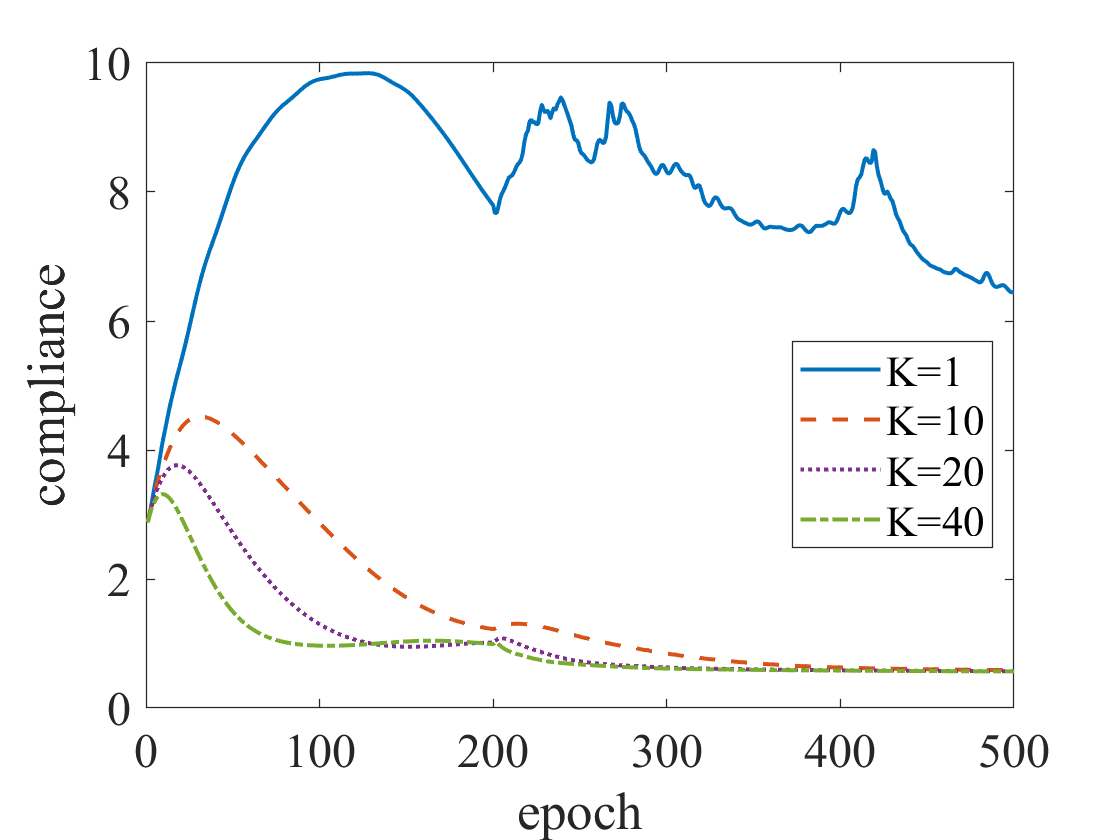}

    \captionsetup{justification=centering,singlelinecheck=false}
    \captionof{figure}{Compliance convergence histories for different numbers of state-network inner iterations $K$ in case (a).}
    \label{fig:state_inner_curve}
  \end{minipage}\hfill
  % ---- Right: Table ----
  \begin{minipage}[t]{0.35\textwidth}
    \vspace{0pt}
    \centering

    \captionsetup{justification=centering,singlelinecheck=false}
    \captionof{table}{Final compliance and runtime for different $K$ in case (a).}
    \label{tab:state_inner_cost}

    \vspace{0.5ex}

    \begin{tabular}{c c c}
      \toprule
      $K$ & Compliance & Time (s) \\
      \midrule
      1  & 0.7511 & 51  \\
      10 & 0.6165 & 124 \\
      20 & 0.5953 & 204 \\
      40 & 0.6021 & 361 \\
      \bottomrule
    \end{tabular}
  \end{minipage}
\end{figure*}

\paragraph{Effect of state-network inner iterations.}
We investigate scalability from the perspective of the optimization loop by varying the number of inner iterations $K$ used to update the state network, while keeping all other parameters fixed. Fig.~\ref{fig:state_inner_curve} shows the corresponding compliance convergence histories. When $K=1$, the optimization exhibits severe instability and fails to reach a low-compliance regime, indicating that insufficiently resolved state fields can fundamentally compromise the design update. Increasing $K$ markedly improves stability and convergence, with $K=20$ achieving the lowest final compliance. Further increasing $K$ yields diminishing returns, as the solution quality saturates while the computational cost increases substantially. Table~\ref{tab:state_inner_cost} summarizes the final compliance values and wall-clock times, highlighting a clear accuracy--efficiency trade-off. These results demonstrate that scalable performance requires sufficiently accurate state updates, but excessive inner iterations are unnecessary and inefficient.

\begin{figure}[htbp]
  \centering
  \begin{minipage}[t]{0.4\linewidth}
    \vspace{0pt}
    \centering
    \includegraphics[width=\linewidth]{ 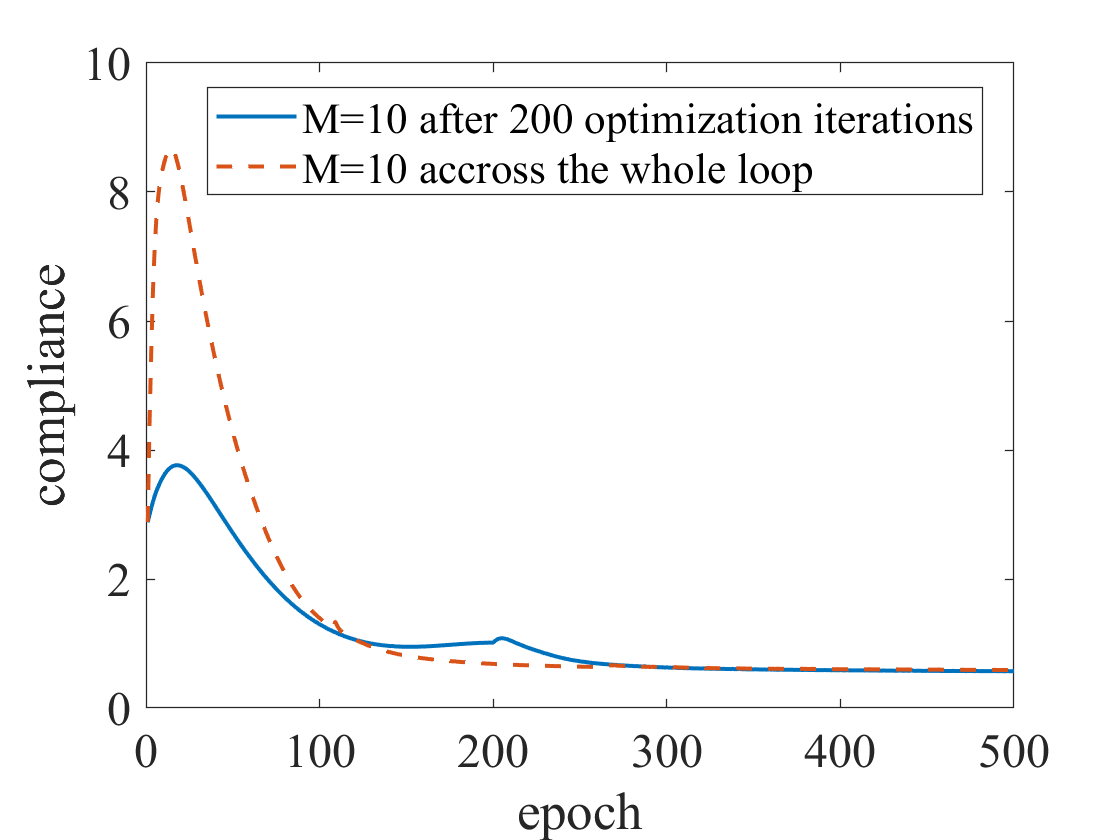}
    \small (i)
  \end{minipage}
  \begin{minipage}[t]{0.4\linewidth}
    \vspace{0pt}
    \centering
    \includegraphics[width=\linewidth]{ 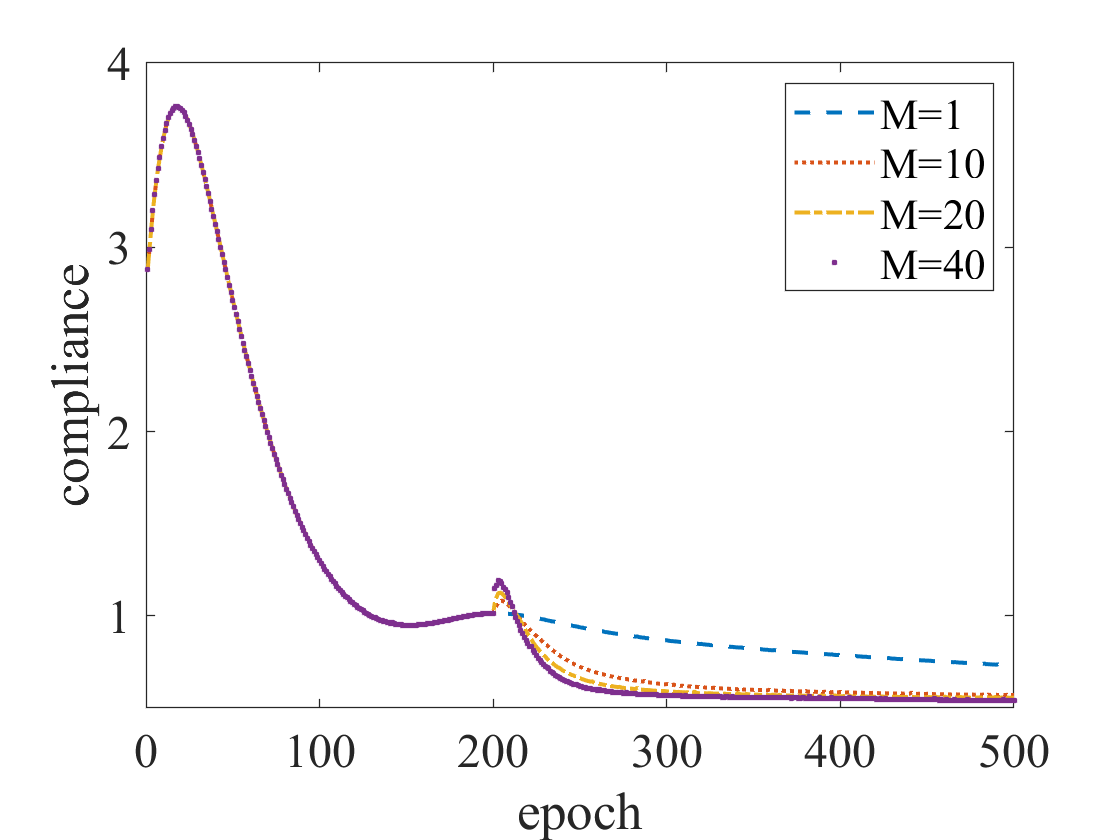}
    \small (ii)
  \end{minipage}

  \caption{
  Effect of scheduling phase-field inner iterations in case (a).
  (i) Topology evolution comparison between a fixed phase-field update budget ($M=10$ from the beginning) and a scheduled strategy (warm-up with $M=1$ for the first 200 epochs, then refinement with $M=10$).
  (ii) Compliance convergence histories with the warm-up stage (single phase-field update for the first 200 epochs) followed by different post-warm-up phase-field inner iterations $M\in\{1,10,20,40\}$.
  }
  \label{fig:phi_schedule_combined}
\end{figure}

\paragraph{Effect of phase-field inner iterations.}
Scalability is further examined from the optimization-loop perspective by scheduling the number of inner iterations of the topology network. Specifically, we perform a warm-up stage in which the phase-field network is updated only once per outer iteration for the first 200 epochs, and then increase the phase-field inner iterations to a larger value for refinement. Fig.~\ref{fig:phi_schedule_combined}(i) compares this scheduled strategy against using aggressive phase-field inner iterations from the beginning, showing that delaying intensive phase-field updates markedly stabilizes early optimization, suppresses fragmented intermediate patterns, and yields a cleaner load-carrying layout without sacrificing final performance. To determine an appropriate refinement budget after warm-up, we further vary the post-200-epoch inner iterations $M\in\{1,10,20,40\}$. As shown in Fig.~\ref{fig:phi_schedule_combined}(ii), $M=1$ results in substantially slower convergence and a much higher final compliance, indicating insufficient phase-field refinement, whereas $M\ge 10$ leads to rapid convergence and very similar final compliance values. The accuracy--efficiency trade-off is further reflected in the quantitative results: increasing $M$ from 1 to 10 yields a large improvement in compliance (0.7660$\rightarrow$0.5953) with only a modest increase in computational time (196$\rightarrow$204\,s), whereas further increasing $M$ beyond 10 leads to diminishing returns despite additional cost. These findings justify the proposed two-stage schedule (conservative warm-up followed by moderate refinement) and indicate that $M=10$--$20$ is sufficient for scalable performance.

% \begin{table}[htbp]
% \centering
% \caption{Final compliance and computational cost for different post-warm-up phase-field inner iterations $M$.}
% \label{tab:phi_Mcost}
% \begin{tabular}{c c c}
% \toprule
% $M$ & Final compliance & Time (s) \\
% \midrule
% 1  & 0.7660 & 196 \\
% 10 & 0.5953 & 204 \\
% 20 & 0.5928 & 210 \\
% 40 & \textbf{0.5894} & 228 \\
% \bottomrule
% \end{tabular}
% \end{table}

\begin{figure*}[t]
  \centering
  \setlength{\tabcolsep}{2pt}

  \begin{tabular}{ccccc}
    \includegraphics[width=0.19\textwidth]{ 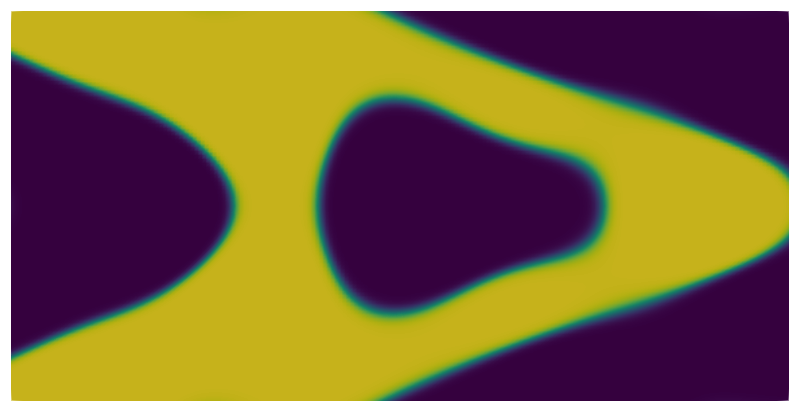} &
    \includegraphics[width=0.19\textwidth]{ 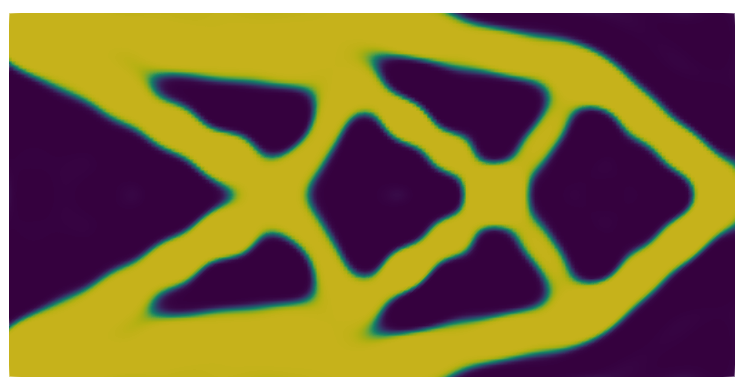} &
    \includegraphics[width=0.19\textwidth]{ 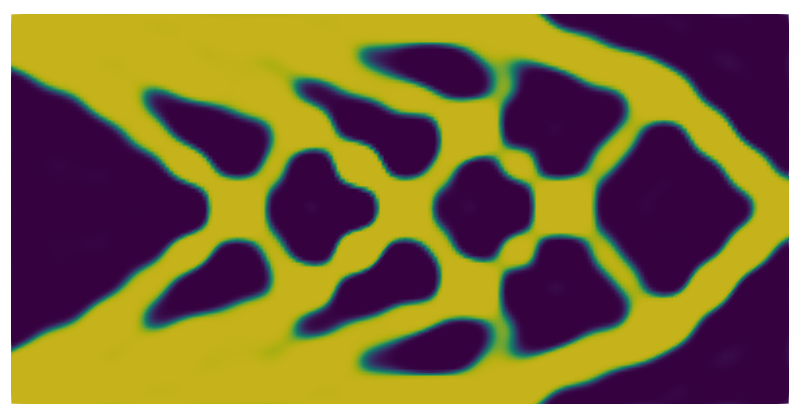} &
    \includegraphics[width=0.19\textwidth]{ 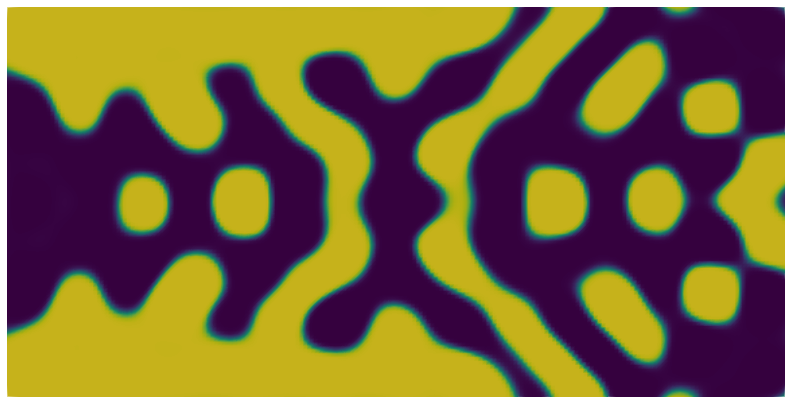} &
    \includegraphics[width=0.19\textwidth]{ 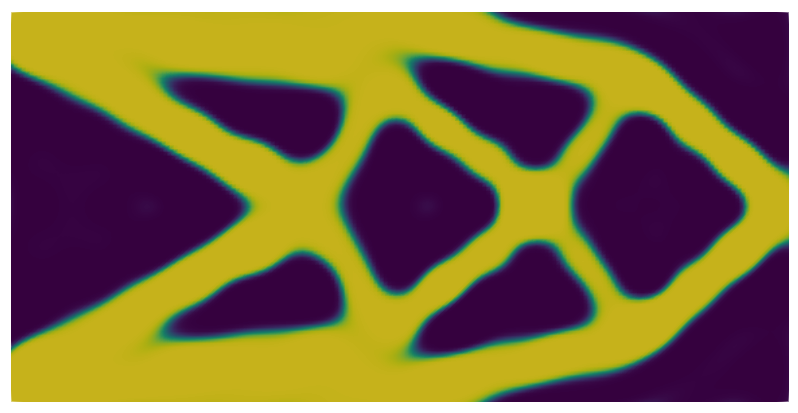} \\
     (i) $(10,35)$ &
     (ii) $(30,35)$ &
     (iii) $(31,35)$ &
     (iv) $(30,10)$ &
     (v) $(30,50)$
  \end{tabular}

  \caption{Effect of the Fourier-feature frequency on the optimized topologies. The tuple $(\omega_{\max}^{\phi},\,\omega_{\max}^{\bm{u}})$ under each subfigure denotes the frequency used in the phase-field and toplogy networks, respectively.}
  \label{fig:highband_tuple}
\end{figure*}

\paragraph{Scalability w.r.t. the spectral bandwidth.}
The Fourier-feature bandwidth $\omega_{\text{max}}$ in Section \ref{subsubsec: fnnsforall} sets the highest spatial frequencies representable by the neural fields, and thus acts as an effective spectral resolution that governs both design complexity and sensitivity quality. Fig.~\ref{fig:highband_tuple} sweeps bandwidth around the reference setting of case (a) ($\omega^{\phi}_{\text{max}}=35$; $\omega^{\bm u}_{\text{max}}=30$). Reducing the phase-field bandwidth yields overly smooth and simplified layouts, whereas increasing it enables richer fine-scale features (cf. Fig.~\ref{fig:highband_tuple}(i)-(iii)). In contrast, the displacement bandwidth primarily affects the fidelity of the state field and the resulting design gradients: a small bandwidth (e.g., 10) under-resolves the displacement solution and degrades structural organization, while a larger bandwidth (e.g., 50) produces crisper, straighter boundaries, indicating improved state representation (cf. Fig.~\ref{fig:highband_tuple}(iv), (ii), (v)). However, excessively large bandwidths may introduce high-frequency oscillations and alter the optimization trajectory, suggesting a practical saturation regime. The results suggest a moderate bandwidth range in which the phase-field bandwidth governs geometric richness while the displacement bandwidth provides accurate and stable sensitivities.

\subsection{Maximizing the first eigenvalue in linear elasticity}

\begin{figure}[htbp] 
    \centering
    \includegraphics[width=0.8\linewidth]{ 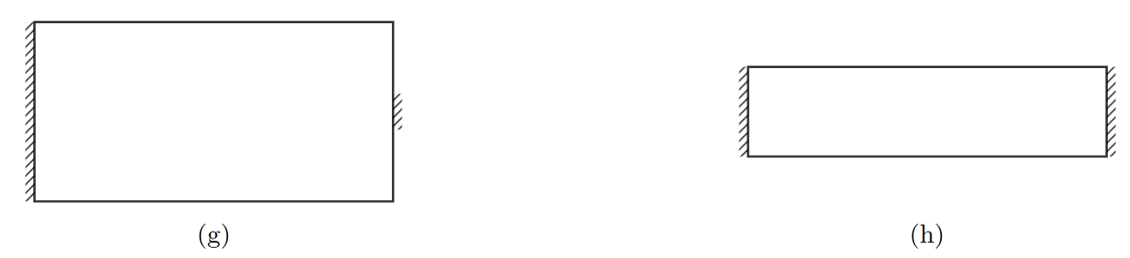} 
    \caption{Problem setups for the two 2D eigenvalue maximization benchmarks: (g) the partially constrained beam and (h) the fully clamped beam.} 
    \label{fig: maxeigen12dsettings}
\end{figure}

We next consider phase–field topology optimization for maximizing the fundamental eigenvalue of a linear–elastic structure, as formulated in Section~\ref{subsubsec:physics-eig}. We consider four benchmark geometries: two planar problems and their 3D extensions. The 2D problems, sketched in Fig.~\ref{fig: maxeigen12dsettings}, are

\begin{itemize}
    \item[(g)] \textit{Partially constrained beam:} $\Omega = (-1, 1) \times (-0.5, 0.5)$. The structure is clamped on the entire left boundary and on a central segment of the right boundary, with Dirichlet boundary $\Gamma_D = \{-1\} \times (-0.5, 0.5) \,\cup\, \{1\} \times (-0.05, 0.05)$. All remaining boundaries are traction-free.

    \item[(h)] \textit{Fully clamped beam:} $\Omega = (-1, 1) \times (-0.25, 0.25)$. The vertical boundaries at $x_{1} = \pm 1$ are fully clamped, while the horizontal boundaries at $x_{2} = \pm 0.25$ are traction-free.
    
    \item[(i)] \textit{3D partially constrained beam:} Extension of case~(g) to the domain $\Omega = (-1, 1) \times (-0.5, 0.5) \times (-0.2, 0.2)$. The structure is clamped on the entire left face and on a central segment of the right face, with $\Gamma_D = \{-1\} \times (-0.5, 0.5) \times (-0.2, 0.2)\,\cup\,\{1\} \times (-0.05, 0.05) \times (-0.2, 0.2)$. All other faces are traction-free.
    
    \item[(j)] \textit{3D fully clamped beam:} Analogous to case~(h), with domain $\Omega = (-1, 1) \times (-0.25, 0.25) \times (-0.3, 0.3)$. The vertical faces at $x_{1} = \pm 1$ are fully clamped, while all other faces are traction-free.
\end{itemize}

For the four eigenvalue benchmarks cases (g)–(j), we adopt a set of implementation choices consistent with the Fourier architectures described in Section~\ref{subsubsec: fnnsforall}. In all tests, the material model, SIMP interpolation and phase-field regularization follow the compliance setting; the target volume fraction is $\beta=0.5$ for the 2D cases (g)–(h) and $\beta=0.4$ for the 3D cases (i)–(j). In the 2D examples we adopt the same phase–field regularization parameters $(\varepsilon,\gamma)$ as in the minimum compliance problems, whereas in the 3D examples we set $(\varepsilon,\gamma)=(0.01,0.1)$. The Rayleigh-quotient formulation is augmented with a mass-normalization penalty of fixed weight $\lambda_{\mathrm{norm}}=50$. Homogeneous Dirichlet conditions on the clamped boundaries in cases (g) and (h) are enforced by a smooth multiplicative factor $g(\bm{x})$ in the eigenmode network, together with a soft boundary penalty of magnitude $\lambda_{\mathrm{D}}=200$ at the corresponding collocation points, while all remaining Dirichlet boundaries in all cases are imposed using boundary encodings based on $\tanh$. A summary of the main hyperparameters for cases (g)–(j) is provided in Table~\ref{tab:eigenvalue_hparams_full}. In all cases, the state frequencies are confined by $\omega^{\bm u}_{\max}=10$, and Fourier mode counts $(n^{\bm u}_1,n^{\bm u}_2)=(32,32)$ in 2D and $(n^{\bm u}_1,n^{\bm u}_2,n^{\bm u}_3)=(8,8,8)$ in 3D. The corresponding topology network is a single Fourier-feature phase-field model with $\omega^{\phi}_{\text{max}}=25$, using $(24,12)$ modes in 2D and $(12,6,6)$ modes in 3D. The eigenmode learning rate is $5\times10^{-2}$ during pre-training and $10^{-4}$ in the alternating loop for the 2D problems, and $5\times10^{-3}$ / $10^{-4}$ in 3D; the topology network is trained with a constant learning rate of $5\times10^{-4}$ in all cases.

\begin{table}[htbp]
\centering
\caption{Hyperparameters in Algorithm \ref{alg: tri-level to} for the elasticity eigenvalue benchmarks cases (g)--(j).}
\label{tab:eigenvalue_hparams_full}
\begin{tabular}{cccccccccc}
\toprule
Case & $N$ & $K$ & $|X|$ & $|X_b|$ & $(\lambda_{\mathrm{penal}}^{(0)},\zeta)$ & $\omega^{\bm u}_{\max}$ & $(n^{\bm u}_1,n^{\bm u}_2,n^{\bm u}_3)$ & $\omega^{\phi}_{\max}$ & $(n^{\phi}_1,n^{\phi}_2,n^{\phi}_3)$ \\
\midrule
(g) & 500  & 20  & 20000 & 1000 & $(1,0.98)$     & 10 & $(32,32,0)$ & 25 & $(24,12,0)$ \\
(h) & 500  & 20  & 20000 & 1000 & $(10,0.99)$    & 10 & $(32,32,0)$ & 25 & $(24,12,0)$ \\
(i) & 600 & 100 & 54000 & 5000 & $(100,1/1.01)$ & 10 & $(8,8,8)$  & 25 & $(12,6,6)$  \\
(j) & 600 & 100 & 32000 & 5000 & $(100,1/1.01)$ & 10 & $(8,8,8)$  & 25 & $(12,6,6)$  \\
\bottomrule
\end{tabular}
\end{table}

\begin{figure}[h!] 
    \centering
    \settowidth{\rowlabelwidth}{ (a))} 

    \begin{tabular}{
        >{}l 
        @{\hspace{1em}} 
        >{\centering\arraybackslash}m{0.22\textwidth} 
        @{\hspace{0.5em}} 
        % >{\centering\arraybackslash}m{0.22\textwidth}
        % @{\hspace{0.5em}}
        >{\centering\arraybackslash}m{0.22\textwidth}
        @{\hspace{0.5em}}
        >{\centering\arraybackslash}m{0.22\textwidth}
    }
        & \multicolumn{1}{c}{Epoch 100} 
        & \multicolumn{1}{c}{Epoch 200} 
        % & \multicolumn{1}{c}{Epoch 300} 
        & \multicolumn{1}{c}{Epoch 500} \\
        \addlinespace[5pt]
        
        % 第一行图片
        (g) &
        \includegraphics[width=\linewidth, height=3.5cm, keepaspectratio]{ 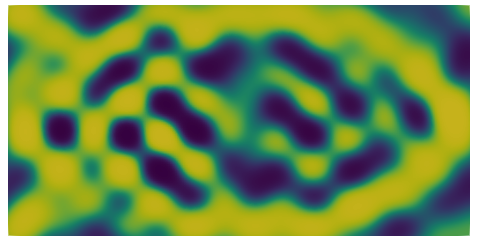} &
        \includegraphics[width=\linewidth, height=3.5cm, keepaspectratio]{ 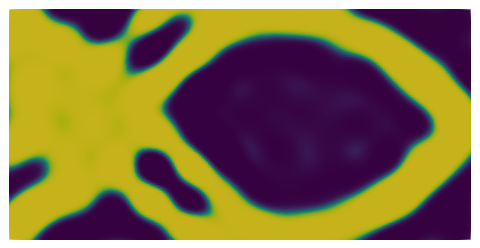} &
        \includegraphics[width=\linewidth, height=3.5cm, keepaspectratio]{ 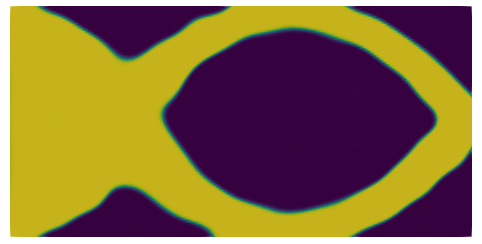} \\
        % \addlinespace[10pt]

        % 第二行图片
        (h) &
        \includegraphics[width=\linewidth, height=3.5cm, keepaspectratio]{ 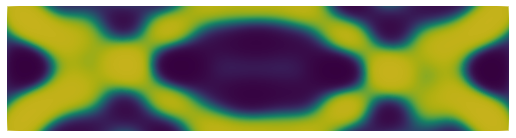} &
        \includegraphics[width=\linewidth, height=3.5cm, keepaspectratio]{ 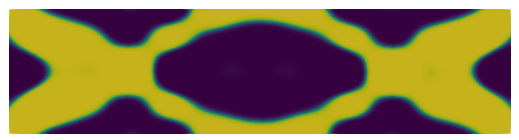} &
        \includegraphics[width=\linewidth, height=3.5cm, keepaspectratio]{ 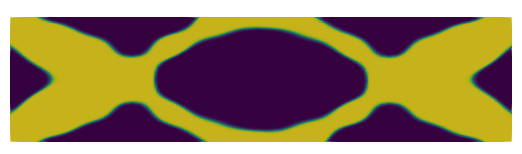} 
    \end{tabular}
    
   \caption{Evolution of the phase-field topology fields for the 2D eigenvalue maximization benchmarks.}
    \label{fig: 2d_eigen1_design}
\end{figure}

The eigenmode network is first pre-trained for $I$ Adam steps using only the state loss: in 2D we use an adaptive window $I\in[10^3,5\times10^3]$ with early stopping once the change in $\mathcal{L}_{\mathrm{state}}$ between successive iterations falls below $10^{-5}$, whereas in 3D we fix $I=5\times10^3$. The subsequent alternating optimization runs for $N=500$ outer iterations in the 2D cases and $N=600$ in 3D. Within each outer iteration, we perform $K$ eigenmode updates and $M$ topology updates. For the 2D examples, we take $K=20$ and start with $M=1$, increasing to $M=5$ once $n>100$; for the 3D cases we use a more aggressive $K=100$ and keep $M=1$ throughout. Interior collocation points are sampled from a tensor-product grid over the design domain, with $|X|=20{,}000$ interior points and $|X_b|=1{,}000$ boundary samples in 2D. In the 3D cases (i) and (j), we use $54{,}000$ interior points on a $60\times30\times30$ grid and $32{,}000$ interior points on an $80\times20\times20$ grid, respectively, together with $|X_b|=5{,}000$ boundary points in each case. The volume-penalty parameter in $\mathcal{P}_{\text{vol}}$ is ramped from an initial value of order $1$ (2D) or $10^2$ (3D) up to $10^3$–$10^4$ using a multiplicative factor $\zeta=0.98$ or $0.99$ in 2D and $\zeta=1/1.01$ in 3D, so that the volume constraint is gradually tightened as the design stabilizes.

\begin{figure}[htbp] 
    \centering
    \includegraphics[width=0.8\linewidth]{ 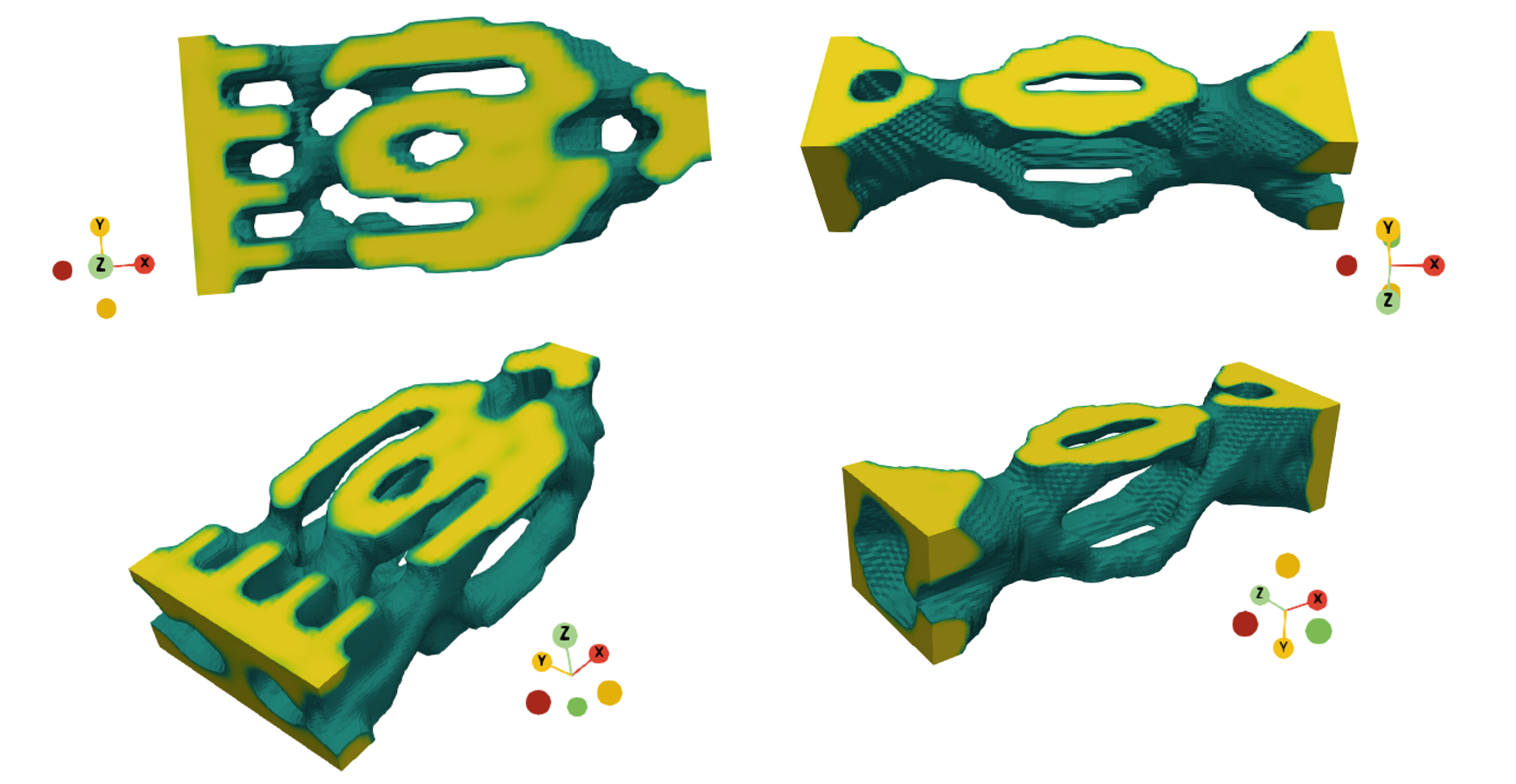}
     \caption{Optimized 3D phase-field topologies for the eigenvalue maximization benchmarks. Left: two views of the partially clamped case (i); right: two views of the fully clamped case (j).}
    \label{fig: 3d_eigen1_design}
\end{figure}

\begin{figure}[htbp] 
    \centering
    \includegraphics[width=0.35\linewidth]{ 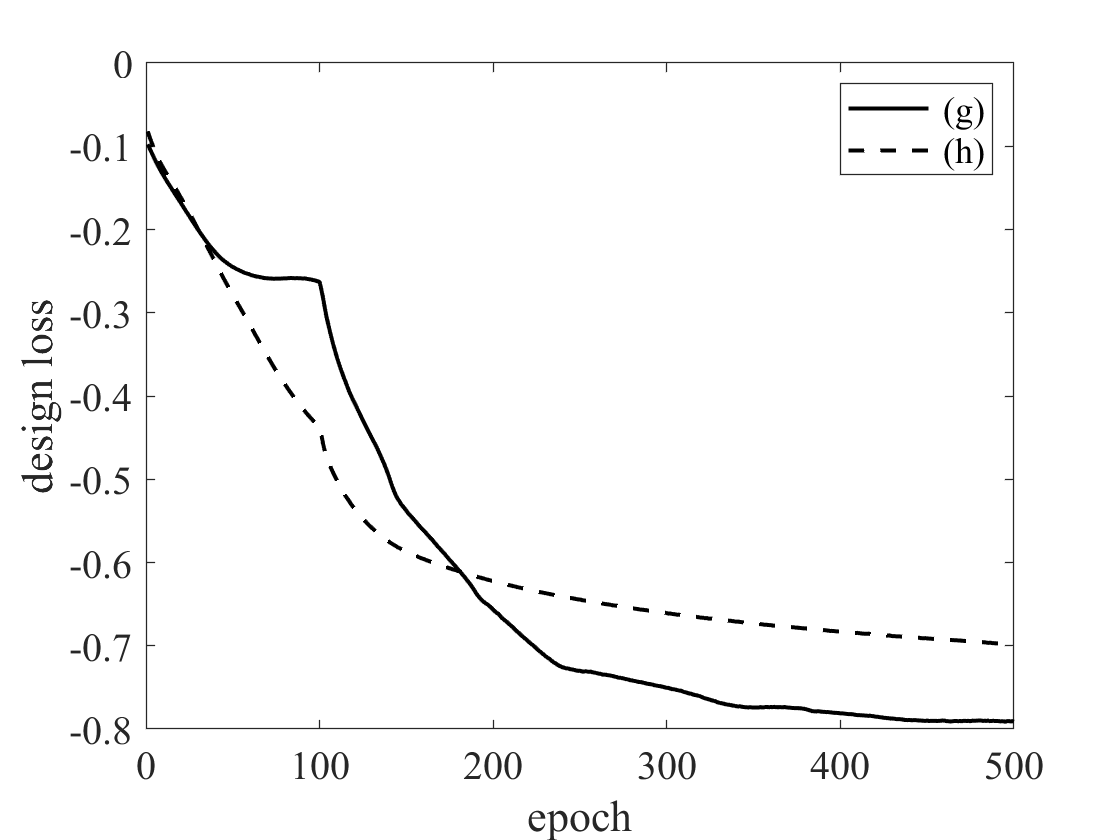}
    \includegraphics[width=0.35\linewidth]{ 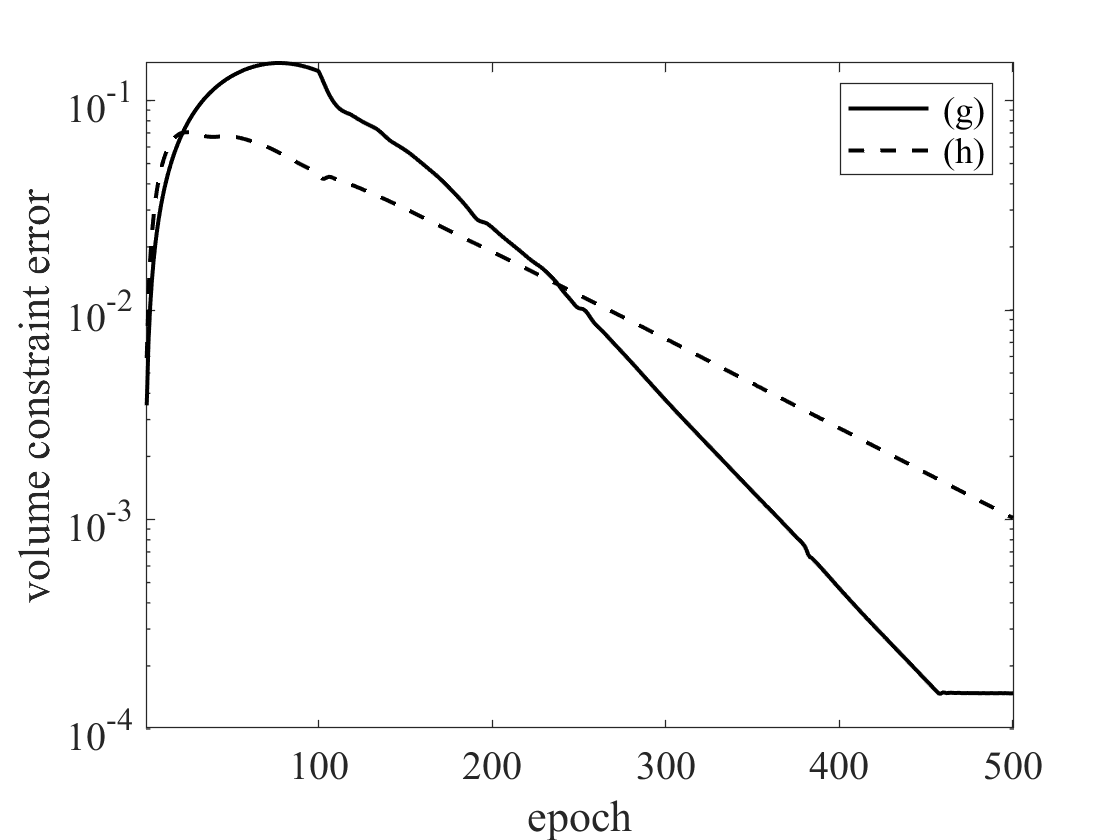}
    \\
    \includegraphics[width=0.35\linewidth]{ 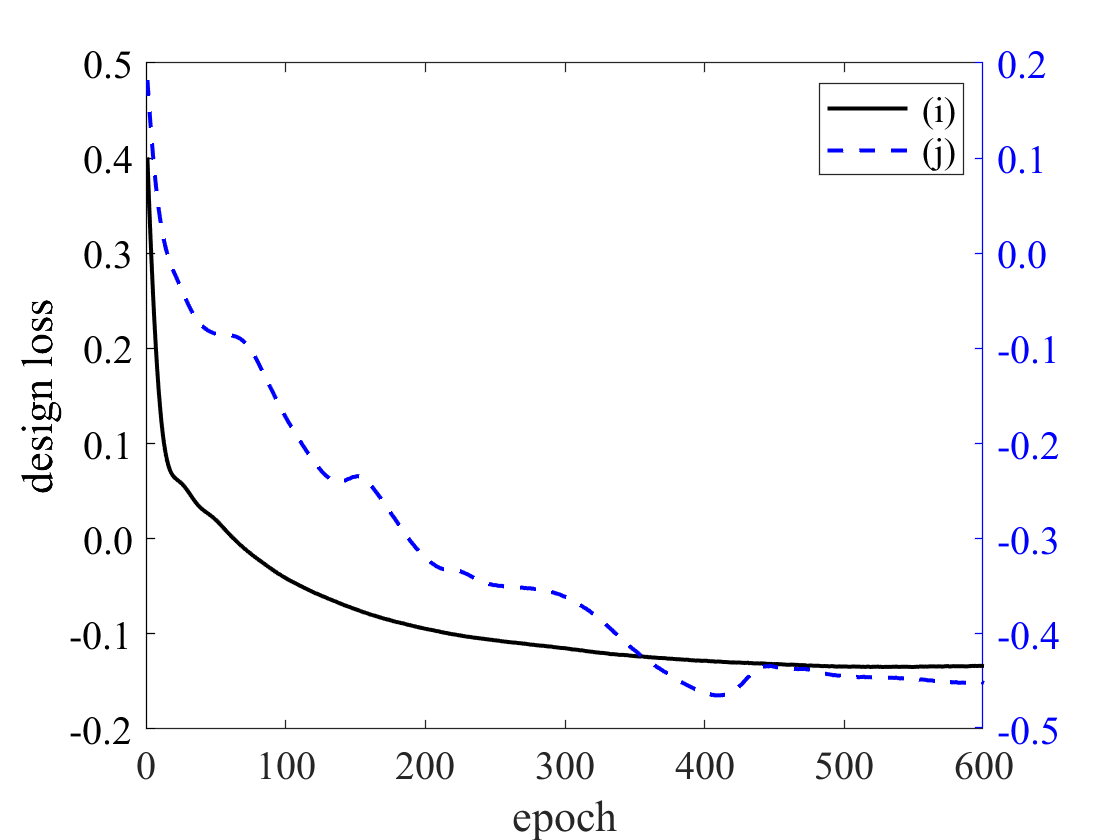}
    \includegraphics[width=0.35\linewidth]{ 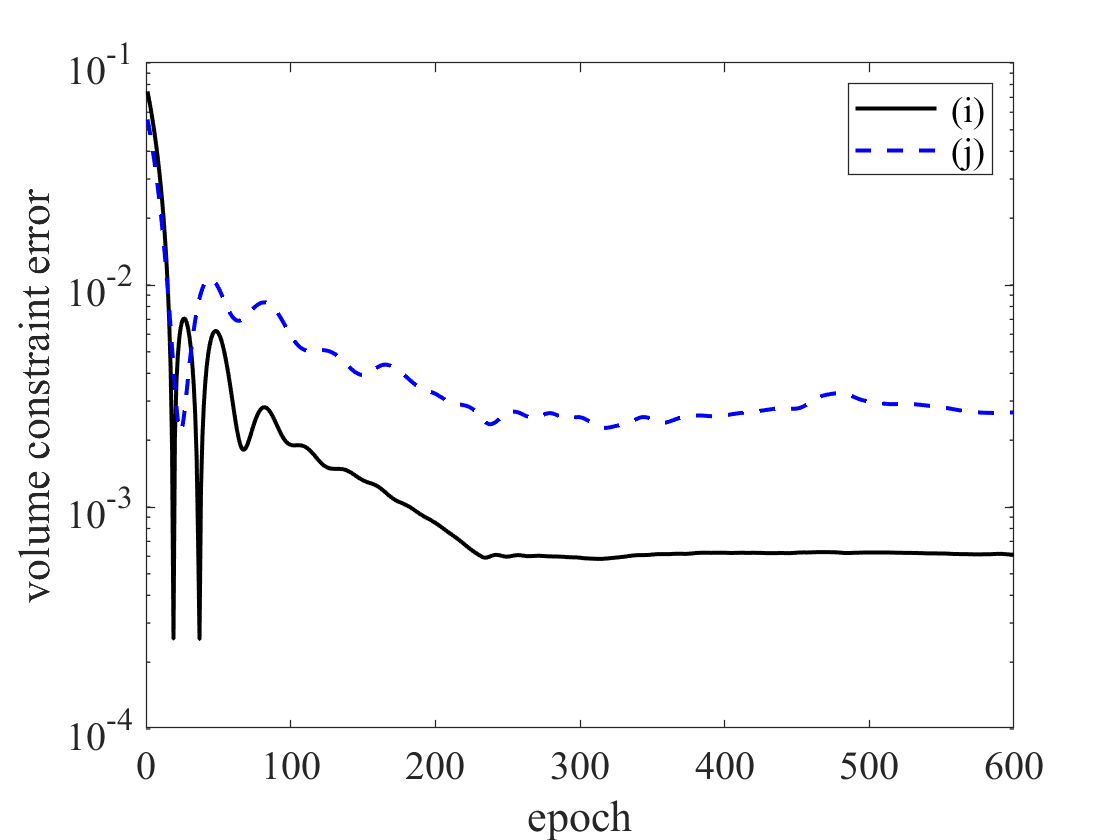}
    \caption{Convergence histories for the 2D (1st row) and 3D (2nd row) eigenvalue maximization benchmarks, showing the design loss (left) and volume constraint error (right) for the cases (g-j). The loss for case (j) is plotted against the right-hand $y$–axis.} 
    \label{fig: eigen1_2d3d_convergence_curve}
\end{figure}

The numerical results for the four eigenvalue benchmarks are summarized in Figs.~\ref{fig: 2d_eigen1_design}--\ref{fig: 3d_eigen1_design}. The qualitative evolution of the phase field is very similar to the compliance cases in Section~\ref{subsec: numerical compliance}: starting from a uniform distribution, the APF--FNNs quickly develop organized structural patterns and converge to almost binary topologies with smooth interfaces. In 2D, the partially clamped beam (g) concentrates material near the single support, whereas the fully clamped beam (h) yields a more symmetric layout attached to both ends. Their 3D counterparts (i)–(j) are direct spatial extensions, with shell– and frame-like configurations that mirror the 2D designs while respecting the same volume fraction.

The convergence histories in Fig.~\ref{fig: eigen1_2d3d_convergence_curve} show that the alternating scheme remains stable under the spectral objective. The design loss $\mathcal{L}_{\text{topology}}$ decreases steadily in all four cases, and the volume constraint error is reduced by several orders of magnitude and
then kept within a narrow band around zero. Compared with the compliance benchmarks, the eigenvalue problems exhibit slightly slower convergence and milder oscillations, reflecting the increased nonconvexity of the objective, but no loss of robustness.

\subsection{Stokes flow optimization}
\label{subsec:stokes_examples}

\begin{figure}[htbp] 
    \centering
    \includegraphics[width=0.8\linewidth]{ 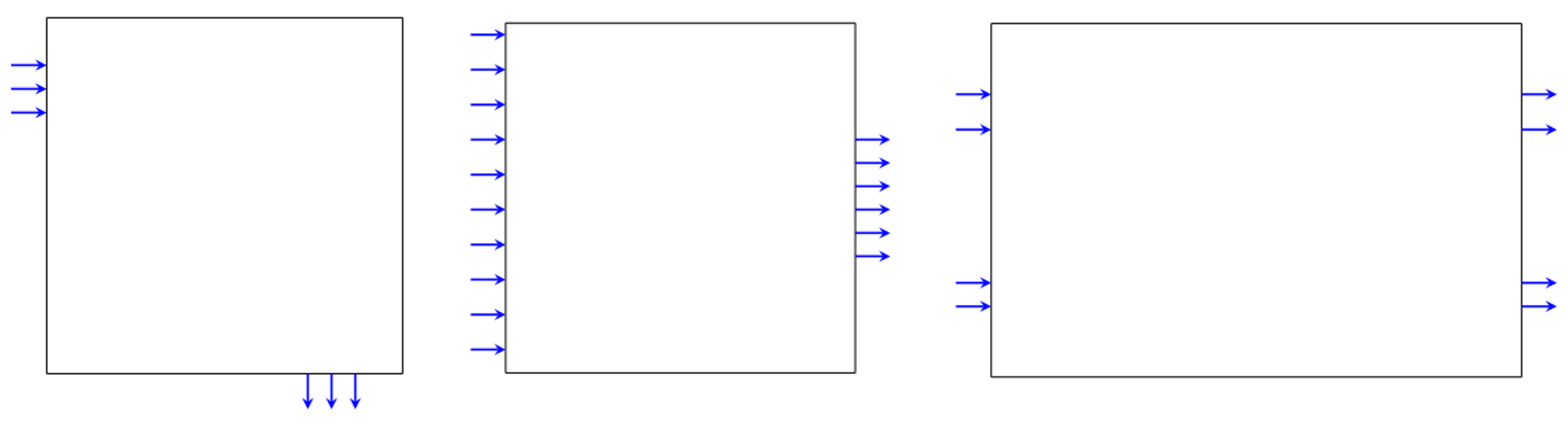} 
    \caption{Boundary conditions for the 2D Stokes benchmarks (k)–(m), shown from left to right. Inward arrows denote prescribed velocity inflows (non-homogeneous Dirichlet), outward arrows zero-traction outflows (homogeneous Neumann). All other boundaries are no-slip walls (homogeneous Dirichlet).}

    \label{fig: stokes2dsettings}
\end{figure}

We next consider phase–field topology optimization for steady incompressible Stokes flow, aiming to minimize the viscous power dissipation under a Darcy penalization of the solid phase, as formulated in Section~\ref{subsubsec:physics-stokes}. In all Stokes and Navier--Stokes benchmarks, the kinematic viscosity is fixed to $\rho = 0.01$, and the phase-field parameters in the Ginzburg--Landau regularization are set to $\varepsilon = 0.01$ and $\gamma = 0.01$. Five benchmarks are studied: three planar configurations and two fully three-dimensional extensions. The 2D boundary conditions are summarized in Fig.~\ref{fig: stokes2dsettings}; all unmarked boundaries in the following are treated as no–slip walls ($\bm{u}_\theta = \bm{0}$).

\begin{itemize}
  \item[(k)] \emph{Single inlet–outlet channel.}
  The domain is a unit square $\Omega = (0,1)^2$. A uniform horizontal inflow $\bm{u}_{\bm{\theta}} = (1,0)^{\top}$ is prescribed on the left boundary segment $\{x_1=0,\;0.7\le x_2\le 0.9\}$, while a homogeneous Neumann (zero-traction) outflow is imposed on the bottom segment $\{x_2=0,\;0.7\le x_1\le 0.9\}$.

  \item[(l)] \emph{Parabolic inlet channel.}
  Again $\Omega = (0,1)^2$. A parabolic profile $\bm{u}_{\bm{\theta}} = (4x_2(1-x_2),0)^{\top}$ is prescribed on the entire left edge $\{x_1=0\}$, and a zero-traction outflow condition is enforced on the right boundary segment $\{x_1=1,\;0.3\le x_2\le 0.7\}$.

  \item[(m)] \emph{Two–inlet mixing channel.}
  The domain is $\Omega=(0,1.5)\times(0,1)$. On both the left $\{x_1=0\}$ and right $\{x_1=1.5\}$ boundaries, two symmetric parabolic inflow profiles are imposed: on the lower ports $1/6\le x_2\le 1/3$, $\bm{u}_{\bm{\theta}} = (1-144(x_2-1/4)^2,0)^{\top}$, and on the upper ports $2/3\le x_2\le 5/6$, $\bm{u}_{\bm{\theta}} = (1-144(x_2-3/4)^2,0)^{\top}$.
\end{itemize}

The 3D tests are defined on variants of the unit cube and are designed to exercise the method on multi–inlet mixers and diffusers:
\begin{itemize}
  \item[(n)] \emph{Four–inlet mixer.}
  The domain is a unit cube $\Omega=(0,1)^3$. Four circular inlets are centered on the side faces $x_1=0$, $x_1=1$, $x_2=0$ and $x_2=1$, with prescribed velocities $\bm{u}_{\bm{\theta}}=(1,0,0)^{\top}$, $(-1,0,0)^{\top}$, $(0,1,0)^{\top}$ and $(0,-1,0)^{\top}$, respectively. A single circular outlet on the bottom face $\{x_3=0\}$ enforces a vertical velocity $\bm{u}_\theta=(0,0,-4)^{\top}$ so that the net flux balances the four inlets.

  \item[(o)] \emph{Diffuser.}
  On the same cubic domain, a uniform inflow $\bm{u}_{\bm{\theta}}=(0,1,0)^{\top}$ is prescribed on the entire bottom face $\{x_2=0\}$. A mass-conserving parabolic outflow profile is imposed on a circular outlet centered on the top face $\{x_2=1\}$, while all remaining faces are no–slip walls.
\end{itemize}

\begin{figure}[htbp!] 
    \centering
    \settowidth{\rowlabelwidth}{ (a))} 

    \begin{tabular}{
        >{}l 
        @{\hspace{1em}} 
        >{\centering\arraybackslash}m{0.22\textwidth} 
        @{\hspace{0.1em}}
        >{\centering\arraybackslash}m{0.22\textwidth}
        @{\hspace{0.1em}}
        >{\centering\arraybackslash}m{0.22\textwidth}
    }
        & \multicolumn{1}{c}{Epoch 100} 
        & \multicolumn{1}{c}{Epoch 200} 
        & \multicolumn{1}{c}{Final design} \\
        \addlinespace[5pt]
        
        % 第一行图片
        (k) &
        \includegraphics[width=\linewidth, height=3.5cm, keepaspectratio]{ 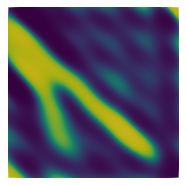} &
        \includegraphics[width=\linewidth, height=3.5cm, keepaspectratio]{ 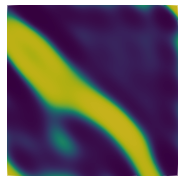} &
        \includegraphics[width=\linewidth, height=3.5cm, keepaspectratio]{ 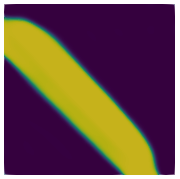} \\
        % \addlinespace[10pt]

        % 第二行图片
        (l) &
        \includegraphics[width=\linewidth, height=3.5cm, keepaspectratio]{ 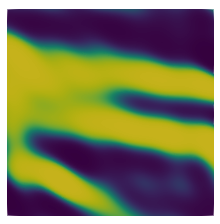} &
        \includegraphics[width=\linewidth, height=3.5cm, keepaspectratio]{ 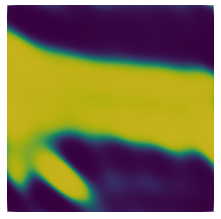} &
        \includegraphics[width=\linewidth, height=3.5cm, keepaspectratio]{ 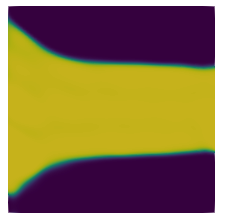}  \\
        (m) &
        \includegraphics[width=\linewidth, height=3.5cm, keepaspectratio]{ 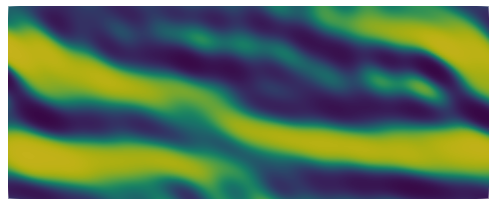} &
        \includegraphics[width=\linewidth, height=3.5cm, keepaspectratio]{ 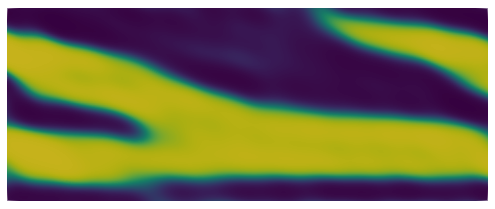} &
        \includegraphics[width=\linewidth, height=3.5cm, keepaspectratio]{ 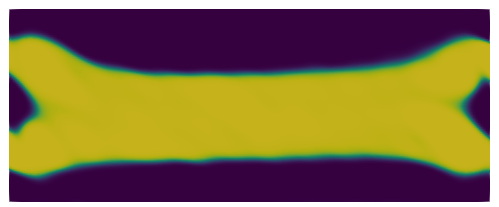} \\
    \end{tabular}
    \caption{Phase-field evolution for the 2D Stokes benchmarks (k)–(m). Rows show cases (k), (l), and (m); columns correspond to epochs 100, 200, and the final design.}
    \label{fig: 2d_stokes_design(k-m)}
\end{figure}

\begin{table}[htbp]
\centering
\caption{Hyperparameters in Algorithm \ref{alg: tri-level to} for the Stokes flow benchmarks cases (k)--(o).}
\label{tab:stokes_hparams_full}
\begin{tabular}{ccccccccccc}
\toprule
Case & $N$ & $K$ & $|X|$ & $|X_b|$ & $\sigma_{\mathrm{Fourier}}$
& $(\lambda_{\text{penal}}^{(0)},\zeta)$
& $\lambda_{\text{Dir}}$
& $\lambda_{\text{Neu}}$
& $(n^{\phi}_1,n^{\phi}_2,n^{\phi}_3)$
& $\omega^{\phi}_{\max}$ \\
\midrule
(k) & 400  & 20 & 20000 & 5000 & 20 & (100,0.98) & 10 & 10 & $(12,12,0)$ & 30 \\
(l) & 500  & 20 & 20000 & 5000 & 20 & (100,0.98) & 10 & 10 & $(12,12,0)$ & 30 \\
(m) & 500  & 20 & 30000 & 8000 & 20 & (100,0.98) & 20 & 20 & $(12,12,0)$ & 30 \\
(n) & 1000 & 10 & 30000 & 6000 & 10 & (100,1/1.02) & -- & 10 & $(6,6,6)$   & 15 \\
(o) & 700 & 10 & 30000 & 6000 & 10 & (100,1/1.02) & -- & 10 & $(6,6,6)$   & 15 \\
\bottomrule
\end{tabular}
\end{table}

For the Stokes benchmarks (k)–(o), we use the Fourier architectures for the velocity–pressure and topology networks described in Section~\ref{subsubsec: fnnsforall}. In all cases, the state network consists of a Gaussian random Fourier projection with $M_{\text{Fourier}}$ features, followed by four fully connected sine-activated hidden layers and a linear $(d+1)$-dimensional output. In 2D, each hidden layer has 64 neurons, whereas in 3D each hidden layer has 128 neurons. The topology field is represented by a single Fourier-feature phase-field network with a sigmoid head. The target volume fractions are $\beta=0.3,0.5,0.5$ for the 2D cases (k)–(m) and $\beta=0.25,0.4$ for the 3D cases (n)–(o). The hyperparameter settings are reported in Tables~\ref{tab:net_summary} and \ref{tab:stokes_hparams_full}. The Fourier feature hyperparameters are kept simple: for the 2D benchmarks we use $M_{\text{Fourier}}=512$ and projection scale $\sigma_{\text{Fourier}}=20$, whereas the 3D cases use $M_{\text{Fourier}}=256$ and $\sigma_{\text{Fourier}}=10$. The topology networks employ $\omega^{\phi}_{\max}=30$ in 2D and $\omega^{\phi}_{\max}=15$ in 3D, together with $(n^{\phi}_1,n^{\phi}_2)=(12,12)$ Fourier modes in 2D and $(n^{\phi}_1,n^{\phi}_2,n^{\phi}_3)=(6,6,6)$ in 3D. The state networks are trained with a learning rate $10^{-3}$ throughout, while the topology networks use a constant learning rate $5\times10^{-3}$ in all Stokes examples.

\begin{figure}[htbp] 
    \centering
    \includegraphics[width=0.6\linewidth]{ 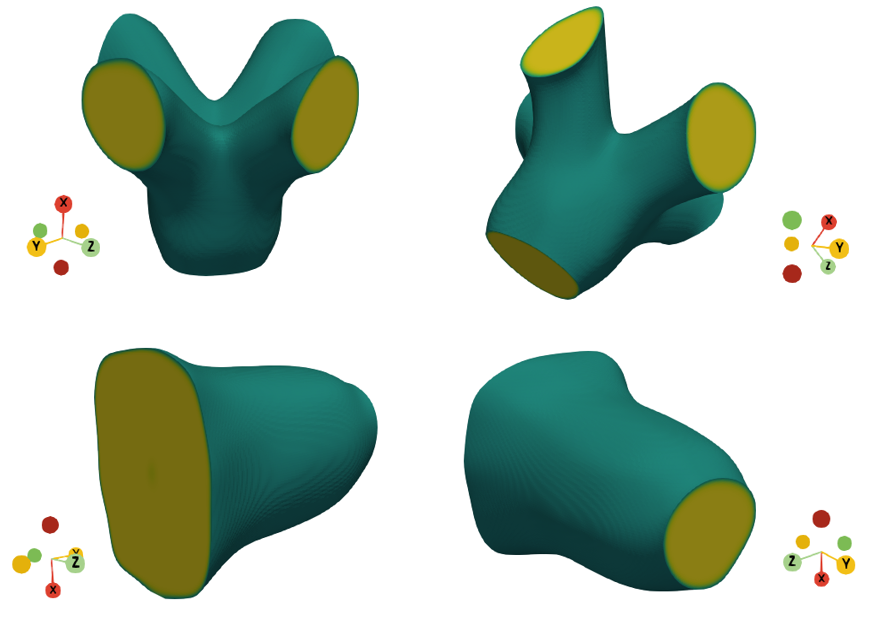}
    \caption{Final optimized topologies for the 3D Stokes benchmarks: two views of the four-inlet mixer case (n) (top) and the diffuser case (o) (bottom).}

    \label{fig: 3d_stokes_design}
\end{figure}

A Darcy penalization of the solid phase is included in all Stokes losses, with a benchmark-dependent target coefficient $\lambda_{\text{Darcy}}$, equal to $100$ in cases (l) and (o), and $50$ in the remaining cases. During the initial pretraining of the state network this coefficient is ramped linearly from $1$ to its target value over $I$ Adam steps ($I=10^4$ for (k)–(m) and $I=5\times10^3$ for (n)–(o)), and then kept fixed in the subsequent alternating optimization. The outer-loop iteration counts are $N=400$ for (k), $N=500$ for (l)–(m), and $N=1000$ and $N=700$ for the 3D cases (n) and (o), respectively. Within each outer iteration we perform $K$ state-network updates and $M$ topology updates: $K=20$ for the planar problems and $K=10$ for the 3D ones, with $M$ increased from $1$ to $10$ once $n>200$ in 2D and once $n>100$ in 3D. Interior collocation points $X$ are drawn from Latin hypercube samples in the flow domain, while boundary points $X_b$ are sampled uniformly on each inlet/outlet or wall segment, as detailed previously. The 2D cases use $|X|=2\times10^4$ (k)-(l) or $3\times10^4$ (m) interior points and $|X_b|=5\times10^3$ (k)-(l) or $8\times10^3$ (m) boundary points; the 3D tests use $|X|=3\times10^4$ and $|X_b|=6\times10^3$. The volume-penalty parameter in $\mathcal{P}_{\text{vol}}$ is increased from $\mathcal{O}(10^2)$ to $\mathcal{O}(10^3\!-\!10^4)$ using a multiplicative factor $\zeta=0.98$ in 2D and $\zeta=1/1.02$ in 3D, while Neumann and Dirichlet boundary penalties in the state loss take values between $10$ and $20$ depending on the case.

\begin{figure}[htbp] 
    \centering
    \includegraphics[width=0.34\linewidth]{ 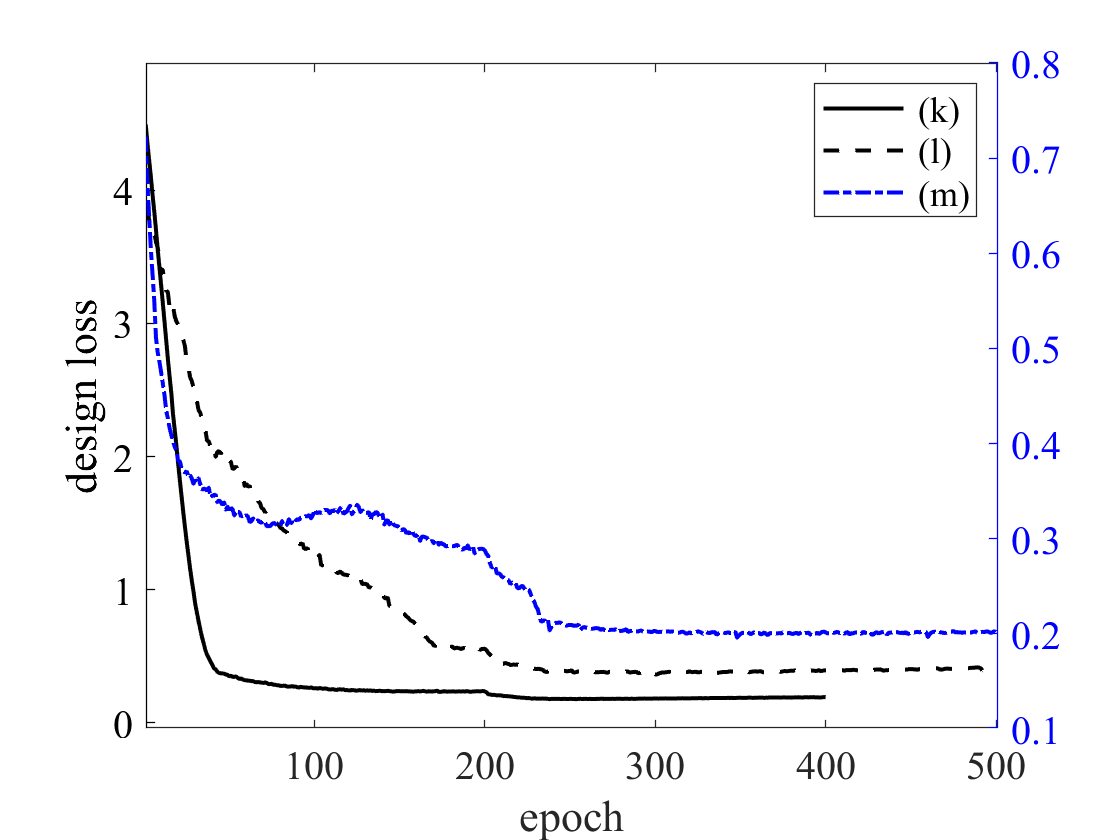}
    \includegraphics[width=0.34\linewidth]{ 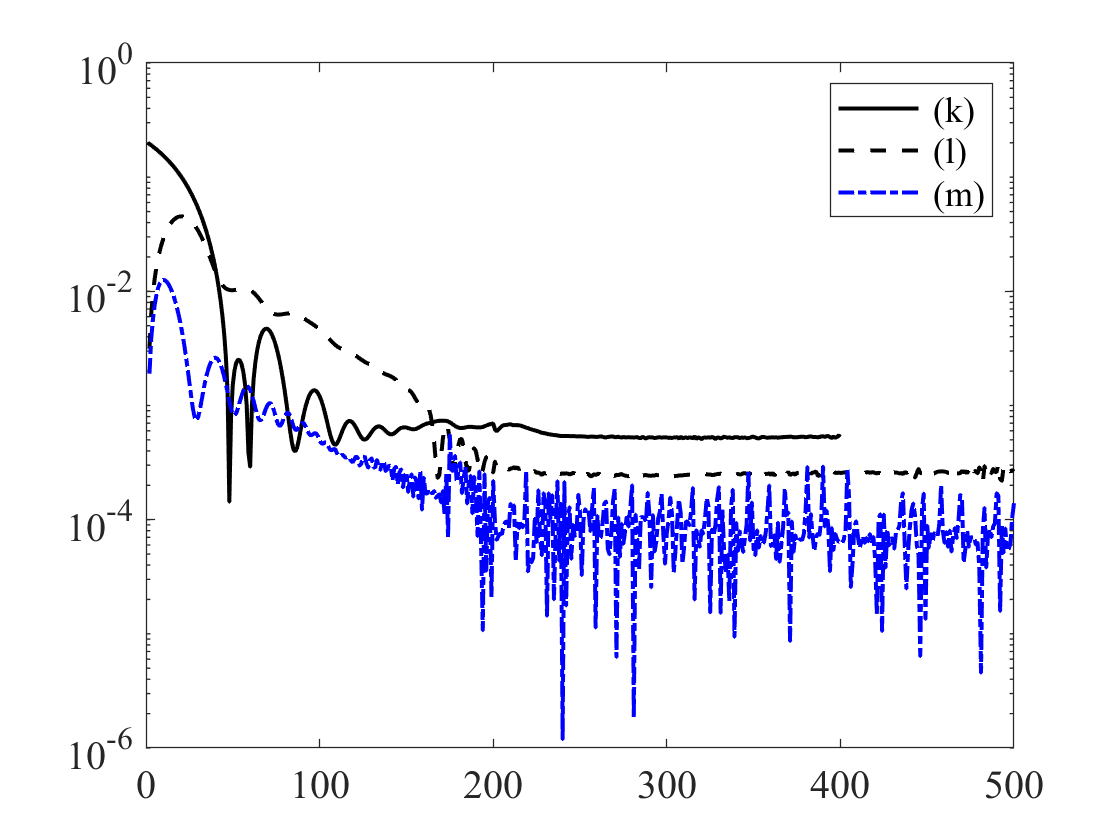}
    \\
    \includegraphics[width=0.34\linewidth]{ 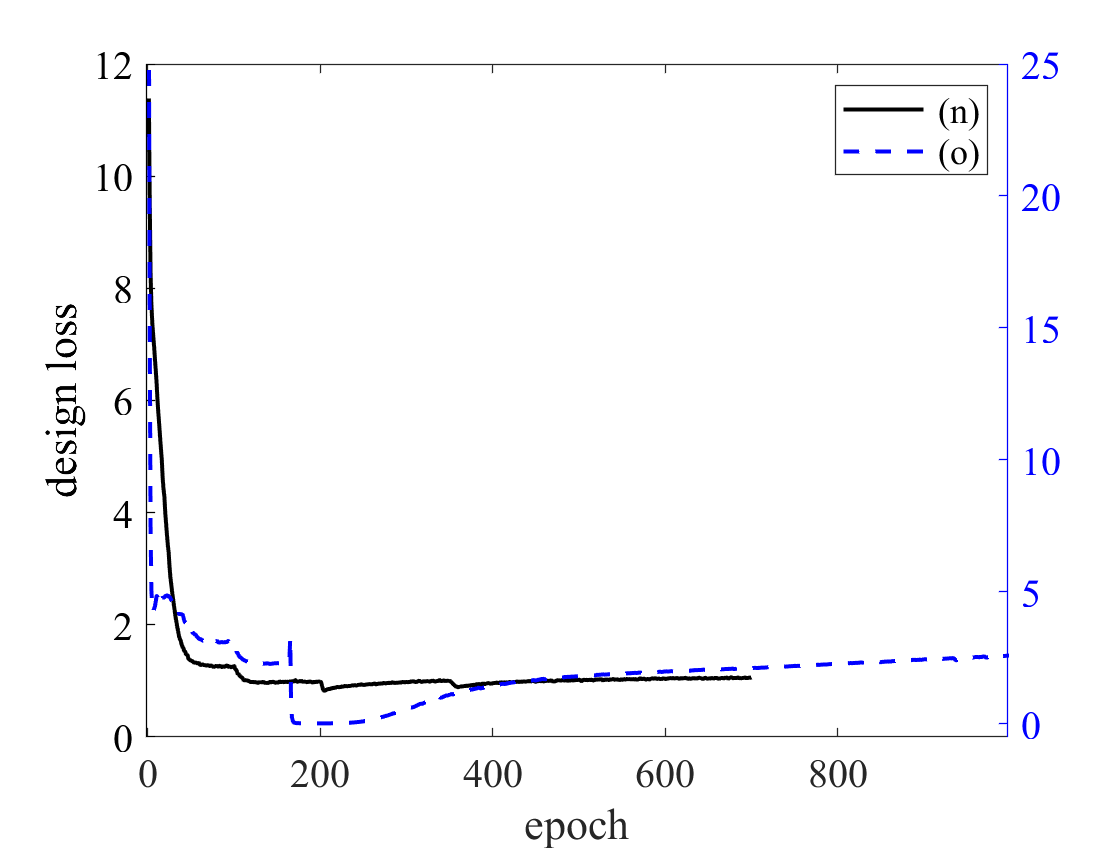}
    \includegraphics[width=0.34\linewidth]{ 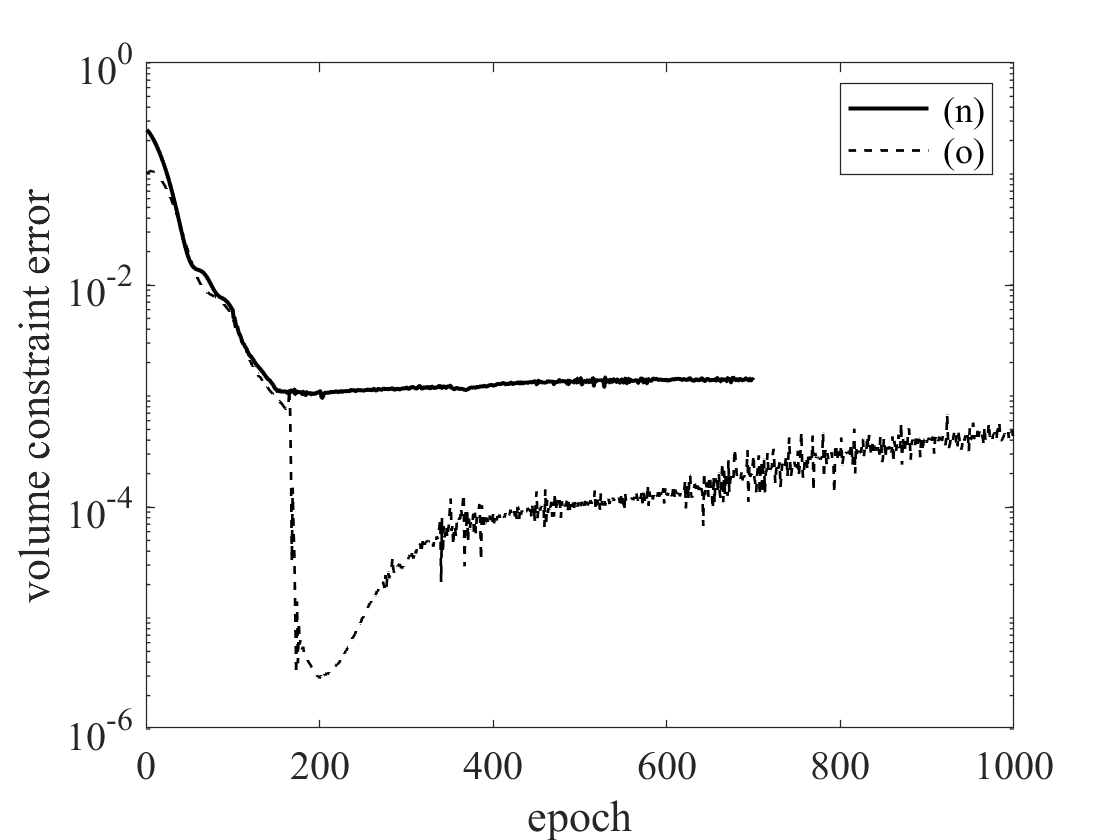}
    \caption{Convergence of design loss and volume-constraint error for the Stokes benchmarks: 2D cases (k)–(m) in the top row and 3D cases (n)- (o) in the bottom row. The loss for cases (m) and (o) are plotted against the right-hand $y$–axis.}
    \label{fig: stokes_2d3d_convergence_curve}
\end{figure}

The Stokes benchmarks further confirm the behavior observed in the elasticity tests. As shown in Figs.~\ref{fig: 2d_stokes_design(k-m)}-\ref{fig: 3d_stokes_design}, the phase field rapidly organizes from a uniform initialization into almost binary fluid/solid layouts that form smooth low–resistance channels in 2D cases (k)–(m) and tubular mixer/diffuser structures in 3D cases (n)–(o), closely aligned with the prescribed inlet–outlet configurations. The convergence histories in Fig.~\ref{fig: stokes_2d3d_convergence_curve} indicate a rapid initial decrease of the design loss followed by gradual refinement, while the volume-constraint violation is reduced by several orders of magnitude and remains small, with only mild oscillations in the more complex mixer cases. Overall, the Stokes results demonstrate that the same APF--FNN architecture and training schedule used for compliance and eigenvalue problems extend robustly to viscous flow optimization.

\subsection{Navier--Stokes flow optimization}
\label{sec:ns_numerics}

We finally assess the APF--FNN framework on non-self-adjoint problems by considering incompressible Navier--Stokes benchmarks, formulated as in Section~\ref{subsubsec:physics-NS}. Three planar tests are posed on the unit square $\Omega = (0,1)^2$ and differ only in their inflow--outflow configurations, while the fourth case is a fully three-dimensional mixer. Unless stated otherwise, all boundaries not designated as inlet or outlet are treated as no-slip walls. And we reuse the same fluid and phase-field parameters as in the Stokes benchmarks, i.e., $\rho = 0.01$ and $\varepsilon = \gamma = 0.01$.

\begin{itemize}
    \item[(p)] On the unit square, a localized horizontal parabolic inflow $\bm{u}_\theta = (1-100(x_2-0.5)^2,\,0)^{\top}$ is prescribed on the left boundary segment $\{x_1=0,\;0.7\le x_2\le 0.9\}$. The fluid leaves through the bottom segment $\{x_2=0,\;0.7\le x_1\le 0.9\}$ with vertical profile $\bm{u}_\theta = (0,\,100(x_1-0.5)^2-1)^{\top}$; the remaining edges are no-slip.

    \item[(q)] Again on $\Omega=(0,1)^2$, a parabolic inflow $\bm{u}_\theta = (1-4(x_2-0.5)^2,\,0)^{\top}$ is imposed along the entire left edge $\{x_1=0\}$. On the right boundary, a higher-speed parabolic outflow $\bm{u}_\theta = (3-108(x_2-0.5)^2,\,0)^{\top}$ is prescribed on the segment $\{x_1=1,\; \tfrac13 \le x_2 \le \tfrac23\}$, while the top and bottom edges are no-slip.

    \item[(r)] On the unit square, two vertically separated parabolic inlets are prescribed on the left boundary $\{x_1=0\}$: on the lower port $1/6 \le x_2 \le 1/3$ the inflow is $\bm{u}_\theta = (1-144(x_2-0.25)^2,\,0)^{\top}$, and on the upper port $2/3 \le x_2 \le 5/6$ it is $\bm{u}_\theta = (1-144(x_2-0.75)^2,\,0)^{\top}$. The right boundary acts as a zero-traction outflow, whereas the top and bottom edges are no-slip walls.

    \item[(s)] The 3D Navier--Stokes benchmark shares the same cubic geometry and inlet/outlet configuration as the four-inlet Stokes mixer (case~(n) in Section~\ref{subsec:stokes_examples}); only the governing equations are changed from Stokes to Navier--Stokes system.
\end{itemize}

\begin{table}[htbp]
\centering
\caption{Hyperparameters in Algorithm \ref{alg: tri-level to} for the Navier--Stokes flow benchmarks cases (p)--(s).}
\label{tab:navier_stokes_hparams_full}
\begin{tabular}{ccccccccccc}
\toprule
Case & $N$ & $K$ & $|X|$ & $|X_b|$ & $\sigma_{\mathrm{Fourier}}$
& $(\lambda_{\text{penal}}^{(0)},\zeta)$
& $\lambda_{\text{Neu}}$
& $\lambda_{\text{Dir}}$
& $(n^{\phi}_1,n^{\phi}_2,n^{\phi}_3)$
& $\omega^{\phi}_{\max}$ \\
\midrule
(p) & 500  & 40 & 20000 & 5000 & 20 & $(10,1.02)$ & 500 & 500 & $(12,12,0)$ & 30 \\
(q) & 500  & 40 & 20000 & 5000 & 20 & $(10,1.02)$ & 500 & 500 & $(12,12,0)$ & 30 \\
(r) & 500  & 40 & 20000 & 5000 & 20 & $(10,1.02)$ & 500 & 500 & $(12,12,0)$ & 30 \\
(s) & 1000 & 10 & 30000 & 6000 & 10 & $(10,1.02)$ & --  & 10  & $(6,6,6)$   & 15 \\
\bottomrule
\end{tabular}
\end{table}

\begin{figure}[htb] 
    \centering
    \settowidth{\rowlabelwidth}{ (a))} 

    \begin{tabular}{
        >{}l 
        @{\hspace{1em}} 
        >{\centering\arraybackslash}m{0.22\textwidth}
        @{\hspace{0.1em}}
        >{\centering\arraybackslash}m{0.22\textwidth}
        @{\hspace{0.1em}}
        >{\centering\arraybackslash}m{0.22\textwidth}
    }
        & \multicolumn{1}{c}{Epoch 100}
        & \multicolumn{1}{c}{Epoch 300} 
        & \multicolumn{1}{c}{Final design} \\
        \addlinespace[5pt]
        
        % 第一行图片
        (p) &
        \includegraphics[width=\linewidth, height=3.45cm, keepaspectratio]{ 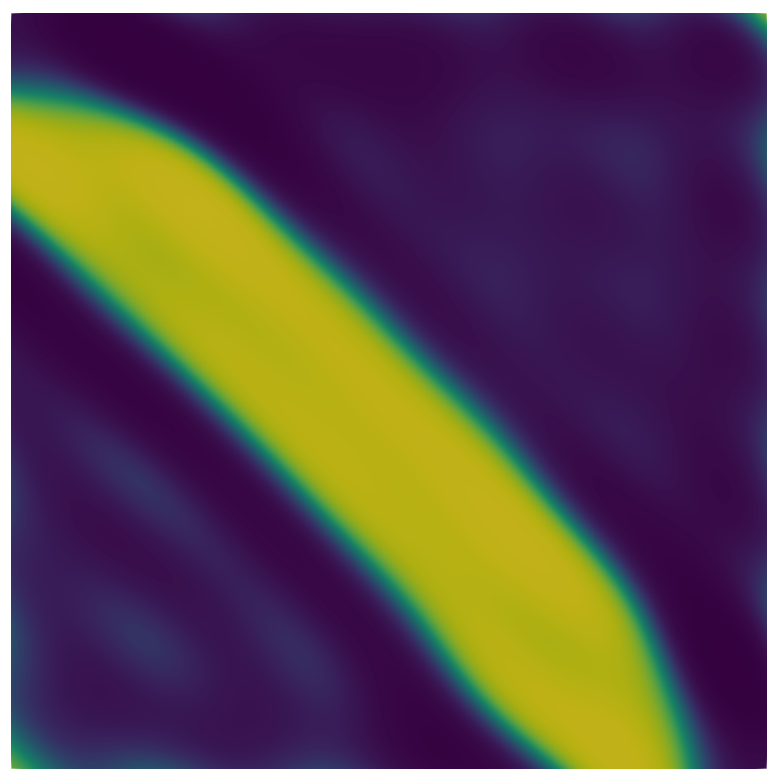} &
        \includegraphics[width=\linewidth, height=3.45cm, keepaspectratio]{ 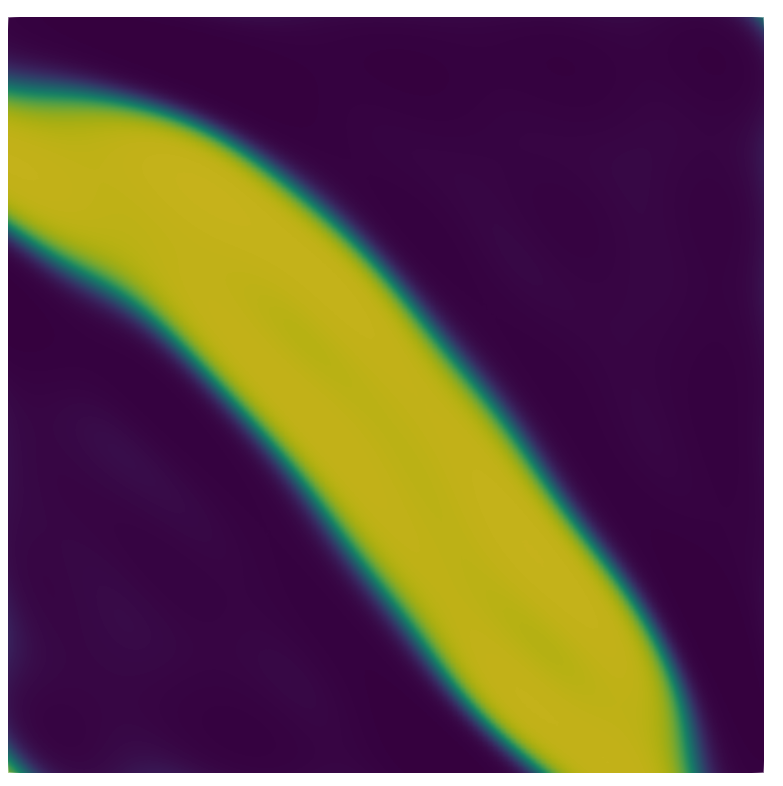} &
        \includegraphics[width=\linewidth, height=3.45cm, keepaspectratio]{ 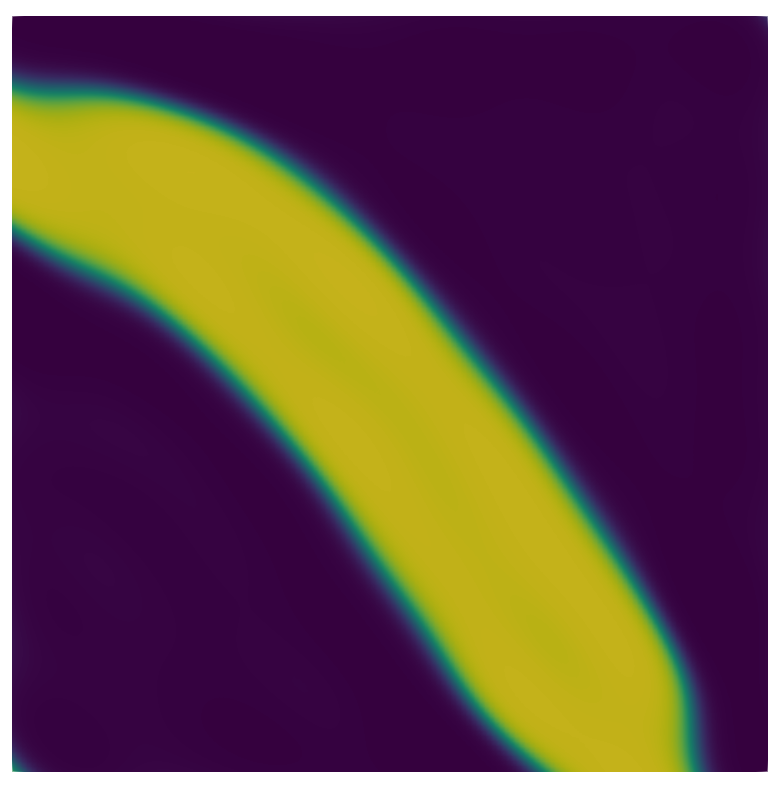} \\
        % \addlinespace[10pt]

        % 第二行图片
        (q) &
        \includegraphics[width=\linewidth, height=3.5cm, keepaspectratio]{ 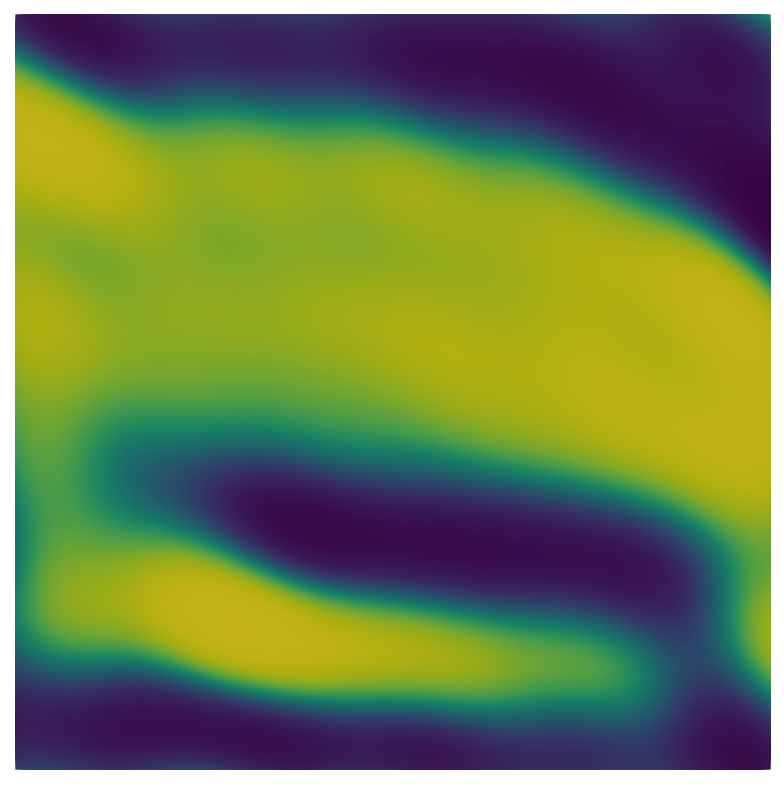} &
        \includegraphics[width=\linewidth, height=3.5cm, keepaspectratio]{ 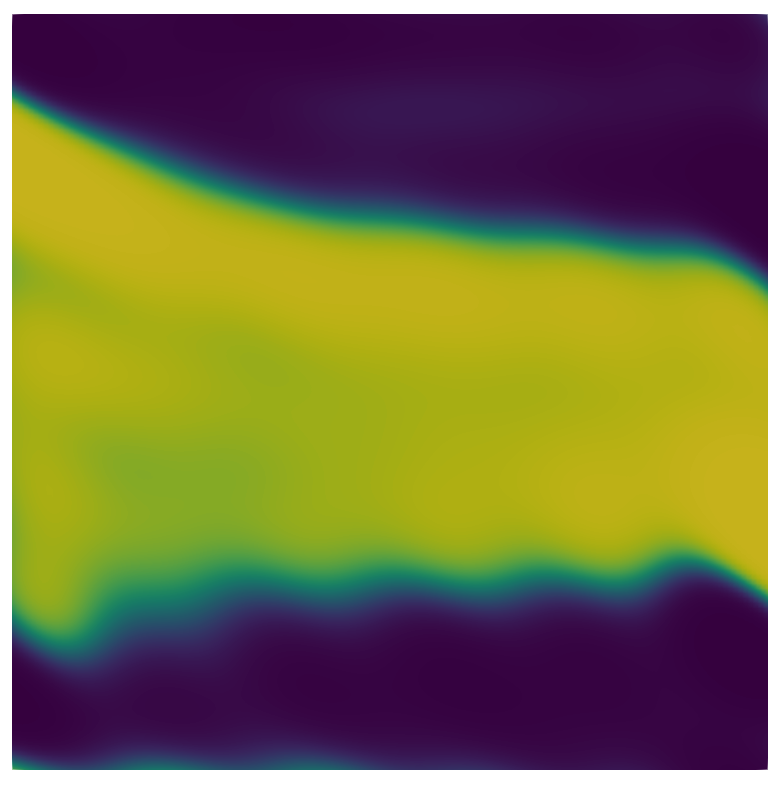} &
        \includegraphics[width=\linewidth, height=3.5cm, keepaspectratio]{ 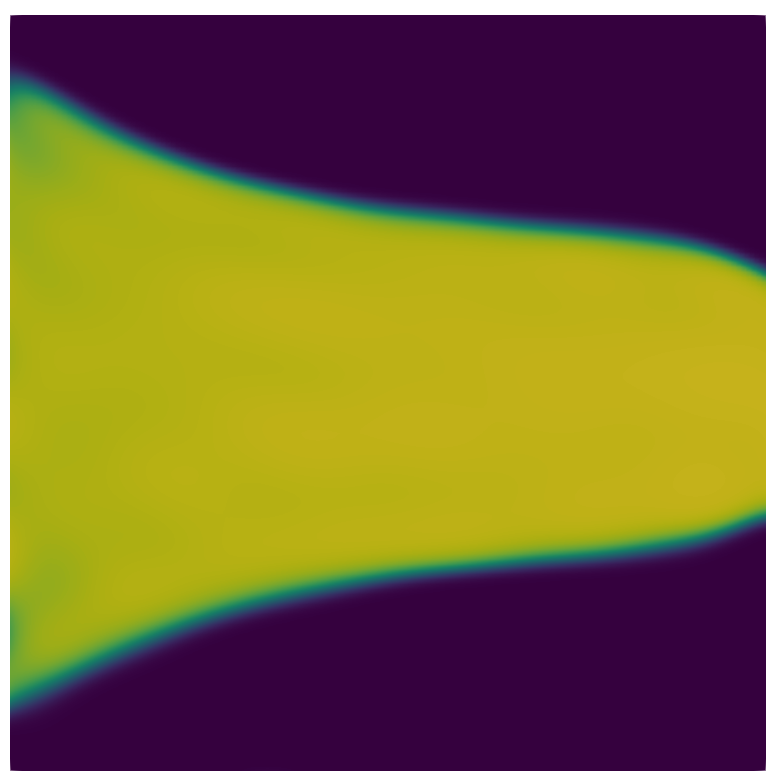} \\

        % 第三行图片
        (r) &
        \includegraphics[width=\linewidth, height=3.5cm, keepaspectratio]{ 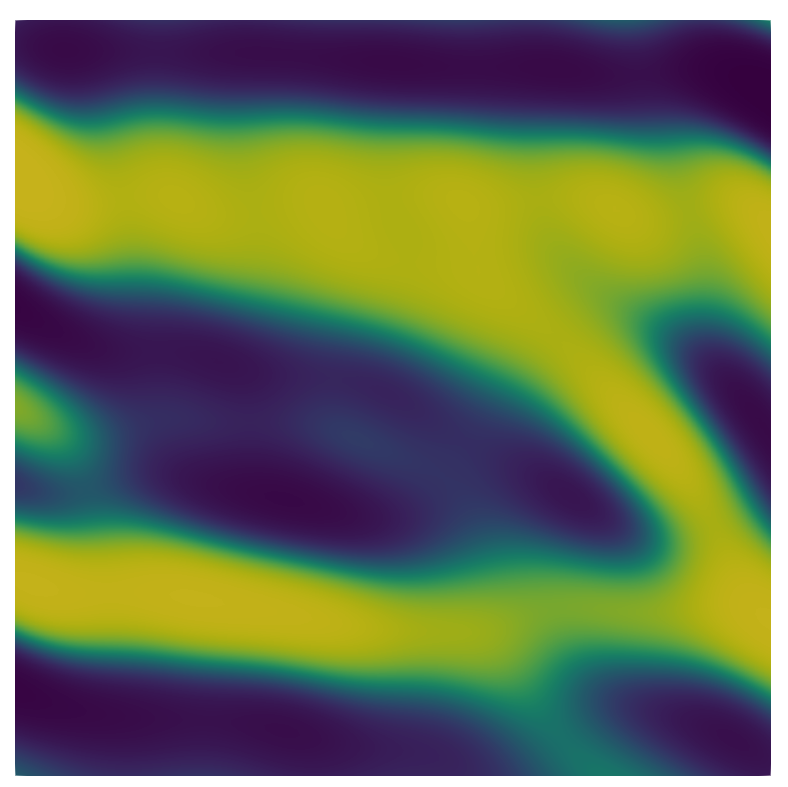} &
        \includegraphics[width=\linewidth, height=3.5cm, keepaspectratio]{ 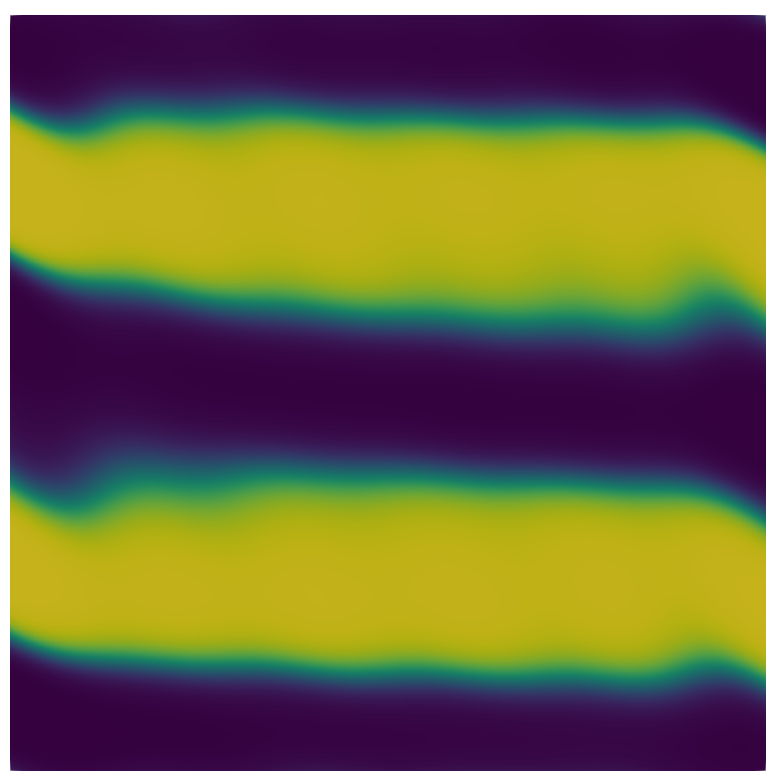} &
        \includegraphics[width=\linewidth, height=3.5cm, keepaspectratio]{ 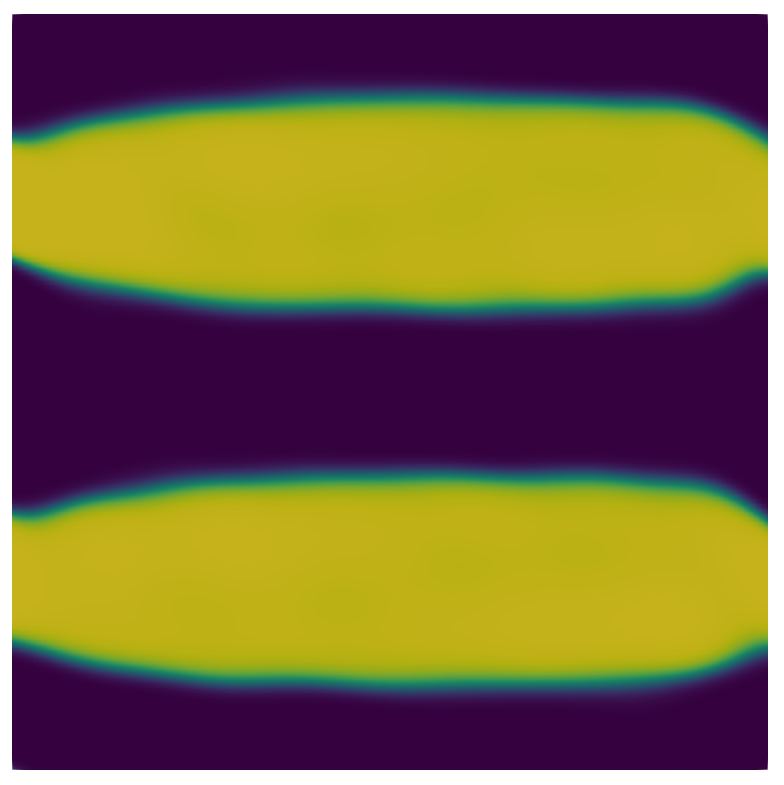}
    \end{tabular}
    
    \caption{Phase-field evolution for the 2D Navier--Stokes cases (p)–(r). Columns correspond to epochs 100, 300, and the final design.}

    \label{fig: 2d_ns_design(q-s)}
\end{figure}

For the Navier--Stokes benchmarks (p)–(s), we employ the same Fourier architectures with Stokes flow optimization, now with both state and adjoint fields represented by networks of identical structure. Each of the four tests uses a Darcy penalization of the solid phase with target coefficient $\lambda_{\text{Darcy}}=50$; during the initial pretraining of the state network this coefficient is ramped linearly from $1$ to $50$ and then kept fixed in the alternating optimization. The target volume fractions are $\beta=0.3,0.5,0.5$ for the 2D cases (p)–(r) and $\beta=0.25$ for the 3D mixer (s). The detailed hyperparameter settings are reported in Tables~\ref{tab:net_summary} and \ref{tab:navier_stokes_hparams_full}. The state and adjoint networks for (p)–(r) share the same hyperparameters: five hidden layers with 64 sine-activated neurons, a Gaussian random Fourier projection with $M_{\text{Fourier}}=512$ features and scale $\sigma_{\text{Fourier}}=20$, and a constant learning rate of $10^{-3}$. In the 3D case (s), we use a slightly wider configuration with 128 neurons per layer and $M_{\text{Fourier}}=256$, $\sigma_{\text{Fourier}}=10$, reflecting the higher complexity of the mixer flow. The topology network is the same single Fourier-layer phase-field model as in the Stokes tests, with frequency $\omega^{\phi}_{\text{max}}=30$ in 2D and $\omega^{\phi}_{\text{max}}=15$ in 3D, $(n^{\phi}_1,n^{\phi}_2)=(12,12)$ for (p)–(r) and $(n^{\phi}_1,n^{\phi}_2,n^{\phi}_3)=(6,6,6)$ for (s), and a learning rate $5\times10^{-3}$ in all Navier--Stokes examples.

\begin{figure}[htbp] 
    \centering
    \includegraphics[width=0.8\linewidth]{ 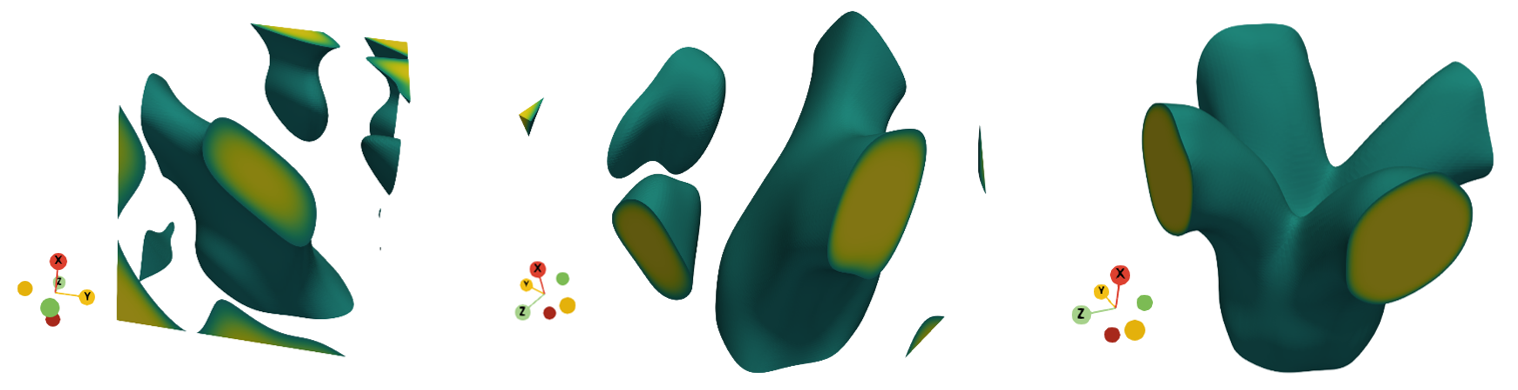}
    \caption{Phase-field evolution for the 3D Navier--Stokes mixer case (s) at epochs 100, 200, and 1000 (from left to right).}

    \label{fig: 3d_ns_design}
\end{figure}

The training schedule follows the Stokes setup with slightly deeper networks. All cases employ $I=10^4$ Adam steps of state-network pretraining, followed by $N=500$ outer iterations for cases (p)–(r) and $N=1000$ for case (s). Within each outer iteration, we perform $K=40$ state (and adjoint) updates in 2D and $K=10$ in 3D. The topology network is updated less frequently at the beginning and more aggressively as the design stabilizes: for cases (p)–(q) we increase the number of topology steps from $M=1$ to $M=40$ once $n>300$, for case (r) from $M=1$ to $M=20$ after $n>300$, and for case (s) from $M=1$ to $M=10$ after $n>100$. Interior collocation points are generated by Latin hypercube sampling, with $|X|=2\times10^4$ and $|X_b|=5\times10^3$ in the 2D problems and $|X|=3\times10^4$, $|X_b|=6\times10^3$ in the 3D mixer. The volume-penalty parameter in $\mathcal{P}_{\text{vol}}$ is gradually increased from
$\mathcal{O}(10)$ to $\mathcal{O}(10^3\!-\!10^4)$ by multiplying it by a fixed factor $\zeta = 1/1.02$ at each outer iteration. Neumann and Dirichlet penalty coefficients in the state and adjoint losses are generally set to $5\times 10^2$ to enforce the inlet/outlet and no-slip conditions, except for the Dirichlet penalty in the state network of case~(s), which is reduced to $10$.

\begin{figure}[htbp] 
    \centering
    \includegraphics[width=0.24\linewidth]{ 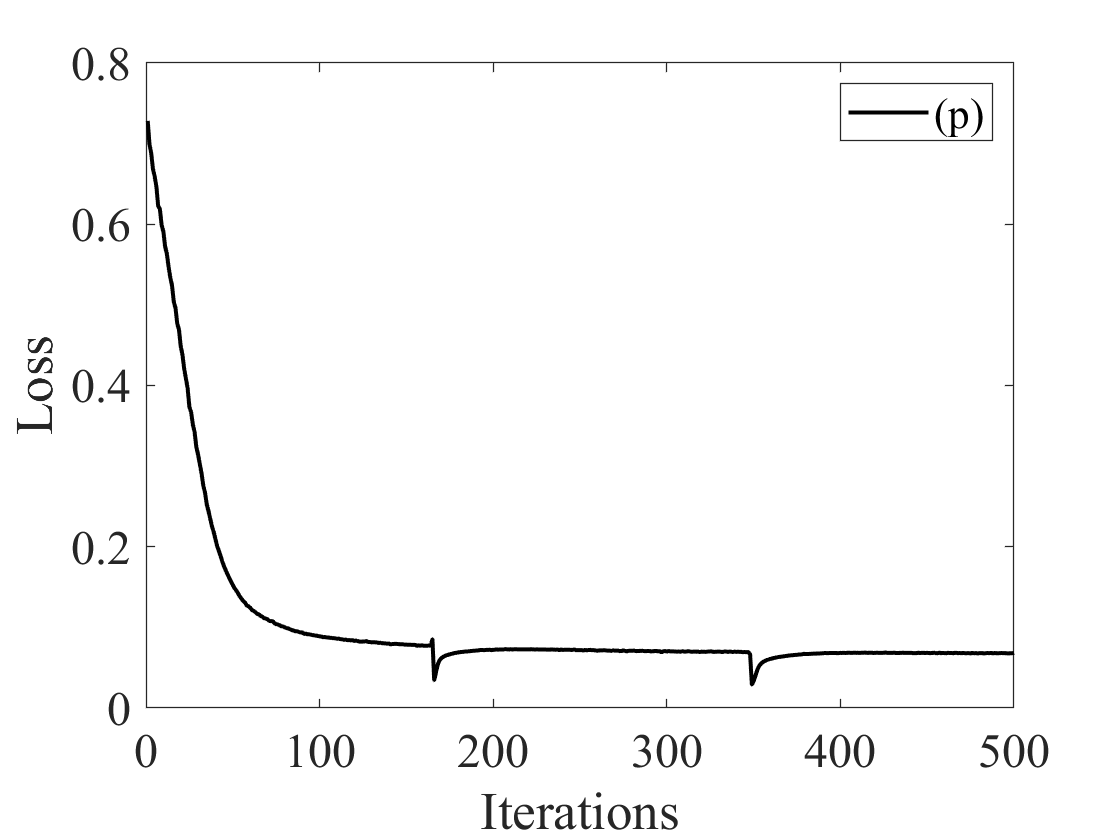}
    \includegraphics[width=0.24\linewidth]{ 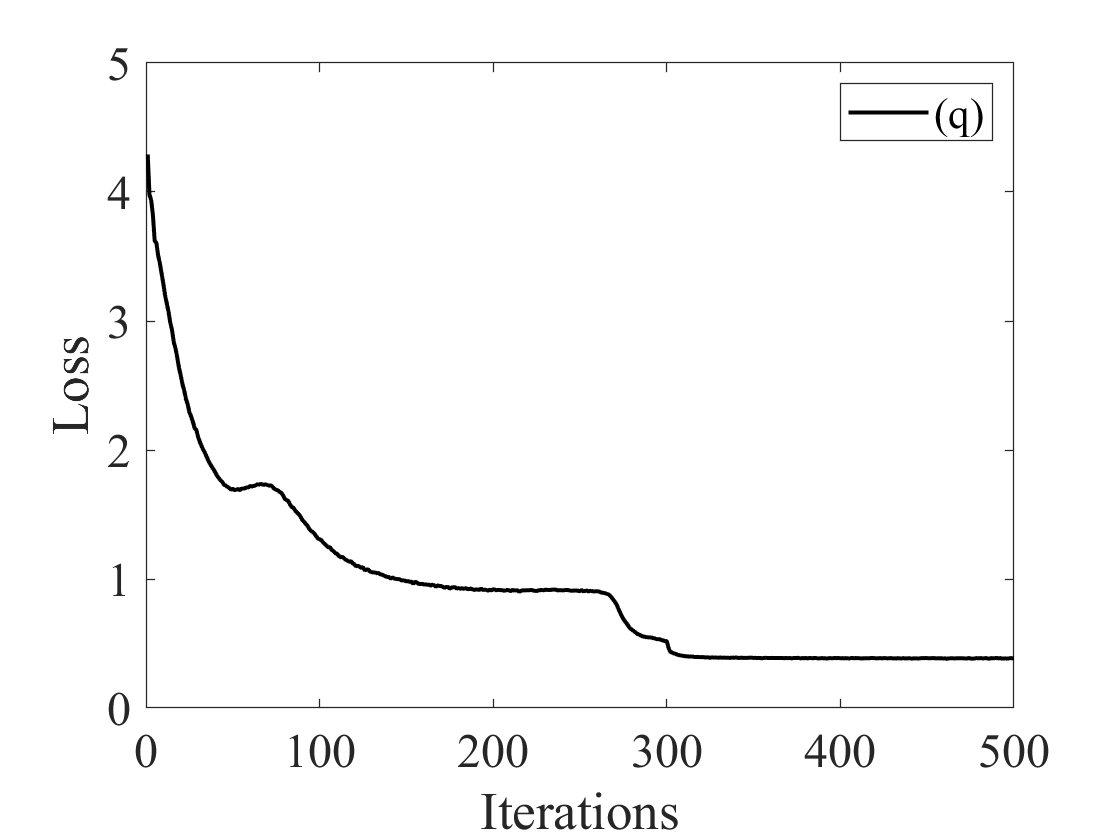}
    \includegraphics[width=0.24\linewidth]{ 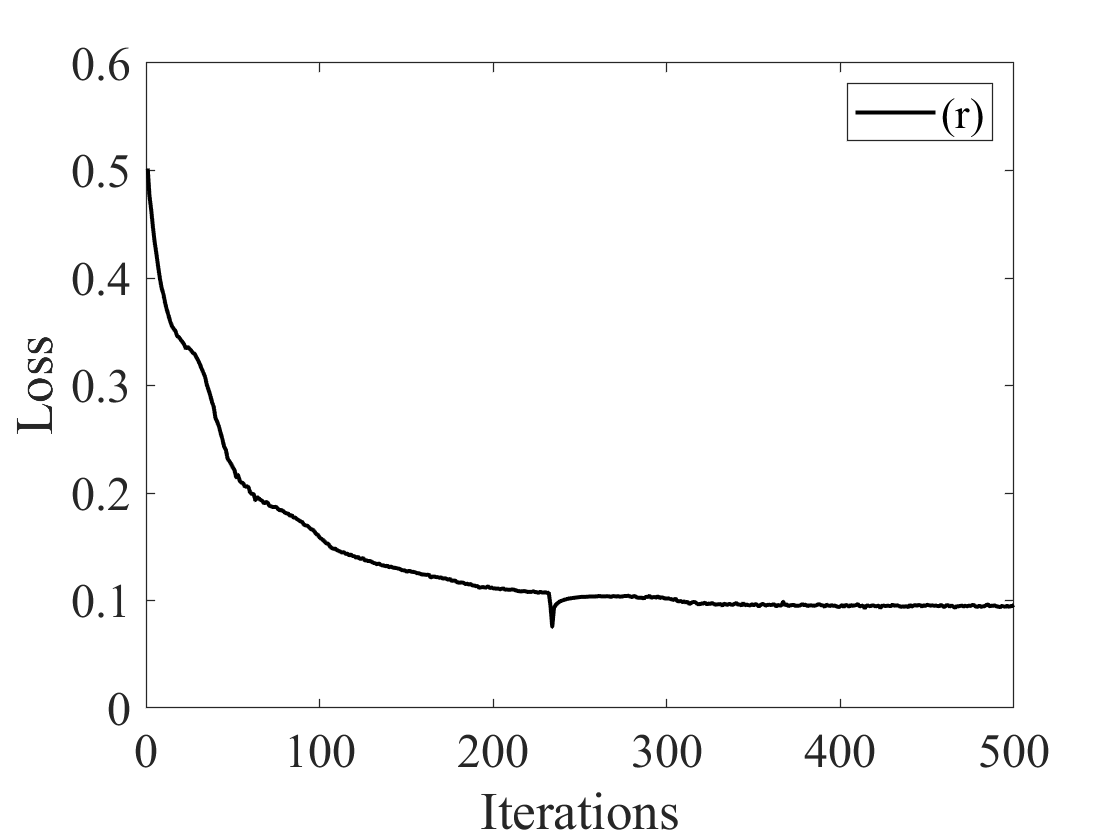}
    \includegraphics[width=0.24\linewidth]{ 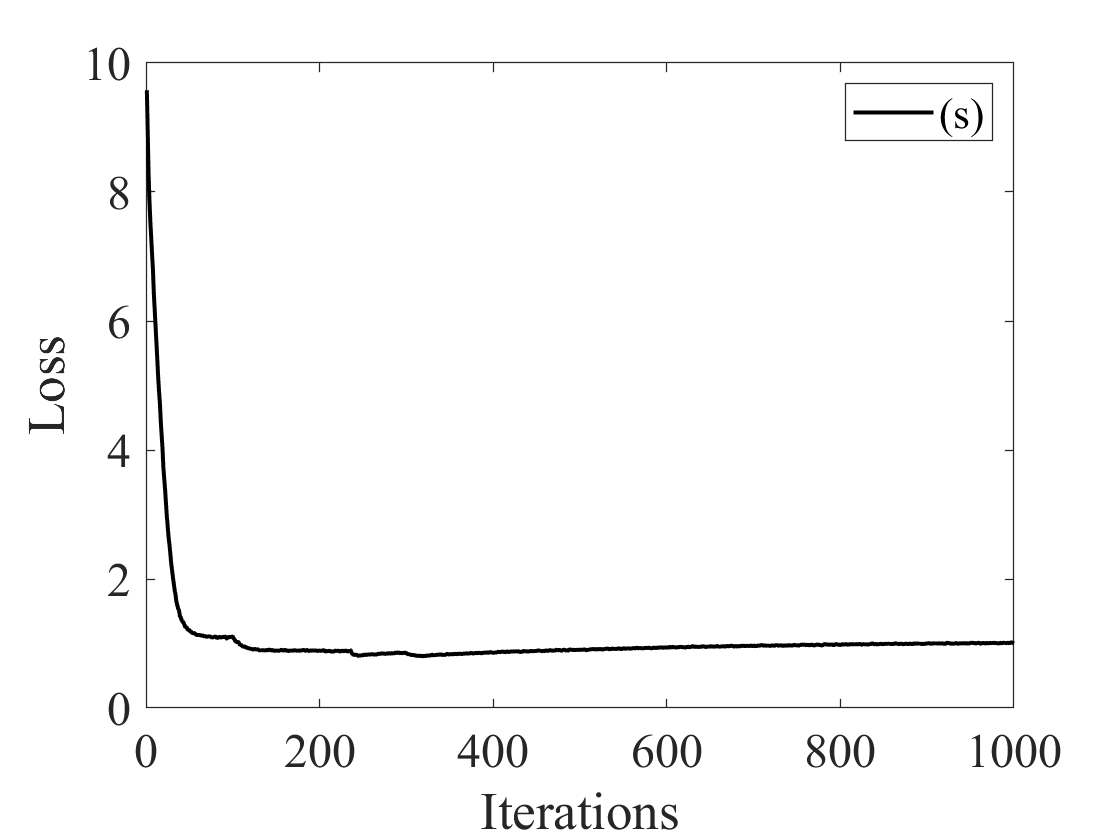}
    \\
    \includegraphics[width=0.24\linewidth]{ 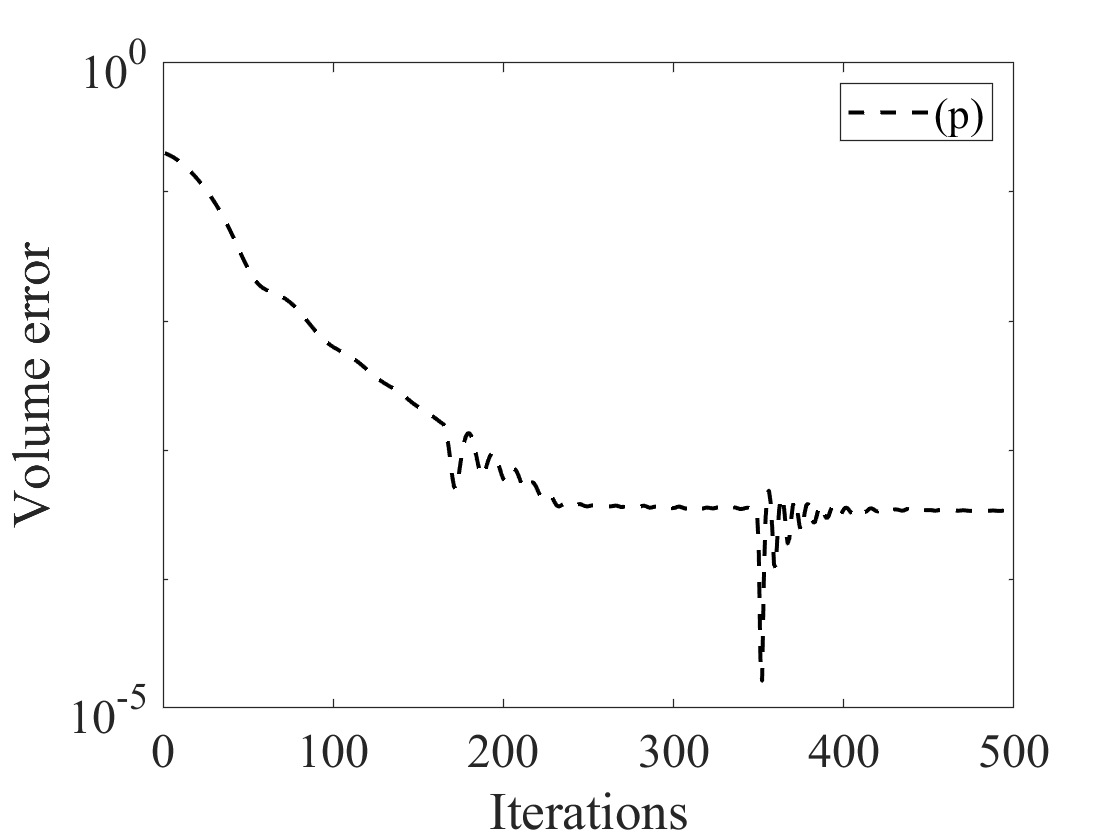}
    \includegraphics[width=0.24\linewidth]{ 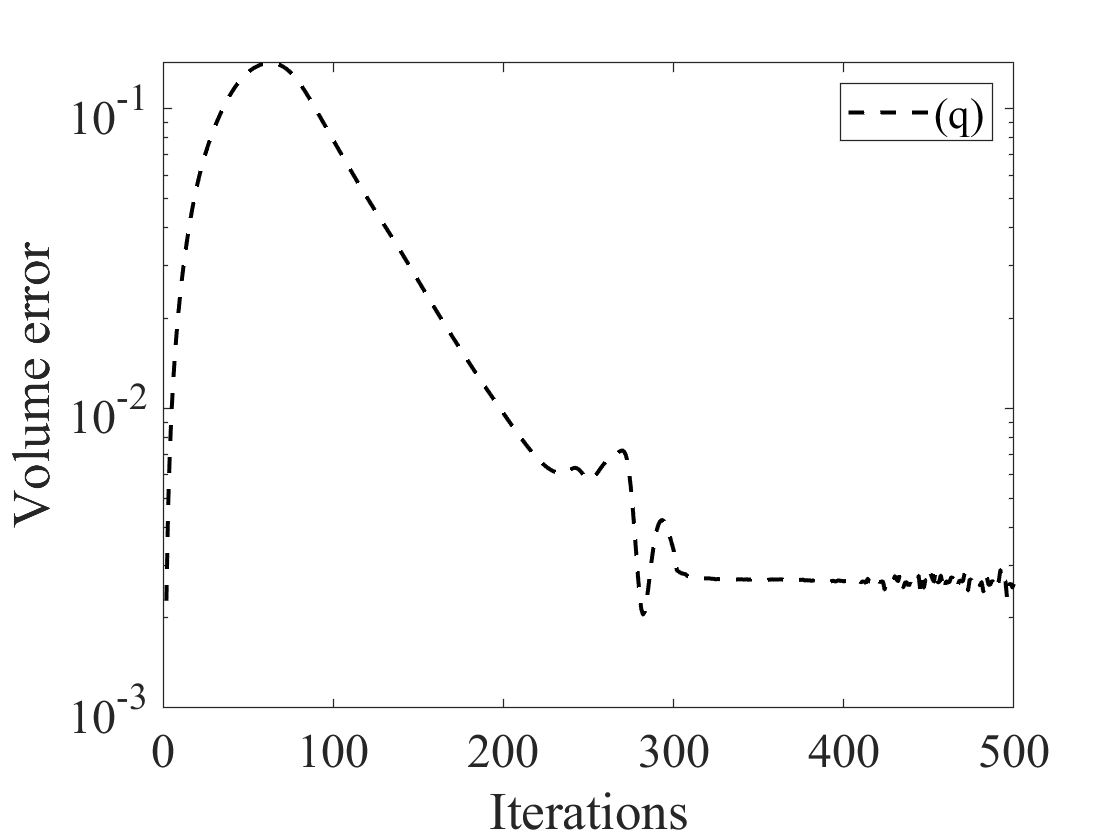}
    \includegraphics[width=0.24\linewidth]{ 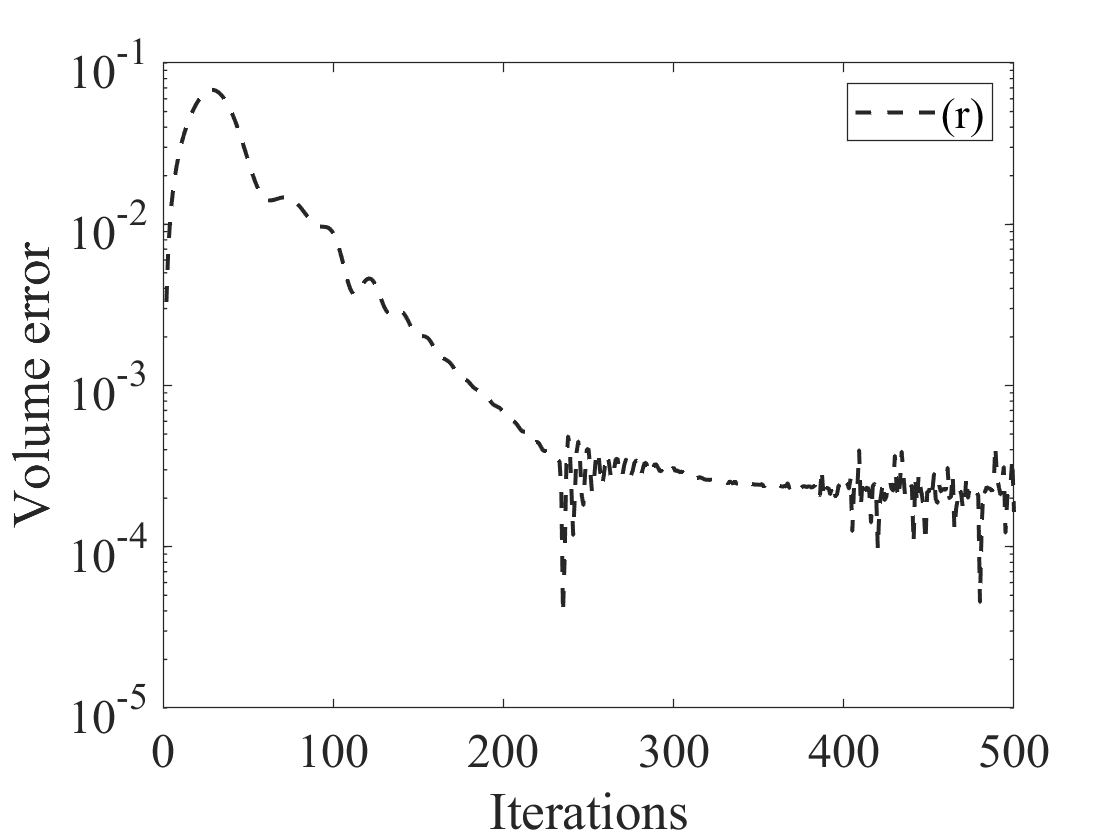}
    \includegraphics[width=0.24\linewidth]{ 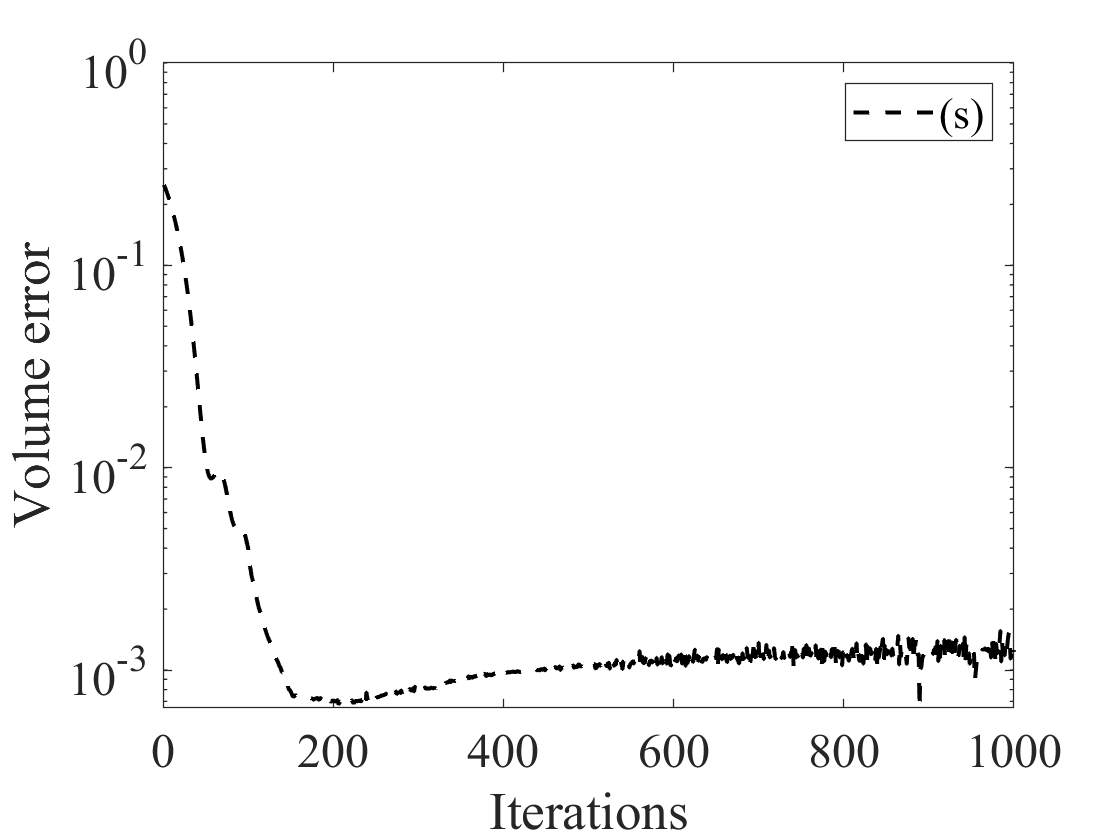}
    \caption{Convergence histories of design loss (top) and volume-constraint error (bottom) for the Navier--Stokes cases (p)–(s).}

    \label{fig: 3d_ns_curves}
\end{figure}

The Navier--Stokes benchmarks probe the fully non-self-adjoint regime, where the convective term couples the state and adjoint networks most strongly. As shown in Figs.~\ref{fig: 2d_ns_design(q-s)}–\ref{fig: 3d_ns_design}, the phase field still evolves from a uniform initialization to clean, almost binary fluid/solid layouts: in 2D cases (p)–(r) the designs form smooth channels and splitter structures that redirect the inflow jets, while the 3D mixer case (s) develops a coherent branching manifold that merges the four inlets into a single outlet jet. The convergence plots in Fig.~\ref{fig: 3d_ns_curves} confirm that, despite slightly stronger oscillations than in the Stokes case, the design loss decreases steadily and the volume constraint remains tightly enforced. Occasional spikes in the loss curves are mainly induced by transient corrections of the volume constraint, when the phase field is re-adjusted to satisfy the target volume more strictly. Overall, these results demonstrate that APF--FNNs can handle both self-adjoint and non-self-adjoint phase-field topology optimization problems within a single unified framework.

\section{Conclusion}
\label{sec:conclusion}

The APF--FNN framework proposed in this work provides a unified and physics-driven approach to phase-field topology optimization. By parameterizing the state, adjoint and topology fields with three dedicated Fourier neural networks and training them in a stable alternating scheme, the method decouples the tri-level optimization problem while retaining a tight coupling to the underlying physics. This architecture offers a flexible foundation for a broad range of design tasks, as demonstrated on four representative benchmark classes: compliance minimization, eigenvalue maximization, and Stokes/Navier--Stokes flow optimization in both 2D and 3D.

A key ingredient of the framework is the topology network update strategy, which combines adjoint-based sensitivity analysis with automatic differentiation to steer the design along efficient descent directions. The inclusion of the Ginzburg--Landau energy functional in the topology loss acts as a physical regularizer, promoting sharp, smooth and well-defined two-phase structures. In contrast to classical phase-field methods that evolve gradient-flow equations in pseudo-time, the APF--FNN directly targets steady-state phase fields satisfying first-order optimality conditions. The numerical experiments indicate that the proposed framework consistently produces high-resolution phase-field representations with competitive computational performance, while tightly enforcing the volume constraint for both self-adjoint and non-self-adjoint systems. Compared with classical FEM-based phase-field approaches, APF--FNNs achieve comparable or improved objective values with substantially better scalability in three-dimensional settings, owing to the avoidance of repeated high-fidelity PDE solves on dense meshes. These results suggest that APF--FNNs form a versatile and scalable basis for future extensions to coupled multiphysics problems and topology optimization under uncertainty.

\section*{Acknowledgments}
This work is supported in part by the National Key Research and Development Program of China under Grant (Nos. 2022YFA1004402 and 2025YFE0150300), the National Natural Science Foundation of China under Grant (No. 12471377), and the Science and Technology Commission of Shanghai Municipality (Nos. 23JC1400503 and 22DZ2229014).

\section*{Data Availability}
The data and code supporting this study are available upon reasonable request.

\bibliographystyle{unsrt}
\bibliography{references}

\end{document}